\input amstex.tex
\documentstyle{amsppt}
\magnification=\magstep1
\voffset-2pc
\refstyle{A}
\nologo
\PSAMSFonts
\loadbold
\def\epsilon{\varepsilon}
\input xy
\xyoption{curve}
\xyoption{matrix}
\xyoption{arrow}
\CompileMatrices

\def\line{\hbox to\hsize}

\def\bfrecc{\ar@{=}|*{\SelectTips{cm}{}\dir{>}}[r]}
\def\cfrecc{\ar@{=}|*{\SelectTips{cm}{}\dir{<}}[r]}

\def\grigia{{\circledast}}
\def\bbianca{{\bigcirc}}
\def\ggrigia{{\ssize\circledast}}

\def\bbianca{{\ssize\bigcirc}}

\def\bianca{\!\circ\!}
\def\nera{\!\bullet\!}

\newcount\t
\newcount\q
\newcount\x
\newcount\y

\t=0 \q=0 \x=0 \y=0 

\def\o{\bold{o}}


\long\def\se#1{\advance\q by1
\t=0 \x=0 \y=0 
\csname head\endcsname{\number\q.} #1\endhead}

\long\def\form#1{\global\advance\t by1%
$$ #1 \tag\number\q.\number\t $$}


\long\def\thm#1{\advance\x by1
\csname proclaim\endcsname {Theorem \number\q.\number\x} #1\endproclaim}

\long\def\lem#1{\advance\x by1
\csname proclaim\endcsname {Lemma \number\q.\number\x} #1\endproclaim}
\long\def\cor#1{\advance\x by1
\csname proclaim\endcsname {Corollary \number\q.\number\x} #1\endproclaim}
\long\def\prop#1{\advance\x by1
\csname proclaim\endcsname {Proposition \number\q.\number\x} #1\endproclaim}

\long\def\esemp#1{\advance\y by1
\csname example\endcsname{Example \number\q.\number\y} #1\endexample}


\def\dimo{\demo{Proof}}
\def\enddimo{\qed\enddemo}

%

%

\def\alweyl{\bold{W}(\Cal{R})}
\def\alweyle#1{\bold{W}(\Cal{R},#1)}
\def\alaffe#1{\bold{A}(\Cal{R},#1)}
\def\anweyl{{\bold{W}(\frak{g},\frak{h})}}
\def\alaff{\bold{A}(\Cal{R})}

\def\Hom{\roman{Hom}}

\def\alaffine#1{\bold{A}_{#1}(\Cal{R})}
\def\alaffineq#1#2{\bold{A}_{#1}(\Cal{R},\Cal{#2})}
\def\rootim{{\Cal{R}}_{\roman{im}}}
\def\rootcp{{\Cal{R}}_{\roman{cp}}}
\def\rootre{{\Cal{R}}_{\roman{re}}}
\def\rootimc{{\Cal{R}_{\bullet}}}
\def\supp#1{{\roman{supp}}_{#1}}

\def\rootimp{{\Cal{R}_{\ast}}}

\def\Interc#1{\bold{Int}_{\Bbb{C}}(\Hat{\frak{#1}})}


\NoBlackBoxes
\topmatter
\title
Orbits of real forms \\
in complex flag manifolds
\endtitle
\author
Andrea Altomani\\
Costantino Medori\\ 
Mauro Nacinovich
\endauthor
\date 24/11/2006 \enddate
\abstract
We study, from the point of view of $CR$ geometry,
the orbits $M$ of a real form $\bold{G}$ of a complex
semisimple Lie group $\Hat{\bold{G}}$ in a complex flag manifold
$\Hat{\bold{G}}/\bold{Q}$.
In particular we
characterize those that are of finite type and satisfy some
Levi nondegeneracy conditions. These properties
are also graphically described by attaching to them
some cross-marked diagrams that generalize those for minimal orbits that
we introduced in a previous paper.
By constructing canonical fibrations 
over real flag manifolds, with simply connected
complex fibers,
we are also able to compute their fundamental group.
\endabstract
\toc
\head 1.
Introduction \endhead
\head 2.
Homogeneous $CR$ manifolds and $CR$ algebras\endhead
\head 3.
Parabolic subalgebras and complex flag manifolds\endhead
\head 4.
Parabolic $CR$ algebras and parabolic $CR$ manifolds\endhead
\head 5. 
Fit Weyl chambers and $CR$ geometry of $M(\frak{g},\frak{q})$
\endhead 
\head 6.
Canonical fibrations over a parabolic $CR$ manifold
\endhead
\head 7.
The isotropy subgroups
\endhead
\head 8.
The fundamental group of parabolic $CR$ manifolds
\endhead
\head 9.
Examples
\endhead
\endtoc 
\address
A.~Altomani:
Dipartimento di Matematica, Universit\`a di Roma ``Tor Vergata'', 
Via della Ricerca Scientifica, 00133 Roma (Italy)
\endaddress
\email altomani\@mat.uniroma2.it\endemail
\address
C.~Medori:
Dipartimento di Matematica, Universit\`a di Parma,
Parco Area delle Scienze~53/A,
43100 Parma (Italy)
\endaddress
\email medori\@prmat.math.unipr.it\endemail
\address
M.~Nacinovich:
Dipartimento di Matematica, Universit\`a di Roma ``Tor Vergata'', 
Via della Ricerca Scientifica, 00133 Roma (Italy)
\endaddress
\email nacinovi\@mat.uniroma2.it\endemail
\subjclassyear{2000}
\endtopmatter
\leftheadtext{A.Altomani\quad C.Medori\quad M.Nacinovich}

\document 
\se{Introduction}
A complex flag manifold  
is a compact complex manifold $\frak{M}$
that is homogeneous
for the action of a semisimple complex Lie group $\Hat{\bold{G}}$.
The orbits $M$ in $\frak{M}$
of a real form $\bold{G}$ of $\Hat{\bold{G}}$ inherit from the complex
structure of $\frak{M}$ a $\bold{G}$-homogeneous $CR$ structure.
In this way we obtain a large class 
of $CR$ manifolds, that we call {\it parabolic $CR$ manifolds}.
They are homogeneous
for the $CR$ action of a real semisimple Lie group. Special examples are
the compact 
standard homogeneous $CR$ manifolds, corresponding to Levi-Tanaka
algebras (see e.g. [MN1], [MN2], [MN3]),
and the symmetric $CR$ manifolds in [LN]. \par
In [AMN] we investigated
the $CR$ structure of all minimal orbits. Here we extend our consideration
to {\it all} orbits $M$ of $\bold{G}$ in $\frak{M}$. 
Several
$CR$ and topological invariants of $M$ are encoded in 
its {\it cross-marked diagrams}.
These are obtained from the Dynkin diagrams of the root system 
$\Cal{R}$ of
the complex Lie algebra $\Hat{\frak{g}}$ of $\Hat{\bold{G}}$
with respect to the complexification of a Cartan subalgebra $\frak{h}$
of the Lie algebra $\frak{g}$ of $\bold{G}$. 
At variance with the systematic use of
cross-marked Satake diagrams in the minimal orbit case,
here we need more general descriptions, involving 
general
Cartan subalgebras $\frak{h}$ of the Lie
algebra $\frak{g}$ of $\bold{G}$ and 
different choices of simple systems
of roots in $\Cal{R}$.\par
The orbits of $\bold{G}$ in $\frak{M}$
were already considered in [W]. 
Among the more recent contributions to the study of this subject, we cite
[Ka] in the context of infinite dimensional representation theory,
[GM], [HW] and [KZ]
for applications to the geometry of symmetric spaces.
\par
In our work we stress 
the point of view of $CR$ geometry. Beyond 
its independent interest, it
also provides a
helpful perspective to investigate their topology. 
The open orbits have been extensively studied (see e.g. [FHW]). 
In particular, they are all
simply connected (see [W]). Also
the topology of the real flag manifolds has been 
thoroughly investigated (see e.g. [CS], [DKV], [Wi]). 
In our paper we show that every parabolic $CR$ manifold $M$ is
the total space of a canonical fibration over a real flag manifold.
The fiber may be disconnected and each connected component is
a simply connected complex manifold, which
can be retracted onto an open orbit.
This essentially 
reduces the computation of 
the fundamental group of $M$ to counting the 
connected components of
the fiber.
\par\smallskip
The paper is organized as follows. 
\par
In \S{2} we review 
the notion of $CR$ algebra from [MN4], [AMN], that was also recently
utilized in [Fe] and [FeK]. We collect here the main results and fix the
notation that will be employed in the rest of the paper, also correcting
a mistake, related to the notion of {\it ideal nondegeneracy}
in [MN4], that
was pointed out to us by G.Fels. \par
In \S{3} we
quickly rehearse parabolic complex Lie subalgebras and complex flag
manifolds. \par
In \S{4} we begin the study of the $CR$ algebras that are
associated to the real orbits $M$ in the complex flag manifolds $\frak{M}$,
also investigating the canonical $\bold{G}$-equivariant maps of [MN4]
in this special situation.
\par
\S{5} is the core of our investigation of the $CR$ properties of $M$.
Through the introduction of adapted Cartan subalgebras and S- and V-fit
Weyl chambers, we associate to $M$ special systems of simple roots.
Weak (i.e. holomorphic
according to [BER])
nondegeneracy and fundamentality (i.e. finite type according to
[BlG]) are proved to be equivalent to properties of these systems of
simple roots. These can be checked from the pattern of
some {\it cross marked diagrams} associated to $M$, 
that generalize those of [AMN], [MN2], [LN].
\par
In \S{6}  we turn to the construction of homogeneous $CR$
manifolds that fiber over our orbit $M$ and that are useful both for
finding the S- and V-fit Weyl chambers and for investigating the topological
properties of $M$ in the following sections. In particular,
we construct the {\it weakest $CR$ model} of $M$, that is a step to build
a chain of fibrations, with simply connected fibers, 
that in some instances coincides with, and 
in general can be considered as
a substitute of, the {\it holomorphic arc
components} of [W].
\par
In \S{7} we investigate the connectivity of the isotropy subgroup of $M$.
This is needed to study the 
connectivity of the fibers of a fibration 
of $M$ over a real flag manifold $M'$,
that we utilize to compute
the fundamental group of $M$.
This is a somehow delicate point: the simply connected fibers
of our construction may be not connected (the
statement of [W, Theorem 8.15(1b)],
that would imply connectedness, is not correct).
We use classical results from [BT1], [BT2] to characterize Cartan subgroups
and isotropy subgroups of connected semisimple real linear
groups in terms of characters.\par
Finally, in \S{8} we discuss the fundamental group of $M$.\par
In the last section (\S{9}) of our
paper we provide 
several examples which show how effective our results are
for the study of $CR$ and topological properties of
the orbits.\nopagebreak


\q=1
\se{Homogeneous $CR$ manifolds and $CR$ algebras}
A $CR$ manifold of type $(n,k)$
is a triple $(M,HM,J)$, consisting of a smooth
paracompact manifold $M$ of real dimension $(2n+k)$,  
a smooth real vector subbundle $HM$ of rank $2n$ of its
real tangent bundle $TM$, and  a smooth complex structure
$J\colon HM@>>>HM$ on the fibers of $HM$.
It is also required that $J$ satisfies
the {\it formal integrability conditions}\,:
 $\left[\Cal{C}^{\infty}(M,T^{0,1}M),\Cal{C}^{\infty}(M,T^{0,1}M)\right]
\subset\Cal{C}^{\infty}(M,T^{0,1}M).$
Here
 $T^{0,1}M=\left\{X+iJX\; \right|\; \left. X\in HM\right\}$ is the
complex subbundle of the complexification $\Bbb C HM$ of $HM$
corresponding to the eigenvalue $-i$ of $J$; we have 
 $T^{1,0}M\cap T^{0,1}M= 0$ and $T^{1,0}M\oplus T^{0,1}M=\Bbb C HM$,
where $T^{1,0}M=\overline{T^{0,1}M}$. When $k=0$, we recover the
abstract definition of a {\it complex} manifold, via the 
Newlander-Nirenberg theorem.\par\smallskip
If $(M,HM,J)$ and $(M',HM',J')$ are $CR$ manifolds, a smooth 
$f\colon M@>>>M'$ is a {\it $CR$ map} if\,: $(i)$ $df(HM)\subset HM'$, 
$(ii)$
$df\circ J=J'\circ df$ on $HM$. Assume that $(M',HM',J')$ is
a $CR$ manifold and $f\colon M@>>>M'$ a smooth immersion. 
For $x\in M$, we set\,:
$$H_xM=[df(x)]^{-1}\left([df(T_xM)\cap H_{f(x)}M']\cap
[J'(df(T_xM)\cap H_{f(x)}M')]\right)\quad
\forall x\in M$$
and 
$$J_x(v)=[df(x)]^{-1}\left(J'([df(x)](v))\right)\quad
\forall v\in H_xM,\;\forall x\in M\,.$$
If the dimension $H_xM$ is a constant integer, independent of 
$x\in M$, then the disjoint union $HM$ of the $H_xM$'s, and the
map $J:HM@>>>HM$, equal to $J_x$ on the fiber $H_xM$, define
a $CR$ manifold $(M,HM,J)$. This is the $CR$ structure on $M$
with the maximal $CR$ dimension among those for which
$f$ is a $CR$ map. 
In this case the map $f:M@>>>M'$ is called a $CR$ {\it immersion}.
If $(M',HM',J')$ is of type $(n',k')$ and      
$(M,HM,J)$ of type $(n,k)$, we always have $n+k\leq n'+k'$.
The immersion is {\it generic} when the equality 
$n+k=n'+k'$ holds. When $f:M@>>>M'$ is also an embedding, we say that
$f$ is
a $CR$ {\it embedding} or a {\it generic $CR$ embedding}, respectively.\par
A $CR$ map $f\colon (M,HM,J)@>>>(M',HM',J')$ is a {\it $CR$ submersion}
if $f:M@>>>M'$ is a smooth submersion and moreover
$df(x)(H_xM)=H_{f(x)}M'$ for all $x\in M$.
If $(M,HM,J)$ is of type $(n,k)$ and $(M',HM',J')$ of type $(n',k')$,
the existence of a $CR$ submersion implies that 
$n\geq n'$ and $k\geq k'$. 
\par
When $f\colon M@>>>M'$ is also a smooth fiber bundle, we say that
$f\colon (M,HM,J)@>>>(M',HM',J')$ is a {\it $CR$ fibration}. The fibers
are $CR$ manifolds of type $(n-n',k-k')$. 
\smallskip
A $CR$ {\it diffeomorphism} of $(M,HM,J)$ onto $(M',HM',J')$ is
a diffeomorphism $f\colon M@>>>M'$ such that both $f$ and $f^{-1}$ are
smooth $CR$ maps. The set of all $CR$ diffeomorphisms of
$(M,HM,J)$ onto itself ({\it $CR$ automorphisms}) is a group with
the composition operation.
\par
We say that $(M,HM,J)$ is a {\it homogeneous $CR$ manifolds} if there
is a Lie group of $CR$ automorphisms that acts transitively on $M$.
\par
Let $(M,HM,J)$ be a $CR$ manifold.
A vector field $X\in\Cal{C}^{\infty}(U,TM)$, defined on an open subset
$U$ of $M$, is an {\it infinitesimal $CR$ automorphism} if
the maps $\varphi_X(t)$ of the local $1$-parameter group of 
local transformations generated by $X$ are $CR$. This is equivalent
to the fact that 
$[X,\Cal{C}^\infty(U,T^{0,1}M)]\subset\Cal{C}^\infty(U,T^{0,1}M)$.
We say that $(M,HM,J;\o)$ is a {\it locally homogeneous $CR$ manifolds} 
at a point $\o\in M$ if, for each $v\in{T}_{\o}M$,
there is an infinitesimal $CR$ automorphism $X$, defined in an open
neighborhood $U$ of $\o$ in $M$, with $v=X(\o)$. 
\par\smallskip
Homogeneous $CR$ manifolds are locally $CR$ homogeneous: a homogeneous
$CR$ manifold $(M,HM,J)$ has a {\it real analytic} $CR$ structure
and therefore (see e.g. [AF]) admits a generic embedding 
$M\hookrightarrow \Hat{M} $ into a complex manifold $\Hat{M}$.
Then the Lie algebra $\frak{g}$
of the Lie group $\bold{G}$ that acts transitively
by $CR$ automorphisms on $M$ can be identified to a Lie algebra
of infinitesimal analytic $CR$ automorphisms defined on $U=M$. 
Each $X^*\in\Cal{C}^\infty(M,TM)$, corresponding to an $X\in\frak{g}$,
is the restriction of the real part of a holomorphic vector field 
$Z^*$,
defined on an open complex neighborhood $\Hat{U}$ of $M$ in $\Hat{M}$
(i.e. $X^*=\left[\roman{Re}\,Z^*\right]\left|_M\right.$;
see e.g. [BER \S{12.4}]).
\par\smallskip
The germs of infinitesimal $CR$ automorphisms of $(M,HM,J)$ at a point
$\o\in M$, with the Lie bracket, form a real Lie algebra 
$\frak{G}=\frak{G}(M,HM,J;\o)$.
We consider its complexification $\Hat{\frak{G}}=\frak{G}\oplus {i}\,
\frak{G}$ and denote by $\frak{Q}=\frak{Q}(M,HM,J;\o)$ 
the complex Lie subalgebra of
$\Hat{\frak{G}}$ consisting of all $Z\in\Hat{\frak{G}}$ with
$Z(\o)\in T^{0,1}_{\o}M$.
The fact that $\frak{Q}$ is actually
a complex Lie subalgebra of $\Hat{\frak{G}}$ is a consequence of the
formal integrability of the partial complex structure $J$.
\par\smallskip
When $(M,HM,J)$ is locally $CR$ homogeneous at $\o\in M$,
the pair $(\frak{G},\frak{Q})=\left(\frak{G}(M,HM,J;\o),\frak{Q}(
M,HM,J;\o)\right)$
completely determines the {\it germ} of 
the $CR$ manifold $(M,HM,J)$ at $\o$. Vice versa, if $\frak{g}$
is a {\it finite dimensional}
real Lie algebra and $\frak{q}$ a complex Lie subalgebra of
its complexification $\Hat{\frak{g}}$, 
the general construction\footnote{If $\frak{g}$ is a finite dimensional
real Lie algebra and $\frak{g}_+$ a real Lie subalgebra of $\frak{g}$,
we can find a germ $(M,\o)$ of analytic real manifold,
unique modulo germs of analytic diffeomorphisms, for which
there is a real Lie algebras homomorphism 
$\imath:\frak{g}@>>>\Cal{C}^\infty_{(\o)}(M,TM)$ with 
$\{\imath(X)(\o)\,|\, X\in\frak{g}\}=T_{\o}M$ and
$\frak{g}_+=\{X\in\frak{g}\,|\,\imath(X)(\o)=0\}$.}
of a germ $(M,\o)$ of homogeneous manifold associated
to the Lie algebra $\frak{g}$ and to its Lie subalgebra
$\frak{g}_+=\frak{g}\cap\frak{q}$ 
(cf. e.g. [Po, Ch.VII])
yields a unique, modulo local $CR$ diffeomorphisms, germ of 
locally homogeneous $CR$ manifold $(M,HM,J;\o)$ at $\o\in M$,
such that the complexification $\Hat{\imath}$ of the correspondence
$\imath:\frak{g}@>>>\Cal{C}^{\infty}_{(\o)}(M,TM)$
yields a homomorphism of complex Lie algebras:
$\Hat{\imath}:\Hat{\frak{g}}@>>>\Hat{\frak{G}}(M,HM,J;\o)$
with 
\form{\imath(\frak{g})\subset\frak{G}(M,HM,J;\o)\,,\quad
\Hat{\imath}(\frak{q})
\subset\frak{Q}(M,HM,J;\o)}\edef\FormBD{\number\q.\number\t}
and for which the induced map on the quotients
$$\frak{g}/(\frak{g}\cap\frak{q}) @>>>
\frak{G}(M,HM,J;\o)/\left(\frak{G}(M,HM,J;\o)
\cap\frak{Q}(M,HM,J;\o)\right)$$ is an isomorphism. 
In this case we say that the germ of $(M,HM,J)$ at $\o$
is {\it associated} to the pair $(\frak{g},\frak{q})$.
\par\smallskip
These remarks led to the introduction in [MN4] of the 
abstract notion of a
$CR$ algebra. 
A {\it $CR$ algebra $(\frak{g},\frak{q})$} is the pair consisting of a
{\it real} Lie algebra $\frak{g}$ and of a {\it complex} Lie
subalgebra $\frak{q}$ of its complexification $\Hat{\frak{g}}$,
such that the quotient $\frak{g}/\left(\frak{g}\cap\frak{q}\right)$
is a finite dimensional real linear space. Note that we do not require
that $\frak{g}$ is finite dimensional.
The intersection $\frak{g}_+=\frak{g}\cap\frak{q}$ is the {\it isotropy}
of $(\frak{g},\frak{q})$. Let $\Cal{H}_+=\{\roman{Re}\,Z\,|\, Z\in\frak{q}\}$
and denote by $\bar{Z}$ the conjugate of $Z\in\Hat{\frak{g}}$
with respect to the real form $\frak{g}$. \par
A $CR$ algebra $(\frak{g},\frak{q})$ is\,:
\par\smallskip\noindent
{\it totally real} \quad if 
$\Cal{H}_+=\frak{g}_+$,
\par\smallskip\noindent
{\it totally complex} \quad if $\Cal{H}_+=\frak{g}$,
\par\smallskip\noindent
{\it fundamental} \quad if $\Cal{H}_+$ generates $\frak{g}$ as a real Lie 
algebra,
\par\smallskip\noindent
{\it transitive, or effective} \quad if $\frak{g}_+$ does not
contain any nonzero ideal of $\frak{g}$,
\par\smallskip\noindent
{\it ideal nondegenerate} \quad if all ideals of $\frak{g}$ 
contained in $\Cal{H}_+$ are contained in $\frak{g}_+$,
\par\smallskip\noindent
{\it weakly nondegenerate} \quad if there is no complex Lie subalgebra 
$\frak{q}'$ of
$\Hat{\frak{g}}$ with\,: \par\centerline{$
\frak{q}\subsetneq \frak{q}'\subset
\frak{q}+\bar{\frak{q}}$,}
\par\smallskip\noindent
{\it strictly nondegenerate} \quad if $\frak{g}_+=\{X\in\Cal{H}_+\,|\,
[X,\Cal{H}_+]\subset\Cal{H}_+\}$.
\par\medskip
Clearly\,:\par\smallskip\centerline{
$\text{strictly nondegenerate}\Longrightarrow
\text{weakly nondegenerate}\Longrightarrow\text{ideal nondegenerate}\,.$}
\par\smallskip
Fundamentality of $(\frak{g},\frak{q})$ is equivalent to the fact
that the associated germ of homogeneous $CR$ manifold
$(M,HM,J;\o)$ is of {\it finite type in the sense
of} [BlG], i.e. that the smallest involutive distribution
of tangent vectors containing $HM$ also contains $T_{\o}M$.\par
Strict and weak nondegeneracy  hold, or do not hold,
for all $CR$ algebras
that are associated to the same germ $(M,HM,J;\o)$ of 
locally homogeneous $CR$ manifold.
They correspond indeed to the nondegeneracy of the
(vector valued) Levi form and of its higher order analog, respectively
(see e.g. [MN4, \S{13}]). In particular,
{\it weak nondegeneracy}  
at a point $\o\in M$ of a (germ of) $CR$ manifold $(M,HM,J)$
means that, for every 
$L\in\Cal C^{\infty}(M,T^{1,0}M)$ with $L(\o)\neq 0$, there exist
finitely many vector fields 
$\bar Z_1,\,\hdots ,\,\, \bar Z_m\in\Cal C^\infty(M,T^{0,1}M)$
such that $[L\,,\,\,\bar Z_1\,,\,\,\hdots\, ,
\,\,\bar Z_m](\o)\notin T_{\o}^{1,0}M\oplus
T_{\o}^{0,1}M$. \par
When $(\frak{g},\frak{q})$ defines at
$\o$ a germ of homogeneous $CR$ manifold, the two notions of 
weak nondegeneracy, the one for $CR$ algebras and the one above for
$CR$ manifolds, coincide and also coincide with the {\it holomorphic
nondegeneracy} of [BER].
\par
We have the following (see [AMN, Proposition 4.1])\,:
\prop{Let $(M,HM,J)$ and $(M',HM',J')$ be $CR$ manifolds. Assume that
$M'$ is locally embeddable and that
there exists a $CR$ fibration 
$(M,HM,J)@>\pi>> (M',HM',J')$, with totally complex
fibers of positive dimension. Then $M$ is weakly degenerate.\qed}

Differently,
both ideal degenerate and ideal nondegenerate
$CR$ algebras
may correspond to the same (weakly degenerate) germ of 
locally homogeneous
$CR$ manifold.
\par\smallskip
From [MN4, Theorem 9.1] we know that {\it if $(\frak{g},\frak{q})$
is fundamental, effective, and ideal nondegenerate, then 
$\frak{g}$ is finite dimensional}.
\par
From this result we deduce the following\,:
\thm{Let $(\frak{g},\frak{q})$ be a fundamental effective $CR$ algebra.
Then there exist a germ of homogeneous complex manifold $(\Hat{M},\o)$ 
at a point $\o$, and a germ of
homogeneous generic $CR$ submanifold $(M,HM,J;\o)$ 
of $(\Hat{M},\o)$ at $\o$,
with associated $CR$ algebra $(\frak{g},\frak{q})$.}
\dimo First we note that the statement holds true 
when
$\frak{g}$ is finite dimensional: by [Po, Ch.VII] there is a germ of
homogeneous complex manifold $(\Hat{M},\o)$ at $\o$ corresponding
to the complex Lie algebra $\Hat{\frak{g}}$ and to its complex
Lie subalgebra $\frak{q}$; the inclusion $\frak{g}\hookrightarrow
\Hat{\frak{g}}$ yields the embedding into $(\Hat{M},\o)$ 
of a germ of
homogeneous $CR$ manifold $(M,HM,J;\o)$ at $\o$, 
corresponding to
the pair $(\frak{g},\frak{q})$.
\par
Consider now the general case. We keep the notation introduced above.
By [MN4, Lemma 7.2] there is a largest
ideal $\frak{a}$ of $\frak{g}$ contained in $\Cal{H}_+$.
By [MN4, Theorem 9.1], $\frak{g}/\frak{a}$ is finite dimensional
and by the first part of the proof there is a germ of complex homogeneous
manifold
$(\Hat{N},\o)$ at $\o$, 
and a germ of generic $CR$ submanifold $(N,HN,J_N;\o)$ of
$\Hat{N}$ at $\o$, associated
to the pair $\left(\frak{g}/\frak{a},\frak{q}/(\frak{q}\cap\Hat{
\frak{a}})\right)$ (here $\Hat{\frak{a}}=\Bbb{C}\otimes_{\Bbb{R}}\frak{a}$
is the complexification of $\frak{a}$ in $\Hat{\frak{g}}$).
If $2d$ is the real dimension of $\frak{a}/(\frak{a}\cap\frak{g}_+)$,
we can take $\Hat{M}=\Hat{N}\times\Bbb{C}^d$ and, likewise,
$M=N\times \Bbb{C}^d$, with $HM=HN\times T(\Bbb{C}^d)$ and
$J=J_N\times J_{\Bbb{C}^d}$. Then 
$\left(M,HM,J;(\o,0)\right)$ is associated to
$(\frak{g},\frak{q})$.\enddimo
Note that the ideal nondegeneracy of $(\frak{g},\frak{q})$ 
implies
that {\it all 
ideals $\frak{x}$ of the complex Lie algebra $\Hat{\frak{g}}$ that
are contained in $\frak{q}$ are contained in $\frak{q}\cap\bar{\frak{q}}$.}
Indeed, if $\frak{x}$ is a (complex) ideal 
of $\Hat{\frak{g}}$ contained in $\frak{q}$, then
$\frak{a}=\left(\frak{x}+\bar\frak{x}\right)\cap\frak{g}$ is an ideal
of $\frak{g}$ contained in 
$\Cal{H}_+=\left(\frak{q}+\bar{\frak{q}}\right)\cap\frak{g}$,
and $\frak{a}\subset\frak{g}_+=\frak{q}\cap\frak{g}$ implies
that $\frak{x}+\bar{\frak{x}}\subset\frak{q}\cap\bar{\frak{q}}$.
\par
\smallskip
Let $(\frak{g},\frak{q})$ be a $CR$ algebra with a
finite dimensional $\frak{g}$. 
We denote by $\Hat{\bold{G}}$ the connected
and simply connected complex Lie group with Lie algebra
$\Hat{\frak{g}}$ and by $\bold{Q}$ its
analytic Lie subgroup, 
generated by $\frak{q}$.
Let $\bold{G}$ be the analytic subgroup of $\Hat{\bold{G}}$
with Lie algebra $\frak{g}$ and set 
$\bold{G}_+=\bold{G}_+(\frak{g},\frak{q}):=
\bold{Q}\cap\bold{G}$. This is a Lie subgroup of
$\bold{G}$ with Lie subalgebra $\frak{g}_+$.
Denote by $\Tilde{\bold{G}}$ a connected and simply connected
Lie group with Lie algebra $\frak{g}$ and by 
$\Tilde{\bold{G}}_+=\Tilde{\bold{G}}_+(\frak{g},\frak{q})$
its analytic Lie subgroup with Lie subalgebra $\frak{g}_+$.
Since $\Hat{\bold{G}}$ is simply connected, the conjugation 
$\sigma$ of
$\Hat{\frak{g}}$ with respect its real form $\frak{g}$ defines
an antiholomorphic involution, still denoted by ${\sigma}$,
in $\Hat{\bold{G}}$. Thus $\bold{G}$, being the connected
component of the identity in the set of fixed points of ${\sigma}$,
is closed in $\Hat{\bold{G}}$. We have the
implications\,:\par\smallskip\centerline{
$\bold{Q}\;\text{closed in}\;\Hat{\bold{G}}\; \Longrightarrow\;
\bold{G}_+\;\text{closed in}\;\bold{G}\,,\quad
\bold{G}_+\;\text{closed in}\;\bold{G} \Longrightarrow\;
\Tilde{\bold{G}}_+\;\text{closed in}\;\Tilde{\bold{G}}\,.$}\par\smallskip
When $\tilde{\bold{G}}_+$ is closed in $\tilde{\bold{G}}$, 
we can uniquely define a $\Tilde{\bold{G}}$-homogeneous $CR$ manifold
$\tilde{M}(\frak{g},\frak{q})=(\Tilde{M},H\Tilde{M},\tilde{J})$,
where the underlying smooth manifold $\tilde{M}$ is the 
$\tilde{\bold{G}}$-homogeneous
space $\tilde{\bold{G}}/\tilde{\bold{G}}_+$, and $(\frak{g},\frak{q})$
is associated to the germ $(\Tilde{M},H\Tilde{M},\tilde{J};\o)$
at the base point $\o=\,e\,\tilde{\bold{G}}_+$.\par
Likewise, for a closed $\bold{G}_+\subset\bold{G}$, we define the
$\bold{G}$-homogeneous $CR$ manifold $M(\frak{g},\frak{q})=(M,HM,J)$
with $M=\bold{G}/\bold{G}_+$ and $(\frak{g},\frak{q})$ associated to
$(M,HM,J;\o)$ for $\o=\,e\,\bold{G}_+$.
\par
If $\bold{Q}$ is closed, 
$\Hat{M}=\Hat{M}(\Hat{\frak{g}},\frak{q}):=\Hat{\bold{G}}/\bold{Q}$ 
is a $\Hat{\bold{G}}$-homogeneous
complex manifold and $M(\frak{g},\frak{q})$ can be identified, its
partial complex structure being that of
a generic $CR$ submanifold of $\Hat{M}$, 
to the orbit of $\bold{G}$ through the base point $\o=\,e\,\bold{Q}$
of $\Hat{M}$\,.

\par
Our canonical choice of $M(\frak{g},\frak{q})$ aims to obtain 
a homogeneous $CR$ manifold with a 
{\it good} embedding in a homogeneous {\it complex} manifold
$\Hat{M}=\Hat{M}(\Hat{\frak{g}},\frak{q})$.
\par\smallskip
A morphism of $CR$ algebras 
$(\frak{g},\frak{q})@>{\varphi}>>(\frak{g}',\frak{q}')$ is a homomorphism
of real Lie algebras
$\varphi:\frak{g}@>>>\frak{g}'$, with $\Hat{\varphi}(\frak{q})\subset\frak{q}'$.
(We shall consistently use "{\it hat}\," to indicate complexification:
e.g.
$\Hat{\varphi}:\Hat{\frak{g}}@>>>\Hat{\frak{g}}'$ 
is the complexification of $\varphi:\frak{g}@>>>\frak{g}'$).
It is called\,:\roster
\item"$\bullet$"
{\it a $CR$ immersion} if the quotient map $[\varphi]:\frak{g}/\frak{g}_+
@>>>\frak{g}'/\frak{g}'_+$ is injective and $\Hat{\varphi}^{-1}(\frak{q}')=
\frak{q}$;
\item"$\bullet$"
{\it a $CR$ submersion} if both $[\varphi]:\frak{g}/\frak{g}_+
@>>>\frak{g}'/\frak{g}'_+$ and $[\Hat{\varphi}]:\frak{q}/\Hat{\frak{g}}_+
@>>>\frak{q}'/\Hat{\frak{g}}_+'$ are onto;
\item"$\bullet$"
{\it a local $CR$ isomorphism} if both $[\varphi]:\frak{g}/\frak{g}_+
@>>>\frak{g}'/\frak{g}'_+$ and $[\Hat{\varphi}]:\frak{q}/\Hat{\frak{g}}_+
@>>>\frak{q}'/\Hat{\frak{g}}_+'$ are isomorphisms;
\item"$\bullet$"
{\it a $CR$ isomorphism} if $\varphi$ is an isomorphism of $CR$ algebras
with $\Hat{\varphi}(\frak{q})=\frak{q}'$.\endroster
We quote from [MN4]\,:
\prop{Let $(\frak{g},\frak{q})@>{\varphi}>>(\frak{g}',\frak{q}')$ be
a morphism of $CR$ algebras, with $\frak{g}$ and $\frak{g}'$ finite
dimensional. Let $(M,HM,J)$ and $(M',HM',J')$ be the germs of homogeneous
$CR$ manifolds at $\o\in M$, $\o'\in M'$, associated to
$(\frak{g},\frak{q})$, $(\frak{g}',\frak{q}')$, respectively.
Then there is a unique germ of smooth $CR$ map
$\Phi:(M,HM,J)@>>>(M',HM',J')$ with $\Phi(\o)=\o'$
such that
$d\Phi_{\o}(\imath(X))=\imath'(\varphi(X))$.
Here $\imath$, $\imath'$ are the homomorphisms
of Lie algebras $\imath:\frak{g}@>>>\frak{G}(M,HM,J)$ and
$\imath':\frak{g}:@>>>\frak{G}(M',HM',J')$ of $(\FormBD)$.
\par
The germ $\Phi$
of smooth $CR$ map is a $CR$ immersion, submersion, diffeomorphism
if and only if the corresponding morphism $\varphi$ of $CR$ algebras is
a $CR$ immersion, a $CR$ submersion, a local $CR$ isomorphism, respectively.
\par
Let $\tilde{\bold{G}}$ and $\tilde{\bold{G}}'$ be the connected and
simply connected real Lie groups with Lie algebras $\frak{g}$ and
$\frak{g}'$, respectively. If the analytic subgroup
$\tilde{\bold{G}}_+$ of $\tilde{\bold{G}}$ with Lie algebra
$\frak{g}_+=\frak{q}\cap\frak{g}$ and the analytic subgroup
$\tilde{\bold{G}}_+'$ of $\tilde{\bold{G}}'$ with Lie algebra
$\frak{g}_+'=\frak{q}'\cap\frak{g}'$
are both closed, then there is a unique smooth $CR$ map
$\tilde{\Phi}:\tilde{M}(\frak{g},\frak{q})@>>>\tilde{M}(\frak{g}',\frak{q}')$ 
that makes the following diagram commutes\,:
$$\CD
\frak{g}@>{\exp}>>\tilde{\bold{G}}@>>>\tilde{M}(\frak{g},\frak{q})\\
@V{\varphi}VV  @VVV  @VV{\tilde{\Phi}}V\\
\frak{g}'@>{\exp}>>\tilde{\bold{G}}'@>>>\tilde{M}(\frak{g}',\frak{q}')
\endCD$$
The map $\Tilde{\Phi}$ is a $CR$ immersion, a $CR$ submersion or a local 
$CR$ diffeomorphism if and only if the corresponding $CR$ morphism of
$CR$ algebras $\varphi$ is
a $CR$ immersion, a $CR$ submersion or a local 
$CR$ isomorphism, respectively.\par
Let $\Hat{\bold{G}}$ and $\Hat{\bold{G}}'$ be the connected and simply
connected complex Lie groups with Lie algebras $\Hat{\frak{g}}$ and
$\Hat{\frak{g}}'$, respectively. Let 
$\bold{G},\bold{Q}\subset\Hat{\bold{G}}$
and $\bold{G}',\bold{Q}'\subset\Hat{\bold{G}}'$ 
be the analytic subgroups with
Lie algebras $\frak{g},\frak{q}$ and $\frak{g}',\frak{q}'$, respectively.
If $\bold{G}_+=\bold{Q}\cap\bold{G}$ and $\bold{G}'=\bold{Q}'\cap
\bold{G}'$ are closed,
then there is a unique 
smooth $CR$ map $\Phi:
{M}(\frak{g},\frak{q})@>>>{M}(\frak{g}',\frak{q}')$
such that the diagram\,:
$$\CD
 \tilde{M}(\frak{g},\frak{q})@>{\Tilde{\Phi}}>>\tilde{M}(\frak{g}',\frak{q}')
\\
@VVV @VVV \\
{M}(\frak{g},\frak{q})@>\Phi>>{M}(\frak{g}',\frak{q}')
\endCD
$$
where the vertical arrows are the natural projections from the universal
coverings, commutes.\par
The map ${\Phi}$ is a $CR$ immersion, a $CR$ submersion or a local 
$CR$ diffeomorphism if and only if the corresponding morphism of
$CR$ algebras $\varphi$ is
a $CR$ immersion, a $CR$ submersion or a local 
$CR$ isomorphism, respectively.\par
If 
$\bold{Q}\subset\Hat{\bold{G}}$ and 
$\bold{Q}'\subset\Hat{\bold{G}}'$
are closed, the map 
${M}(\frak{g},\frak{q})@>\Phi>>{M}(\frak{g}',\frak{q}')$
is the restriction of the holomorphic map
$\Hat{\Phi}:\Hat{M}=\Hat{\bold{G}}/\bold{Q} @>>>
\Hat{M}'=\Hat{\bold{G}}'/\bold{Q}'$ 
defined by the commutative diagram\,:
$$\CD
\Hat{\frak{g}} @>{\exp}>> \Hat{\bold{G}} @>>> \Hat{M}\\
@V{\Hat{\varphi}}VV @VVV @VV{\Hat{\Phi}}V \\
\Hat{\frak{g}}' @>{\exp}>> \Hat{\bold{G}}' @>>> \Hat{M}'
\endCD$$
where the central vertical arrow is the homomorphism of complex
connected simply connected Lie algebras defined by the homomorphism
$\Hat{\varphi}:\Hat{\frak{g}}@>>>\Hat{\frak{g}}'$ of their Lie algebras.
\qed}\par\medskip
To discuss, later on, the structure of the fibers of some $CR$ fibrations,
we briefly introduce the notion of
{\it semidirect sum}
of $CR$ algebras.\par
Let $(\frak{g}_1,\frak{q}_1)$, $(\frak{g}_2,\frak{q}_2)$ be
$CR$ algebras, and assume that $\frak{g}_2$ has a structure of
$\frak{g}_1$-module and that $\frak{q}_2$ is
a ${\frak{q}}_1$-module for the restriction of the complexification of
the action of $\frak{g}_1$ on $\frak{g}_2$.
Then $\frak{q}=\frak{q}_1\rtimes\frak{q}_2$ (semidirect sum) is a complex Lie
subalgebra of the complexification of the semidirect sum
$\frak{g}=\frak{g}_1\rtimes\frak{g}_2$, and the $CR$ algebra
$(\frak{g},\frak{q})= 
(\frak{g}_1\rtimes\frak{g}_2,\frak{q}_1\rtimes\frak{q}_2)$
is called the {\it semidirect sum} of the
$CR$ algebras $(\frak{g}_1,\frak{q}_1)$ and $(\frak{g}_2,\frak{q}_2)$ \,:
\form{(\frak{g},\frak{q})=
(\frak{g}_1\rtimes\frak{g}_2,\frak{q}_1\rtimes\frak{q}_2)=
(\frak{g}_1,\frak{q}_1)\rtimes (\frak{g}_2,\frak{q}_2)
\quad\text{(semidirect sum)}\,.}
We shall assume that
$\frak{g}_1$ and $\frak{g}_2$ are finite dimensional.
Denote by\,:\par
$\Hat{\bold{G}}$, $\Hat{\bold{G}}_1$, $\Hat{\bold{G}}_2$
the connected and simply connected complex Lie groups with Lie algebras
$\Hat{\frak{g}}$, $\Hat{\frak{g}}_1$, $\Hat{\frak{g}}_2$, respectively;
\par
${\bold{G}}$, ${\bold{G}}_1$, ${\bold{G}}_2$
the analytic real subgroups of the corresponding complex connected
Lie groups $\Hat{\bold{G}}$, $\Hat{\bold{G}}_1$, $\Hat{\bold{G}}_2$,
with Lie algebras
${\frak{g}}$, ${\frak{g}}_1$, ${\frak{g}}_2$, respectively;
\par
$\bold{Q}\subset\Hat{\bold{G}}$, 
$\bold{Q}_1\subset\Hat{\bold{G}}_1$, $\bold{Q}_2\subset\Hat{\bold{G}}_2$
the Lie subgroups corresponding to the Lie subalgebras
$\frak{q}\subset\Hat{\frak{g}}$, 
$\frak{q}_1\subset{\frak{g}}_1$, 
$\frak{q}_2\subset{\frak{g}}_2$, respectively; 
\par
$\Tilde{\bold{G}}$, $\Tilde{\bold{G}}_1$, $\Tilde{\bold{G}}_2$
connected and simply connected real Lie groups
with Lie algebras
${\frak{g}}$, ${\frak{g}}_1$, ${\frak{g}}_2$, respectively;
\par
${\bold{G}}_+=\bold{Q}\cap{\bold{G}}$,
${\bold{G}_1\,}_+=\bold{Q}\cap{\bold{G}}_1$,
${\bold{G}_2\,}_+=\bold{Q}\cap{\bold{G}}_2$; \par
$\Tilde{\bold{G}\,}_+\subset
\Tilde{\bold{G}}$,
$\Tilde{\bold{G}_1\,}_+\subset
\Tilde{\bold{G}}_1$, 
${\Tilde{\bold{G}}_2\,}_+\subset\Tilde{\bold{G}}_2$ 
the analytic subgroups corresponding
to the Lie subalgebras
$\frak{g}_+=\frak{q}\cap\frak{g}$, 
${\frak{g}_1}_+=\frak{q}_1\cap{\frak{g}}_1$, 
${\frak{g}_2}_+=\frak{q}_2\cap{\frak{g}}_2$, respectively.\par\smallskip
Let
$(M,HM,J;\o)$,
$(M_1,HM_1,J_1;\o_1)$,$(M_2,HM_2,J_2;\o_2)$ be the
germs of locally homogeneous $CR$ manifolds associated to
the $CR$ algebras
$(\frak{g},\frak{q})$, $(\frak{g}_1,\frak{q}_1)$,
$(\frak{g}_2,\frak{q}_2)$,
respectively.
\par
We obtain\,:
\thm{The diffeomorphism 
$\bold{G}_1\times\bold{G}_2\ni (g_1,g_2)@>>> g_1g_2\in\bold{G}_1\rtimes
\bold{G}_2$ defines a germ of $CR$ diffeomorphism\,:
$$(M_1,HM_1,J_1;\o_1)\times
(M_2,HM_2,J_2;\o_2)@>>>(M,HM,J;\o)\,.$$
If ${\Tilde{\bold{G}}_1\,}_+$ and 
${\Tilde{\bold{G}}_2\,}_+$ are closed, then
$\Tilde{\bold{G}}_+$ is closed and we obtain a
global $CR$ diffeomorphism\,:
$$\Tilde{M}(\frak{g}_1,\frak{q}_1)\times
\Tilde{M}(\frak{g}_2,\frak{q}_2) @>>> \Tilde{M}(\frak{g},\frak{q})\,.
 $$
If ${\bold{G}_1\,}_+$ and ${\bold{G}_2\,}_+$ are closed, then
$\bold{G}_+$ is closed and we obtain a
global $CR$ diffeomorphism\,:
$${M}(\frak{g}_1,\frak{q}_1)\times
{M}(\frak{g}_2,\frak{q}_2) @>>> {M}(\frak{g},\frak{q})\,.\tag $*$ $$
When $\bold{Q}_1$ and $\bold{Q}_2$ are closed, also
$\bold{Q}$ is closed and the map $(*)$ is the restriction of a
biholomorphic map
$$\left(\Hat{\bold{G}}_1/\bold{Q}_1\right)\times
\left(\Hat{\bold{G}}_2/\bold{Q}_2\right)@>>>
\Hat{\bold{G}}/\bold{Q}\, .\qquad\qquad\qed$$}
Let $\frak g$ be a real Lie algebra and $\frak q,\frak q'$ 
complex subalgebras of its complexification $\Hat{\frak g}$,
with $\frak q\subset\frak q'$.
Then the identity map in $\frak g$ and the inclusion $\frak q\hookrightarrow
\frak q'$ define a 
$\frak{g}$-equivariant $CR$ submersion of $CR$ algebras 
\form{(\frak g,\frak q)@>>>(\frak g,\frak q')\, .}
\edef\formlacaa{\number\q.\number\t}
\par
If $(M,HM,J;\o)$ and $(M',HM',J';\o')$ are germs
of locally homogeneous $CR$ manifolds with associated $CR$ algebras
$(\frak{g},\frak{q})$ and $(\frak{g},\frak{q}')$, respectively,
then the identity map $\frak{g}@>>>\frak{g}$ defines, by passing to
the quotients, the differential of a $CR$ map 
$\pi_{(\o)}:(M,HM,J;\o)@>>>(M',HM',J';\o')$ that is
locally $\bold{G}$-equivariant for a (connected) real Lie group with
Lie algebra $\frak{g}$.
\par 
Let $\bold{G}$ be a 
connected Lie group with Lie algebra $\frak{g}$ and assume
that there are two $\bold{G}$-homogeneous
$CR$ manifolds
$(M,HM,J)$, $(M',HM',J')$ that are
associated to $(\frak{g},\frak{q})$ 
at some
$\o\in M$ and to $(\frak{g},\frak{q}')$ 
at some
$\o'\in M'$, respectively.
Then there is a unique $\bold{G}$-equivariant $CR$ map
$\pi:(M,HM,J)@>>>(M',HM',J')$ with $\pi(\o)=\o'$.
\par
In general $\pi_{(\o)}$ (and $\pi$, when defined) are smooth,
but not $CR$, $\bold{G}$-equivariant fibrations:
a necessary and sufficient condition for $(\formlacaa)$ to be
a $\frak g$-equivariant $CR$ fibration, and hence for $\pi_{(\o)}$ 
(and $\pi$, when defined) to be $\bold{G}$-equivariant local 
(resp. global) fibrations
is that (see [MN4, Lemma 5.1])
\form{\frak q'=\frak q+\Hat{\frak g}'_+\, .}

We call $(\frak g,\frak q')$ the {\it basis} of the fibration
$(\formlacaa)$. 
The {\it fiber} of $(\formlacaa)$ is the $CR$ algebra
$(\frak g'_+,\frak q'')$ where 
$\frak q''=\Hat{\frak g}'_+\cap\frak q$. It is a $CR$ algebra
associated to the germ 
$\left(F,HF,J\left|_{HF}\right.;\o\right)$, where
$\left(F;\o\right)=\left(\pi^{-1}_{(\o)}(\o');\o\right)$, and
the germ of partial complex structure 
$\left(HF,J\left|_{HF}\right.; \o\right)$
that is characterized by
requiring that the smooth embedding $(\pi^{-1}(\o');\o)
\hookrightarrow (M;\o)$ is a $CR$ immersion.
\par
We know that $(\formlacaa)$ is always a $CR$ fibration, with a totally
complex fiber, when
$\frak q\subset\frak q'\subset\frak q+\overline{\frak q}$\,: indeed in this
case $\frak q'=\frak q+\Hat{\frak g}_+'$.\par
From [MN4, \S{5},] we have\,:
\prop{Let $(\frak{g},\frak{q})$ be a $CR$ algebra. Then there exist\,:
\par\smallskip\centerline{$\matrix\format\l\\
\text{a largest ideal $\frak{i}$ 
of $\frak{g}$ with $\frak{i}\subset\frak{g}_+$}\\
\text{a largest ideal $\frak{a}$ 
of $\frak{g}$ with $\frak{a}\subset\Cal{H}_+$}\\
\text{a largest complex Lie subalgebra $\frak{q}'$ 
of $\Hat{\frak{g}}$ with $\frak{q}\subset
\frak{q}'\subset\frak{q}+\bar{\frak{q}}$}\\
\text{a smallest complex Lie subalgebra $\frak{q}''$ 
of $\Hat{\frak{g}}$ with $\frak{q}+\bar{\frak{q}}
\subset \frak{q}''\,.$}
\endmatrix$}\par\smallskip
We have $\frak{i}\subset\frak{a}\subset\frak{q}'\subset\frak{q}''$ and
$\frak{q}''=\bar{\frak{q}}\,''=\Hat{\frak{g}}\,''$
for a real Lie subalgebra $\frak{g}''$ of $\frak{g}$.\par
The identity in $\frak{g}$ defines
$\frak{g}$-equivariant $CR$ fibrations
$(\frak{g},\frak{q})@>>>(\frak{g},\frak{q}')@>>>(\frak{g},\frak{q}\,'')$,
where $(\frak{g},\frak{q}')$ is weakly nondegenerate and
$(\frak{g},\frak{q}\,'')$ is totally real. For all complex Lie subalgebras
$\frak{f}$ of $\Hat{\frak{g}}$ with 
$\frak{q}\subset\frak{f}\subsetneq\frak{q}'$, the $CR$ algebra
$(\frak{g},\frak{f})$ is weakly degenerate. 
For all complex Lie subalgebras $\frak{f}$ with 
$\frak{q}\subset\frak{f}\subset\frak{q}'$, the $\frak{g}$-equivariant
map $(\frak{g},\frak{q})@>>>(\frak{g},\frak{f})$ is a $CR$ fibration
with a totally complex fiber.\par
The $CR$ algebra
$(\frak{g}'',\frak{q})$ is fundamental and, for all
real Lie subalgebras
$\frak{l}$ of $\frak{g}$ with 
$\frak{g}''\subsetneq{\frak{l}}\subset\frak{g}$
the $CR$ algebra $(\frak{l},\frak{q})$ is not fundamental.\qed}
\par\edef\PropBD{\number\q.\number\x}

\edef\paragB{\number\q}

\se{Parabolic subalgebras and complex flag manifolds}
In this section we collect the notions on complex parabolic
subalgebras and fix that notation
that will be utilized throughout the paper.
This is mostly a review of classical results, for which  
general references are [B1, Ch.IV,\S{2}.6, Ch.VI,\S{1}], 
[B2, Ch.VIII,\S{3}],
[Kn, Ch.VII] [Wa, Ch.1], [W].\par
Let $\Hat{\frak{g}}$ be a complex Lie algebra. A 
maximal solvable complex Lie subalgebra $\frak{b}$ of 
$\Hat{\frak{g}}$ is called a {\it Borel}, or {\it minimal parabolic}
complex Lie subalgebra of $\Hat{\frak{g}}$. A complex Lie subalgebra
$\frak{q}$ of $\Hat{\frak{g}}$ is {\it parabolic} if it contains a
complex Borel subalgebra $\frak{b}$ of $\Hat{\frak{g}}$.
\par
For our purposes, it will be sufficient to consider the case
of a semisimple $\Hat{\frak{g}}$. Thus from now on we shall assume
that $\Hat{\frak{g}}$ is a semisimple complex Lie algebra.\par
A parabolic subalgebra $\frak{q}$ of $\Hat{\frak{g}}$ contains
a complex Cartan subalgebra $\Hat{\frak{h}}$ of 
$\Hat{\frak{g}}$. Let $\Cal{R}=\Cal{R}(\Hat{\frak{g}},\Hat{\frak{h}})$
be the root system of $\Hat{\frak{g}}$ with respect to 
$\Hat{\frak{h}}$. 
We denote by $\frak{h}_{\Bbb{R}}$ the real form
of $\Hat{\frak{h}}$ on which all roots
are real
valued. Thus $\Cal{R}$ is a subset of the real dual
space $\frak{h}^*_{\Bbb{R}}$.
The Killing form $\kappa_{\Hat{\frak{g}}}$ of
$\Hat{\frak{g}}$ restricts to a real positive scalar product
in $\frak{h}_{\Bbb{R}}$. We shall write 
$(A|B)=\kappa_{\Hat{\frak{g}}}(A,B)$ for $A,B\in\frak{h}_{\Bbb{R}}$.
We set also
$(\xi|\eta)=(T_{\xi}|T_{\eta})$ for $\xi,\eta\in\frak{h}_{\Bbb{R}}^*$
and $(T_{\xi}|A)=\xi(A)$, $(T_{\eta}|A)=\eta(A)$ for all 
$A\in\frak{h}_{\Bbb{R}}$ (dual scalar product in~$ \frak{h}_{\Bbb{R}}^*$).
Roots $\alpha,\beta\in\Cal{R}$ for which $\alpha\pm\beta\notin\Cal{R}$
are called {\it strongly orthogonal}. Note that strongly orthogonal roots
are also orthogonal for the scalar product  in $\frak{h}_{\Bbb{R}}^*$\,.
\par
An element $H\in\Hat{\frak{h}}$ is {\it regular}
if $\alpha(H)\neq 0$ for all $\alpha\in\Cal{R}$. Denote by
$\frak{C}(\Cal{R})$ the set of the Weyl chambers
of $\Cal{R}$. They are the connected components of the set
of regular elements of $\frak{h}_{\Bbb{R}}$. For
$C\in\frak{C}(\Cal{R})$, and $H\in C$, the set
$\Cal{R}^+(C)=\{\alpha\in\Cal{R}\,|\, \alpha(H)>0\}$ is independent
of the choice of $H\in{C}$\,: 
it is called the set of {\it positive roots}
with respect to $C$. The set 
$\Cal{R}^-(C)=\Cal{R}^+({C}^{\roman{opp}})$, 
for ${C}^{\roman{opp}}=\{-H\,|\, H\in C\}$,
is the complement of $\Cal{R}^+(C)$ in $\Cal{R}$ and is called the
set of {\it negative roots} with respect to $C$. A Weyl chamber $C$
also defines a partial order relation ''$\prec_C$'' in $\frak{h}^*_{\Bbb{R}}$,
by\,: \par\smallskip\centerline{
$\eta\prec_C \xi$ if $\eta(A)<\xi(A)$ for all $A\in C$.}\par\smallskip
 Thus
$\Cal{R}^+(C)=\{\alpha\in\Cal{R}\,|\, \alpha\succ_C 0\}$.
\par\smallskip
With
$\Hat{\frak{g}}^{\,\alpha}=\{X\in\Hat{\frak{g}}\,|\,
[H,X]=\alpha(H)\,X\;\forall H\in\Hat{\frak{h}}\}$, we set
\form{\Cal{Q}=\{\alpha\in\Cal{R}\,|\,\Hat{\frak{g}}^{\alpha}
\subset\frak{q}\}\,.}
This is the {\it parabolic set} associated to $\frak{q}$ and
$\Hat{\frak{h}}$. Parabolic sets of roots are abstractly defined by
the two conditions\,:
$$\align
&\alpha,\beta\in\Cal{Q}\,,\; \alpha+\beta\in\Cal{R}\,\;
\Longrightarrow \alpha+\beta\in\Cal{Q}\quad\text{(closedness)}\tag i \\
&\Cal{Q}\,\cup\, {\Cal{Q}}^{\roman{opp}}\,=\,\Cal{R} \quad\text{where}\quad
{\Cal{Q}}^{\roman{opp}}
=\{-\alpha\,|\,\alpha\in\Cal{Q}\}\,. \qquad\qquad\tag ii
\endalign
$$
Condition $(\roman{ii})$ is equivalent to the fact that 
$\Cal{Q}\supset\Cal{R}^+(C)$ for some Weyl chamber $C\in\frak{C}(\Cal{R})$.
We have
\form{\frak{q}=\frak{q}_{\Cal{Q}}=
\Hat{\frak{h}}\oplus\dsize\sum_{\alpha\in\Cal{Q}}{
\Hat{\frak{g}}^{\,\alpha}}}
and the correspondence $\Cal{Q} \longleftrightarrow \frak{q}_{\Cal{Q}}$
is one-to-one between parabolic subsets of $\Cal{R}$
and parabolic subalgebras of $\Hat{\frak{g}}$ containing $\Hat{\frak{h}}$.
\par
Given a parabolic set $\Cal{Q}\subset\Cal{R}$ we set\,:
\form{\Cal{Q}^r=\{\alpha\in\Cal{Q}\,|\, -\alpha\in\Cal{Q}\}\quad
\text{and}\quad \Cal{Q}^n=\Cal{Q}\setminus\Cal{Q}^r=\{\alpha\in\Cal{Q}\,|\,
-\alpha\notin\Cal{Q}\}\,.}
Then 
\form{\frak{q}^n=\dsize\sum_{\alpha\in\Cal{Q}^n}{\Hat{\frak{g}}^{\alpha}}}
is the {\it nilradical} of $\frak{q}$, i.e. the set of the elements
$Z$ of its radical $\frak{r}(\frak{q})$ for which
$\roman{ad}_{\Hat{\frak{g}}}(Z)$ is nilpotent, and
\form{\frak{q}^r=\Hat{\frak{h}}\oplus\dsize\sum_{\alpha\in\Cal{Q}^r}{
\Hat{\frak{g}}^{\alpha}}} \edef\FormCE{\number\q.\number\t}
a reductive complement of $\frak{q}^n$ in $\frak{q}$. 
The complex parabolic Lie subalgebra $\frak{q}$ of $\Hat{\frak{g}}$
is its own normalizer and
the normalizer of its nilradical $\frak{q}^n$\,:
\form{\frak{q}=\{Z\in\Hat{\frak{g}}\,|\, [Z,\frak{q}]\subset\frak{q}\}\,=\,
\{Z\in\Hat{\frak{g}}\,|\, [Z,\frak{q}^n]\subset\frak{q}^n\}\,.
}
\par\smallskip
If $A\in\frak{h}_{\Bbb{R}}$, then the set
\form{\Cal{Q}_A=\{\alpha\in\Cal{R}\,|\, \alpha(A)\geq 0\}}
is parabolic, with $\Cal{Q}_A^r=\{\alpha\in\Cal{R}\,|\, \alpha(A)=0\}$
and $\Cal{Q}^n_A=\{\alpha\in\Cal{R}\,|\, \alpha(A)>0\}$.
Vice versa, if $\Cal{Q}$ is parabolic, set
$\delta=\sum\{\alpha\,|\,\alpha\in\Cal{Q}^n\}$, and define
$T_{\delta}\in\frak{h}_{\Bbb{R}}$ by
$(T_{\delta}|A)=\delta(A)$ for all $A\in\frak{h}_{\Bbb{R}}$.
Then $\Cal{Q}=\Cal{Q}_{T_{\delta}}=\{\alpha\in\Cal{R}\,|\,
(\alpha|\delta)\geq 0\}$. The set of $A\in\frak{h}_{\Bbb{R}}$
for which $\Cal{Q}=\Cal{Q}_A$ is in fact a relatively open convex cone
in $\frak{h}_{\Bbb{R}}$. \par
When $\Cal{Q}=\Cal{Q}_A$ for some $A\in\frak{h}_{\Bbb{R}}$,
we shall also write $\frak{q}_A$ for $\frak{q}_{\Cal{Q}_A}$.\par
\smallskip
The sets $\Cal{Q}^n$ associated to parabolic $\Cal{Q}$'s are called
{\it horocyclic} (see [Wa, \S{1.1}]).
The correspondence $\Cal{Q}^n\longleftrightarrow \frak{q}^n=\sum_{
\alpha\in\Cal{Q}^n}{\Hat{\frak{g}}^{\,\alpha}}$ is one-to-one between
horocyclic sets of roots in $\Cal{R}$ and nilradicals of
complex parabolic Lie subalgebras containing $\Hat{\frak{h}}$.
\par
Given a parabolic subset $\Cal{Q}\subset\Cal{R}$, 
we use the notation
$\Cal{Q}^{-n}$ for its {\it opposite horocyclic set\/}
$[{\Cal{Q}^n}]^{\roman{opp}}=\Cal{R}\setminus\Cal{Q}$\,:
the corresponding nilpotent ideal $\frak{q}^{-n}=\sum_{\alpha\in\Cal{Q}^{-n}}
{\Hat{\frak{g}}^{\alpha}}$ is a complement of $\frak{q}$ in
$\Hat{\frak{g}}$.
\smallskip
To a parabolic $\Cal{Q}\subset\Cal{R}$ 
we associate the set of Weyl chambers\,:
\form{\frak{C}(\Cal{R},\Cal{Q})=\{C\in\frak{C}(\Cal{R})\,|\,
\Cal{R}^+(C)\subset\Cal{Q}\}=\{C\in\frak{C}(\Cal{R})\,|\,
\Cal{R}^+(C)\supset\Cal{Q}^n\}\,.}
We denote by $\Cal{B}(C)$ the simple roots of $\Cal{R}^+(C)$,
for $C\in\frak{C}(\Cal{R})$. Every $\alpha\in\Cal{R}$ can be written
in a unique way as a linear combination, with integral coefficients
(either all $\geq 0$ or all $\leq 0$),
of the simple roots in $\Cal{B}(C)$:
\form{\alpha=\dsize\sum_{\beta\in\Cal{B}(C)}{k_{\alpha}^{\beta}(C)\beta}\,.}
We set 
\form{\supp{C}(\alpha)=\{\beta\in\Cal{B}(C)\,|\,
k_{\alpha}^\beta(C)\neq 0\}\,. }
If $\Cal{Q}$ is parabolic, $C\in\frak{C}(\Cal{R},\Cal{Q})$, and
$\Phi_C(\Cal{Q})=\Cal{B}(C)\cap\Cal{Q}^n$, then
\form{\Cal{Q}^n=\{\alpha\in\Cal{R}^+(C)\,|\, \supp{C}(\alpha)\cap
\Phi_C(\Cal{Q})\neq\emptyset\}\,.}
The correspondence
\form{\Cal{B}(C)\supset\Phi_C \longleftrightarrow \frak{q}=\Hat{\frak{h}}
\;\oplus\!\!
\dsize\sum_{\alpha\in\Cal{R}^+(C)}{\Hat{\frak{g}}^{\,\alpha}}\;\oplus
\!\!\dsize\sum_{\smallmatrix
\alpha\in\Cal{R}^-(C)\\
\supp{C}(\alpha)\cap\Phi_C=\emptyset\endsmallmatrix}{
\Hat{\frak{g}}^{\,\alpha}}}\edef\Fformcorr{\number\q.\number\t}
is one-to-one between subsets $\Phi_C$ of $\Cal{B}(C)$ and
complex parabolic Lie subalgebra of $\Hat{\frak{g}}$ that contain
$\Hat{\frak{h}}$ and have an associated parabolic set
$\Cal{Q}$ with $C\in\frak{C}(\Cal{R},\Cal{Q})$.\par
Having fixed a Weyl chamber $C\in\frak{C}(\Cal{R})$
and $\Phi_C\subset\Cal{B}(C)$, we shall denote by
$\frak{q}_{\Phi_C}$ the complex parabolic Lie subalgebra of
$\Hat{\frak{g}}$ defined by 
the right hand side of $(\Fformcorr)$ and by
$\Cal{Q}_{\Phi_C}$ the corresponding parabolic set.\par\smallskip
\smallskip
We denote by $\alweyl$ the Weyl group of $\Cal{R}$, (i.e. the group
of isometries of $\frak{h}^*_{\Bbb{R}}$ generated by the symmetries
$\xi@>>>s_{\alpha}(\xi)=\xi-2[(\xi|\alpha)/\|\alpha\|^2]\,\alpha$ for
$\alpha\in\Cal{R}$)
and by $\alaff$ the group of all isometries of $\frak{h}^*_{\Bbb{R}}$
(with respect to the scalar product defined above) that transform
$\Cal{R}$ into itself. For $C\in\frak{C}(\Cal{R})$ we denote by
$\alaffine{C}$ the subgroup of $w\in\alaff$ for which
$w(\Cal{R}^+(C))=\Cal{R}^+(C)$. Then $\alaff=\alaffine{C}\rtimes
\alweyl$. \par
We define $\alweyle{\Cal{Q}}$
and $\alaffe{\Cal{Q}}$ as the subgroups of
$\alweyl$ and $\alaff$, respectively, that transform $\Cal{Q}$
into itself. Then we have Chevalley's Lemma 
(see e.g. [Wa, Theorem 1.1.2.8])\,:
\lem{The group $\alweyle{\Cal{Q}}$ is generated by the symmetries
$s_{\alpha}$ with $\alpha\in\Cal{Q}^r$. If $C\in\frak{C}(\Cal{R},\Cal{Q})$,
then the symmetries $s_{\alpha}$ with 
$\alpha\in\Cal{B}(C)\setminus\Phi_C(\Cal{Q})$ generate 
$\alweyle{\Cal{Q}}$ and $\alaffe{\Cal{Q}}$ is a semidirect product
$\alaffe{\Cal{Q}}=\alaffineq{C}{Q}\rtimes \alweyle{\Cal{Q}}$,
with $\alaffineq{C}{Q}=\alaffine{C}\cap\alaffe{\Cal{Q}}$.
\qed}
\par\smallskip
Let $\Hat{\bold{G}}$ be a connected complex Lie group with Lie algebra
$\Hat{\frak{g}}$. If $\bold{Q}$ is any Lie subgroup of
$\Hat{\bold{G}}$ with complex Lie subalgebra $\frak{q}$ that is parabolic
in $\Hat{\frak{g}}$, then $\bold{Q}$ is closed, connected and 
coincides with its normalizer
in $\Hat{\bold{G}}$ and is the normalizer of
its Lie algebra for the adjoint representation\,:
\form{\bold{Q}=\{g\in\Hat{\bold{G}}\,|\, 
g\,\bold{Q}\,g^{-1}
=\bold{Q}\}=\{g\in\Hat{\bold{G}}\,|\, \roman{Ad}_{\Hat{\frak{g}}}(g)
(\frak{q})=\frak{q}\}}
The homogeneous space $\frak{M}=\Hat{\bold{G}}/\bold{Q}$ is compact and
simply connected. Since the center $\bold{Z}(\Hat{\bold{G}})$ 
of $\Hat{\bold{G}}$ is
contained in all its parabolic subgroups, the choice of
different connected complex 
Lie groups $\Hat{\bold{G}}$
yields the same $\frak{M}$. 
Hence
we can consider the
complex flag manifold $\frak{M}=\Hat{M}(\Hat{\frak{g}},\frak{q})$
as an object 
associated simply to the pair $(\Hat{\frak{g}},\frak{q})$.
We also recall 
(see [W, 2.7]) that
the integral cohomology $H^*(\frak{M},\Bbb{Z})$
is torsion free and $0$ in odd degrees. \par
If $\Hat{\frak{h}}$ is a Cartan subalgebra of $\Hat{\frak{g}}$
contained in $\frak{q}$ and $\Cal{Q}\subset\Cal{R}(\Hat{\frak{g}},
\Hat{\frak{h}})$ the parabolic set of roots associated to $\frak{q}$,
the complex dimension of 
$\frak{M}=\Hat{M}(\Hat{\frak{g}},\frak{q})$ 
equals the number of roots in $\Cal{Q}^n$.\par

\edef\paragC{\number\q}

\se{Parabolic $CR$ algebras and parabolic $CR$ manifolds}

A $CR$ algebra $(\frak g,\frak q)$ is called {\it parabolic}
if $\frak g$ is finite dimensional, and $\frak q$
is a {\it parabolic} subalgebra of its complexification $\hat\frak g$.\par
By the results quoted at the end of \S\paragC, all the 
homogeneous spaces
$\tilde{M}=\tilde{M}(\frak{g},\frak{q})$, 
$M={M}(\frak{g},\frak{q})$, and
$\frak{M}=\Hat{M}(\Hat{\frak{g}},\frak{q})$ 
of {\S\paragB}
are well defined and
$M$ is an orbit, in the complex flag manifold
$\frak{M}$ of a connected complex semisimple Lie group $\Hat{\bold{G}}$,
of a connected real form $\bold{G}$ of 
$\Hat{\bold{G}}$. We say that $M$ is a {\it parabolic $CR$ manifold}.
\par
Vice versa, if $\bold{G}$ is a connected
real form of the complex semisimple
Lie group $\Hat{\bold{G}}$, then all $\bold{G}$-orbits in the complex
flag manifolds $\frak{M}=\Hat{M}(\Hat{\frak{g}},\frak{q})$ are homogeneous
$CR$ manifolds of the form $M=M(\frak{g},\frak{q}')$,
for some parabolic complex Lie subalgebra $\frak{q}'$
of $\Hat{\frak{g}}$, conjugated to $\frak{q}$ by
an inner automorphism.
\par 
It is worth noticing that, in the definition of the homogeneous
$CR$ manifold
$M(\frak{g},\frak{q})=\bold{G}/\bold{G}_+$,
we can define the isotropy $\bold{G}_+=\bold{G}_+(\frak{g},\frak{q})$
by
\form{\bold{G}_+=\{g\in\bold{G}\,|\,\roman{Ad}_{\Hat{\frak{g}}}(
g)(\frak{q})=\frak{q}\}\,.}
Since the center of $\bold{G}$ is always contained in $\bold{G}_+$,
we obtain an equivalent definition of $M(\frak{g},\frak{q})$
if we substitute to $\Hat{\bold{G}}$ any connected complex Lie group
$\Hat{\bold{G}}'$
with the same Lie algebra $\Hat{\frak{g}}$
and to $\bold{G}$ the  
analytic subgroup $\bold{G}'$ of $\Hat{\bold{G}}'$
with Lie algebra $\frak{g}$.
However, it is more convenient to fix a simply connected $\Hat{\bold{G}}$,
since in this case, by [BT2, Corollaire 4.7], we have\,:
\form{\bold{G}=\Hat{\bold{G}}^{\sigma}=\{g\in\Hat{\bold{G}}\,|\,
\sigma(g)=g\}\,,} \edef\FRMLDB{\number\q.\number\t}
where $\sigma:\Hat{\bold{G}}@>>>\Hat{\bold{G}}$ is the 
anti-holomorphic involution
of $\Hat{\bold{G}}$ corresponding to the conjugation 
$\sigma$ of $\Hat{\frak{g}}$
defined by the real form $\frak{g}$. 
\par\smallskip
We begin by recalling some general facts about parabolic $CR$ algebras,
and their associated $CR$ manifolds, that were explained in
[AMN, \S{5}].
\prop{A parabolic $CR$ algebra $(\frak g,\frak q)$ is effective if
and only if\,:
$(i)$ $\frak g$  is semisimple,
$(ii)$ no simple ideal of $\hat\frak g$ is 
contained in $\frak q\cap\overline{\frak q}$.\par
An effective parabolic $CR$ algebra $(\frak{g},\frak{q})$ with $\frak{g}$ 
simple is either totally complex or ideal nondegenerate.\qed}
\prop{Let $(\frak g,\frak q)$ be an effective parabolic $CR$ algebra
and let $\frak g =\frak g_1\oplus \cdots \oplus\frak g_\ell$ be the
decomposition of $\frak g$ into the direct sum of its simple ideals.
Then\,:
\roster
\item"($i$)" $\frak q=\frak q_1\oplus \cdots \oplus\frak q_\ell$
where $\frak q_j=\frak q\cap\hat\frak g_j$ for $j=1,\hdots,\ell$;
\item"($ii$)" for each $j=1,\hdots,\ell$, $(\frak g_j,\frak q_j)$
is an effective parabolic $CR$ algebra;
\item"($iii$)" $(\frak g,\frak q)$ is ideal (resp. weakly, strictly)
nondegenerate if and only if for each $j=1,\hdots,\ell$, the
$CR$ algebra $ (\frak g_j,\frak q_j)$ is  ideal (resp. weakly, strictly)
nondegenerate;
\item"($iv $)" $(\frak g,\frak q)$ is fundamental
if and only if for each $j=1,\hdots,\ell$, the
$CR$ algebra $ (\frak g_j,\frak q_j)$ is fundamental.
\item"($v$)" We have ($\cong$ meaning biholomorphic
or $CR$ equivalence)\,:
$$\align
\Hat{M}(\Hat{\frak{g}},\frak{q})&\cong\Hat{M}(\Hat{\frak{g}}_1,\frak{q}_1)
\times\cdots\times \Hat{M}(\Hat{\frak{g}}_{\ell},\frak{q}_{\ell})\,,\\
\Tilde{M}({\frak{g}},\frak{q})&\cong\Tilde{M}({\frak{g}_1},\frak{q}_1)
\times\cdots\times \Tilde{M}({\frak{g}_{\ell}},\frak{q}_{\ell})\,,\\
{M}({\frak{g}},\frak{q})&\cong{M}({\frak{g}_1},\frak{q}_1)
\times\cdots\times {M}({\frak{g}_{\ell}},\frak{q}_{\ell})\,.\qquad\qed
\endalign
$$
\endroster}
When $\frak q$ is parabolic in $\hat\frak g$,
its conjugate $\bar{\frak q}$
with respect to the real form $\frak g$ is also parabolic in
$\hat\frak g$. Therefore the intersection $\frak q\cap\bar{\frak q}$
contains a Cartan subalgebra $\hat{\frak h}$ that is invariant under
conjugation. The intersection $\frak h=\hat{\frak h}\cap\frak g$ is
a Cartan subalgebra of $\frak g$, contained in 
$\frak{g}_+=\frak{q}\cap\frak{g}$. \par
A Cartan subalgebra $\frak{h}$ of $\frak{g}$ contained in
$\frak{g}_+=\frak{q}\cap\frak{g}$ is said to be {\it adapted} to
$(\frak{g},\frak{q})$.
We also recall [AMN, Proposition 5.4]\,:
\prop{Let $(\frak g,\frak q)$ be an effective parabolic $CR$ algebra,
with isotropy subalgebra $\frak{g}_+=\frak{q}\cap\frak{g}$.
The elements $A$ of the radical
$\frak r(\frak g_+)$ of $\frak{g}_+$ for which
$\roman{ad}_{\frak g}(A):\frak g@>>>\frak g$ is nilpotent,
form a nilpotent ideal $\frak{n}$ of $\frak g_+$.  It admits
a reductive complement $\frak{g}_0$ 
in $\frak g_+$\,:
\form{\frak{g}_+=\frak{n}\oplus\frak{g}_0\, .}
The reductive subalgebra $\frak{g}_0$ is uniquely determined modulo
inner automorphisms of $\frak{g}_+$ from 
the subgroup generated 
by those of the form $\exp\left(\roman{ad}_{\frak g_+}(X)\right)$
with $X\in\frak{n}$.\qed}\edef\proposizioneAE{\number\q.\number\x}
\edef\formulaAA{\number\q.\number\t}\edef\FormDA{\number\q.\number\t}
\medskip
Let $\frak z_0$ be the center and
$\frak s_0=[\frak{g}_0,\frak{g}_0]$ the semisimple ideal
of $\frak{g}_0$. Then
\form{\frak{g}_0\, = \, \frak z_0\oplus\frak s_0\, .}
Thus, a Cartan subalgebra $\frak h\subset\frak g_+$ of $\frak g$ 
can be taken as the direct sum 
\form{\frak h \, = \, \frak z_0\oplus\frak h_0}
\noindent
of the center\edef\cartdec{\number\q.\number\t} 
$\frak z_0$ of $\frak{g}_0$ and a Cartan subalgebra $\frak h_0$
of $\frak s_0$. Vice versa, 
every Cartan subalgebra $\frak h$ of $\frak g$
adapted to $(\frak{g},\frak{q})$
has the form (\cartdec) for some reductive
subalgebra $\frak{g}_0$ of $\frak g_+$.\par\smallskip
 
It is also convenient to consider 
a Cartan decomposition (see e.g. [B2])\,:
\form{\frak g=\frak{k}\oplus\frak{p}} 
of $\frak{g}$, \edef\decompocartan{\number\q.\number\t}
corresponding to a Cartan
involution $\vartheta$ . 
The set $\frak{k}=\{X\in\frak{g}\,|\,
\vartheta(X)=X\}$ of fixed points of $\vartheta$ is
a maximal compact Lie subalgebra of $\frak{g}$ and 
$\frak{p}=\{X\in\frak{g}\,|\,
\vartheta(X)=-X\}$ 
its orthogonal for the Killing form
$\kappa_{\frak{g}}$ of $\frak{g}$. 
Any $\vartheta$-invariant Cartan subalgebra $\frak{h}$
of $\frak g$ decomposes into the direct sum
$\frak{h}=\frak{h}^+\oplus\frak{h}^-$ of 
its {\it compact} (or {\it toroidal}) part 
$\frak{h}^+=\frak{h}\cap\frak{k}\subset\frak{k}$ 
and its noncompact (or {\it vector part})
$\frak{h}^-=\frak{h}\cap\frak{p}\subset\frak{p}$.\par
We say that the Cartan decomposition (\decompocartan)
is {\it adapted} to 
the effective parabolic  $CR$ algebra $(\frak{g},\frak{q})$ if
$\frak{k}$ contains a maximal compact Lie subalgebra
of $\frak{g}_+$. Then\,:
\lem{If (\decompocartan) is adapted to the parabolic $CR$ algebra
$(\frak{g},\frak{q})$, then 
every Cartan subalgebra $\frak{h}$ of $\frak{g}$ that is adapted
to $(\frak{g},\frak{q})$ is conjugate, modulo an inner automorphism of
$\frak{g}_+$, to a $\vartheta$-invariant Cartan subalgebra
$\frak{h}_0$ of $\frak{g}$ that is adapted to $(\frak{g},\frak{q})$.
\par
Vice versa, if $\frak{h}$ is a Cartan subalgebra of $\frak{g}$
adapted to $(\frak{g},\frak{q})$, then there exists a Cartan decomposition
(\decompocartan), adapted to $(\frak{g},\frak{q})$,
 such that $\frak{h}=(\frak{h}\cap\frak{k})\oplus
(\frak{h}\cap\frak{p})$. \par
In particular, if $\{\frak{q}_i\,|\, i\in I\}$ 
is a family of complex parabolic Lie subalgebras
of $\Hat{\frak{g}}$ such that $\bigcap_{i\in I}{\frak{q}_i}$ is
parabolic in $\Hat{\frak{g}}$, then there exist both a Cartan
decomposition (\decompocartan) and a Cartan subalgebra $\frak{h}$
of $\frak{g}$, compatible with (\decompocartan), that are adapted
to all the $(\frak{g},\frak{q}_i)$'s.\qed}
We say that $(\vartheta,\frak{h})$ is an {\it adapted Cartan pair} for
$(\frak{g},\frak{q})$ if\,: 
\par\noindent
$\;(i)$ $\vartheta$ is the Cartan involution
of a Cartan decomposition (\decompocartan) adapted to $(\frak{g},\frak{q})$\,;
\par\noindent
$(ii)$ $\frak{h}$ is a $\vartheta$-invariant Cartan subalgebra of
$\frak{g}$ contained in $\frak{g}_+=\frak{g}\cap\frak{q}$.
\par\smallskip
Being
$\sigma:\Hat{\frak{g}} \ni X @>>> \bar X \in\Hat{\frak{g}}$ 
the conjugation in $\Hat{\frak{g}}$ associated to the real form $\frak g$,
and having fixed (\decompocartan),
we also consider the conjugation 
$\tau:\Hat{\frak{g}}@>>>\Hat{\frak{g}}$
of $\Hat{\frak{g}}$ with respect to
its compact real form $\frak{u}=\frak{k}\oplus i\,\frak{p}$
and use the same symbol
$\vartheta$ to denote the $\Bbb{C}$-linear extension to
$\Hat{\frak{g}}$ of the
Cartan involution $\vartheta$ of $\frak{g}$. 
We obtain in this way three commuting involutions
$\sigma,\,\tau,\,\vartheta$ of $\Hat{\frak{g}}$,
each being the composition product of the other two\,:
\form{\tau=\vartheta\circ\sigma=\sigma\circ\vartheta\,,\quad
\sigma=\vartheta\circ\tau=\tau\circ\vartheta\,,\quad
\vartheta=\sigma\circ\tau=\tau\circ\sigma\,.}
In particular $\frak u$ is invariant under $\sigma$:
$\sigma(\frak u)=\frak u$.\par\smallskip
Let $(\vartheta,\frak{h})$ be a Cartan pair adapted to $(\frak{g},\frak{q})$.
 Then $\frak{h}_{\Bbb{R}}$ decomposes as\,:
$\frak{h}_{\Bbb{R}}=\frak{h}^- \oplus i\,\frak{h}^+$ 
and $\frak{h}_{\Bbb{R}}$ is a
Cartan subalgebra 
of a split real form $\frak{g}_{\Bbb{R}}$ of
$\Hat{\frak{g}}$. 
The involutions
$\sigma$, $\tau$  and $\vartheta$ transform $\frak{h}_{\Bbb{R}}$ into
itself. Hence, by transposition, they define involutions on
$\frak{h}_{\Bbb{R}}^*$, that we still denote by the same symbols
$\sigma$, $\tau$  and $\vartheta$, and that transform 
the set of roots
$\Cal{R}$ into itself. We set $\bar\alpha=\sigma(\alpha)$
for all $\alpha\in\frak{h}_{\Bbb{R}}^*$. We have\,: 
\form{\tau(\alpha)=-\alpha\,,\quad \vartheta(\alpha)=-\bar\alpha
\quad\forall \alpha\in\frak{h}_{\Bbb{R}}^*\,.}
The reductive complement $\frak{q}^r$ of
$\frak{q}^n$ in $\frak{q}$ of $(\FormCE)$ is $\frak{q}\cap\tau(\frak{q})$,
while the reductive complement 
$\frak{g}_0$ of $\frak{n}$ in $\frak{g}_+$ of
$(\FormDA)$ can be taken equal to $\frak{g}_+\cap\vartheta(\frak{g}_+)$. 
\medskip
A real Lie subalgebra 
$\frak t$ of $\frak g$
is {\it triangular} if all linear maps 
$\roman{ad}_{\frak g}(X)\in\frak{gl}_{\Bbb R}(\frak g)$ with $X\in\frak t$
can be simultaneously represented by triangular matrices in a suitable
basis of $\frak g$. All maximal triangular subalgebras of $\frak g$ are
conjugate by an inner automorphism (cf. [Mo, \S5.4]). 
A real Lie subalgebra of $\frak g$ containing a maximal triangular
subalgebra of $\frak g$ is called a $t$-{\it subalgebra}.\par
An effective parabolic $CR$ algebra $(\frak g,\frak q)$ will be
called {\it minimal} if $\frak g_+=\frak q\cap\frak g$ is
a $t$-subalgebra of $\frak g$. \par
We observe that a maximal triangular subalgebra of $\frak g$ contains
a maximal Abelian subalgebra of semisimple elements 
of $\frak{g}$ having real eigenvalues.
We have (see [AMN, \S{5}]):
\prop{An effective parabolic 
minimal $CR$ algebra $(\frak{g},\frak{q})$
admits an adapted Cartan pair $(\vartheta,\frak{h})$ in which
$\frak{h}$ is a
maximally noncompact Cartan subalgebra of $\frak{g}$.\par
Let $\frak g$ be a semisimple real Lie algebra and 
$\frak q$ a parabolic
subalgebra of its complexification $\Hat{\frak g}$. Then, up to
$CR$ isomorphisms, there is a unique parabolic 
minimal $CR$ algebra
$(\frak g',\frak q')$ with $\frak g'$ isomorphic to $\frak g$ and
$\frak q'$ isomorphic to $\frak q$.\par
A totally real effective parabolic $CR$ algebra is minimal.\qed}

By [AMN, Theorem 8.6 and Corollary 9.2] we have\,:
\thm{Let $(\frak{g},\frak{q})$ be a parabolic minimal 
fundamental $CR$ algebra. Then 
$\tilde{M}(\frak{g},\frak{q})\cong M(\frak{g},\frak{q})$ is
compact and simply connected.\qed}
Next we consider the $CR$ fibrations of Proposition {\PropBD}
in the special case of a parabolic $CR$ algebra.
\thm{Every effective parabolic $CR$ algebra $(\frak g,\frak q)$
admits a unique $\frak g$-equivariant $CR$ fibration 
$(\frak{g},\frak{q})@>>>(\frak{g},\Hat{\frak{g}}')$,
where $\Hat{\frak{g}}'\supset\frak{q}$ is the complexification of a real
parabolic subalgebra $\frak{g}'$ of $\frak{g}$,
and the fiber is fundamental. 
The basis $(\frak{g},\Hat{\frak{g}}')$ is a totally
real parabolic $CR$ algebra and also the 
fiber $(\frak{g}',\frak{q})$ is parabolic.
\par 
This yields a $\bold{G}$-equivariant
$CR$ fibration
$\pi:M(\frak{g},\frak{q})@>>>M(\frak{g},\Hat{\frak{g}}')$
with compact basis.
Each connected component of the fiber is $CR$ diffeomorphic to
$M(\frak{g}',\frak{q})$.}
\dimo Let $(\frak g,\frak q)$ be an effective parabolic $CR$
algebra. 
The complex subalgebra $\frak{q}''$ generated by 
$\frak q+\bar{\frak{q}}$ is parabolic in $\Hat{\frak{g}}$
because contains $\frak{q}$, and 
is the complexification
of a real parabolic subalgebra
$\frak{g}'$ of $\frak{g}$ because
$\bar{\frak{q}}\,''=\frak{q}''$. Then $(\formlacaa)$ yields
a $\frak g$-equivariant $CR$ fibration with a totally real basis.
The fiber is $(\frak g',\frak q)$. This is parabolic because
$\frak{q}$, being parabolic in $\Hat{\frak{g}}$, is also parabolic
in $\Hat{\frak{g}}'\subset\Hat{\frak{g}}$. \par
The final statement follows from the commutative diagram\,:
$$\CD
M(\frak{g},\frak{q}) @>>> \Hat{M}(\Hat{\frak{g}},\frak{q})\\
@V{\pi}VV @V{\Hat{\pi}}VV\\
M(\frak{g},\Hat{\frak{g}}') @>>> \Hat{M}(\Hat{\frak{g}},\Hat{\frak{g}}')
\endCD
$$
that yields an embedding of each fiber of
$\pi:M(\frak{g},\frak{q})@>>>M(\frak{g},\Hat{\frak{g}}')$
into a fiber of $\Hat{\pi}:\Hat{M}(\Hat{\frak{g}},\frak{q})
@>>>\Hat{M}(\Hat{\frak{g}},\Hat{\frak{g}}')$.
The basis $M(\frak{g},\Hat{\frak{g}}\,')$ 
is compact because $(\frak{g},\Hat{\frak{g}}\,')$
is totally real and hence minimal.\enddimo
\medskip
\esemp{The fiber is not necessarily connected.
An example of
a $\bold{G}$-equivari\-ant fibration
$M(\frak{g},\frak{q})@>{\pi}>>M(\frak{g},\frak{q}')$ of parabolic
$CR$ manifolds with a disconnected fiber is the following.
\par
Let $\frak{g}=\frak{sl}(3,\Bbb{R})\subset\Hat{\frak{g}}
=\frak{sl}(3,\Bbb{C})$, and
$$\align
\frak{q}&=\left\{\left(\smallmatrix
z_{1,1}&z_{1,2}&z_{1,3}\\
0&z_{2,2}&z_{2,3}\\
0&0&z_{3,3}\endsmallmatrix\right)\in\frak{sl}(3,\Bbb{C})\right\}\, ,\\
\frak{q}'&=\left\{\left(\smallmatrix
z_{1,1}&z_{1,2}&z_{1,3}\\
0&z_{2,2}&z_{2,3}\\
0&z_{3,2}&z_{3,3}\endsmallmatrix\right)\in\frak{sl}(3,\Bbb{C})\right\}\,.
\endalign
$$
Let $e_1,e_2,e_3$ be the canonical basis of $\Bbb{C}^3$.
We have
$\bold{G}=\bold{SL}(3,\Bbb{R})\subset\Hat{\bold{G}}=\bold{SL}(3,\Bbb{C})$.
The isotropy subgroup
$\bold{G}_+=\{g\in\bold{SL}(3,\Bbb{R})\,|\,
\roman{Ad}_{\Hat{\frak{g}}}(g)(\frak{q})=\frak{q}\}$
of the point $\o=e\,\bold{Q}\in M(\Hat{\frak{g}},\frak{q})$ is 
the stabilizer in $\bold{SL}(3,\Bbb{R})$ 
of the complete flag
$\langle e_1\rangle\subset\langle e_1,e_2\rangle$. It has $4$ connected
components. Likewise, the isotropy subgroup
$\bold{G}'=\{g\in\bold{SL}(3,\Bbb{R})\,|\,
\roman{ad}_{\Hat{\frak{g}}}(g)(\frak{q}')=\frak{q}'\}$, 
at $\o'=e\,\bold{Q}'\in M(\frak{g},\frak{q}')$,
is the stabilizer in $\bold{SL}(3,\Bbb{R})$ 
of the line $\langle e_1\rangle$, and has
$2$ connected components. The fiber $\bold{G}'/\bold{G}_+$
is isomorphic to the stabilizer of $\langle e_1\rangle$ in
$\bold{SL}_{\Bbb{R}}(\Bbb{R}\,e_1+\Bbb{R}\,e_2)$
and has
$2$ connected components.}\par\medskip
\thm{Let $(\frak g,\frak q)$ be an effective parabolic $CR$ algebra.
Then there is a unique $\frak{g}$-equivariant $CR$ fibration
$(\frak{g},\frak{q})@>>>(\frak{g},\frak{q}')$ with a weakly nondegenerate
basis $(\frak{g},\frak{q}')$ and a totally complex fiber.
The basis
$(\frak g,\frak q')$ is a parabolic $CR$ algebra.} 
\dimo We recall from [MN4, \S5], that $\frak q'$ is the unique
maximal subalgebra of $\Hat{\frak g}$ that contains $\frak q$
and is contained in $\frak q+\overline{\frak q}$. Clearly
$\frak q'$ is parabolic because it contains the parabolic
subalgebra $\frak q$. \enddimo
\smallskip
Next we investigate the general structure of the fiber of
a $\bold{G}$-equivariant submersion $M(\frak{g},\frak{q})@>>>
M(\frak{g},\frak{q}')$ for a pair of complex parabolic 
subalgebras $\frak{q}\subset\frak{q}'$ of $\Hat{\frak{g}}$.
\thm{Let $\frak{q}\subset\frak{q}'$ be complex parabolic Lie subalgebras
of $\Hat{\frak{g}}$. With
$\frak{g}'=\frak{g}\cap\frak{q}'$, the $CR$ algebra
$(\frak{g}',\bar{\frak{q}}'\cap\frak{q})$ is the fiber
of the $\frak{g}$-equivariant fibration
$(\frak{g},\frak{q})@>>>(\frak{g},\frak{q}')$ (see \S{\paragB}).\par
The nilradical $\frak{n}'$ of $\frak{g}'$, consisting of
the $\roman{ad}_{\frak{g}}$-nilpotent elements of
the radical $\frak{r}(\frak{g}')$ of $\frak{g}'$, has a reductive
complement $\frak{g}'_0$ in $\frak{g}'$ such that\,:
\roster
\item"$(i)$" The $CR$ algebra
$(\frak{g}'_0,\Hat{\frak{g}}'_0\cap\frak{q})$ is parabolic.
\item"$(ii)$" 
The fiber $(\frak{g}',\Hat{\frak{g}}'\cap\frak{q})$ is the semidirect
sum of the parabolic $CR$ algebra $(\frak{g}'_0,\Hat{\frak{g}}'_0\cap\frak{q})$
and of the nilpotent $CR$ algebra $(\frak{n}',\Hat{\frak{n}}'\cap\frak{q})$.
\item"$(iii)$" The nilpotent $CR$ algebra 
$(\frak{n}',\Hat{\frak{n}}'\cap\frak{q})$ is totally complex.
\item"$(iv)$" The connected components of the fibers of the
$\bold{G}$-equivariant fibration 
$\pi:M(\frak{g},\frak{q})@>>>M(\frak{g},\frak{q}')$ 
are $CR$ diffeomorphic to the Cartesian product of a parabolic
$CR$ manifold
$M(\frak{g}'_0,\frak{q}\cap\Hat{\frak{g}}'_0)$ 
and of a Euclidean
complex manifold ($\cong \Bbb{C}^{\ell}$).
\endroster}\edef\TeoremaID{\number\q.\number\x}
\dimo Fix a Cartan pair $(\vartheta,\frak{h})$, that is adapted for both
$(\frak{g},\frak{q})$ and $(\frak{g},\frak{q}')$. With the notation
introduced above, set
$\frak{q}^r=\frak{q}\cap\tau(\frak{q})$ and
$\frak{q}'\,^r=\frak{q}'\cap\tau(\frak{q}')$.
Since $\frak{q}\subset\frak{q}'$  and 
$\frak{q}$ is parabolic, we have the inclusions\,:
$\frak{q}^r\subset{\frak{q}'}\,^r$ and
$\frak{q}^n\supset{\frak{q}'}\,^n$\,. 
The complexification of the fiber $\frak{g}'$ is\,:
$$\Hat{\frak{g}}'=\frak{q}'\cap\bar{\frak{q}}\,'=
\Hat{\frak{g}}'_0\rtimes\Hat{\frak{n}}'\,,$$
where\,:
$$\cases
\Hat{\frak{g}}'_0={\frak{q}'}\,^r\cap{\bar{\frak{q}}'}\,^r,\\
\Hat{\frak{n}}'=\left({\frak{q}'}\,^r\cap{\bar{\frak{q}}'}\,^n\right)\oplus
\left({\frak{q}'}\,^n\cap{\bar{\frak{q}}'}\,^r\right)\oplus
\left({\frak{q}'}\,^n\cap{\bar{\frak{q}}'}\,^n\right)=
\left(\frak{q}'\cap{\bar{\frak{q}}'}\,^n\right)\,+\,
\left({\frak{q}'}\,^n\cap{\bar{\frak{q}}'}\right)\,.
\endcases
$$
Thus $\frak{g}'=\frak{g}'_0\rtimes\frak{n}'$, where
$\frak{n}'=\Hat{\frak{n}}'\cap\frak{g}$ is a real form of the
nilradical $\Hat{\frak{n}}'$ of $\frak{q}'\cap\bar\frak{q}'$  
and $\frak{g}'_0:=\Hat{\frak{g}}'_0\cap\frak{g}$ 
a reductive complement
of $\frak{n}'$ in $\frak{g}'$. \par
We have\,: 
$$\Hat{\frak{g}}'\cap\frak{q}=
\left(\Hat{\frak{g}}'_0\rtimes\Hat{\frak{n}}'\right)
\cap\frak{q}\,=\,\left(\Hat{\frak{g}}'_0\cap\frak{q}\right)\rtimes
\left(\Hat{\frak{n}}'\cap\frak{q}\right)\,,$$
so that\,: 
$$(\frak{g}',\Hat{\frak{g}}'\cap\frak{q})=
(\frak{g}'_0,\Hat{\frak{g}}'_0\cap\frak{q})
\rtimes (\frak{n}',\Hat{\frak{n}}'\cap\frak{q})\,.$$ 
The complex
Lie subalgebra $\Hat{\frak{g}}'_0\cap\frak{q}$ is parabolic
in $\Hat{\frak{g}}'_0$, because $\Hat{\frak{g}}'_0$ is reductive,
$\frak{q}$ is parabolic in $\Hat{\frak{g}}$, 
and 
$\Hat{\frak{g}}'_0\cap\frak{q}$ 
contains a Cartan subalgebra of $\Hat{\frak{g}}'_0$
and  
of $\Hat{\frak{g}}$.\par
Note that $\Hat{\frak{n}}'\cap\frak{q}$ is contained in, but in general not
equal to, the nilradical $\Hat{\frak{n}}$ of $\frak{q}\cap\bar{\frak{q}}$.
We have\,:
$$\Hat{\frak{n}}'\cap\frak{q}
\supset {\frak{q}'}\,^n\cap\bar{\frak{q}}\,'\,,$$
so that\,:
$$\Hat{\frak{n}}'\cap\frak{q}\,+\,\overline{\Hat{\frak{n}}'\cap\frak{q}}
\supset {\frak{q}'}\,^n\cap\bar{\frak{q}}\,'\,+\,
\frak{q}'\cap{\bar{\frak{q}}'}\,^n=\Hat{\frak{n}}'$$
shows that actually \,:
$$\Hat{\frak{n}}'\cap\frak{q}\,+\,\overline{\Hat{\frak{n}}'\cap\frak{q}}
=\Hat{\frak{n}}'$$
and the nilpotent $CR$ algebra $(\frak{n}',\Hat{\frak{n}}'\cap\frak{q})$
is totally complex. \par\smallskip
The $\bold{G}$-equivariant
fibration $\pi:M(\frak{g},\frak{q})@>>>M(\frak{g},\frak{q}')$ is the
restriction of the $\Hat{\bold{G}}$ equivariant fibration
$\Hat{\pi}:\Hat{M}(\Hat{\frak{g}},\frak{q})@>>>
\Hat{M}(\Hat{\frak{g}},\frak{q}')$.
The typical fiber $\frak{F}$ of $\Hat{\pi}$ is $\bold{Q}'/\bold{Q}$. Since
$\frak{q}$ is a parabolic complex Lie subalgebra of
$\frak{q}'$,  the fiber $\frak{F}$ is 
a complex flag manifold and, in particular, is 
compact, connected and
simply connected. 
Thus the typical fiber $F$ of 
$\pi:M(\frak{g},\frak{q})@>>>M(\frak{g},\frak{q}')$
is a submanifold of a complex flag manifold.\par
Denote still 
by ${\tau}:\Hat{\bold{G}}@>>>\Hat{\bold{G}}$ the involution
of $\Hat{\bold{G}}$ associated to the conjugation 
$\tau:\Hat{\frak{g}}@>>>\Hat{\frak{g}}$ with respect to the compact
real form $\frak{u}=\frak{k}\oplus i\,\frak{p}$. 
Since ${\frak{q}'}^n\subset\frak{q}^n$, the fiber
$\frak{F}$ can be viewed also as
a flag manifold of the reductive
complex closed connected Lie subgroup
${\bold{Q}'}\,^r={\bold{Q}'}\cap\tau({\bold{Q}'})$ of $\Hat{\bold{G}}$.
The fiber $F$ is contained in the orbit $\frak{F}_0\subset\frak{F}$ 
of the closed complex Lie
subgroup $\Hat{\bold{G}}':=\bold{Q}'\cap{\sigma}(\bold{Q}')$
of ${\bold{Q}'}$. The group $\Hat{\bold{G}}'$ is connected because
it contains a Cartan subgroup of $\Hat{\bold{G}}$,
and decomposes into the semidirect product
\form{\Hat{\bold{G}}'=\Hat{\bold{G}}_0'\rtimes\Hat{\bold{N}}'\,,}
where $\Hat{\bold{G}}_0'$ and $\Hat{\bold{N}}'$ are the analytic 
complex Lie subgroups
of $\Hat{\bold{G}}$ generated by the Lie subalgebras $\Hat{\frak{g}}'_0$
and $\Hat{\frak{n}}'$, respectively. We have 
$\Hat{\bold{G}}_0'={\bold{Q}'}\,^r\cap\sigma({\bold{Q}'}\,^r)$, so that
$\Hat{\bold{G}}_0'$ is closed in $\Hat{\bold{G}}$. Moreover, since
$\roman{ad}_{\Hat{\frak{g}}}(Z)$ is nilpotent for all $Z\in\Hat{\frak{n}}$,
by Engel's  theorem and the semisimplicity of $\Hat{\frak{g}}$, we obtain
that $\exp:\Hat{\frak{n}}'@>>>\Hat{\bold{N}}'$ is an analytic diffeomorphism,
and $\Hat{\bold{N}}'$ is Euclidean.\par
The validity of $(iv)$ is then a consequence of the
next Proposition.\enddimo
\prop{Let $\Hat{\bold{N}}'$ be a connected
nilpotent complex Lie group with
complex Lie algebra $\Hat{\frak{n}}'$ and
$\frak{n}'$ a real form
of $\Hat{\frak{n}}'$.
Let $\bold{N}'$ be the real analytic 
Lie subgroup of $\Hat{\bold{N}}'$ with Lie algebra
$\frak{n}'$, and $\bold{Q}_0$ a closed connected complex Lie subgroup
of $\Hat{\bold{N}}'$, with Lie algebra
$\frak{q}_0\subset\Hat{\frak{n}}'$,
and set $\bold{N}=\bold{Q}_0\cap\bold{N}'$. Assume that the $CR$ algebra
$(\frak{n}',\frak{q}_0)$ is totally complex. 
With $E=\bold{N}'/\bold{N}$
and $\Hat{E}=\Hat{\bold{N}}'/\bold{Q}_0$, the natural map
$E @>>>\Hat{E}$ obtained from the inclusion
$\bold{N}'\hookrightarrow \Hat{\bold{N}}'$ is a diffeomorphism.}

\dimo
The condition that $(\frak{n}',\frak{q}_0)$ is totally complex 
is equivalent to the equality
$\frak{n}'+\frak{q}_0=\hat{\frak{n}}'$. 
Since $\Hat{\frak{n}}'$ is nilpotent,
this equality implies (see the proof below)
that the map 
$\bold{N}'\times\bold{Q}_0\ni (n,q) @>>> n\cdot q\in\Hat{\bold{N}}'$ is 
onto,
and hence the inclusion $\bold{N}'\hookrightarrow\Hat{\bold{N}}'$ yields,
by passing to the quotients, a smooth one-to-one 
map $f:E@>>>\Hat{E}$. We note that
$E=\bold{N}'/\bold{N}$ is 
a complex manifold with the homogeneous
$CR$ structure defined by the $CR$ algebra $(\frak{n}',\frak{q}_0)$.
With this complex structure on $E$, and with the complex structure
that $\Hat{E}$ 
inherits from $\Hat{\bold{N}}'$, the map
$f$ is holomorphic. Being one-to-one, $f$ is  
a biholomorphism.
\par\smallskip
We give here a simple argument to prove that 
$\Hat{\bold{N}}'=\bold{N}'\bold{Q}_0$.\par
Consider the lower central series
$$\matrix\format\r&\l\\
\Hat{\frak{n}}=\roman{C}^{(0)}(\Hat{\frak{n}}')\supset
\roman{C}^{(1)}(\Hat{\frak{n}}')=\left[\Hat{\frak{n}}',\Hat{\frak{n}}'\right]
\supset &\;\cdots
\supset \roman{C}^{(h)}(\Hat{\frak{n}}')=
\left[\roman{C}^{(h-1)}(\Hat{\frak{n}}'),\Hat{\frak{n}}'\right]\supset \cdots\\
\vspace{2\jot}
&\;\cdots \supset \roman{C}^{(m)}(\Hat{\frak{n}}')=
\left[\roman{C}^{(m-1)}(\Hat{\frak{n}}'),\Hat{\frak{n}}'\right]
=\{0\}\,.\endmatrix$$
Since $\Hat{\frak{n}}'$ is nilpotent and $\Hat{\bold{N}}'$ is
connected, the exponential map
$\exp:\Hat{\frak{n}}'@>>>\Hat{\bold{N}}'$ is surjective.
Let $Z\in\Hat{\frak{n}}'$. We want to prove that there is
$g\in\bold{N}'$ such that $g^{-1}\cdot\exp(Z)\in\bold{Q}_0$. 
To this aim, let $X\in\frak{n}'$ and $W\in\frak{q}_0$ be such that
$Z=X+W$. Let $Z_1\in\Hat{\frak{n}}'$ be such that
$\exp(Z_1)=\exp(-X)\exp(Z)\exp(-W)$.
We claim that, if $Z\in\roman{C}^{h}(\Hat{\frak{n}}')$, then
$Z_1\in\roman{C}^{h+1}(\Hat{\frak{n}}')$.
\par
While proving this claim, we can assume that $\Hat{\bold{N}}'$ is also
simply connected, so that all Lie subgroups 
$\Hat{\bold{N}}'_h=\exp\left(\roman{C}^{(h)}(\Hat{\frak{n}}')\right)$ 
are closed and simply connected. For
each integer $h\geq 0$ we have a commutative diagram\,:
$$\CD
\Hat{\frak{n}}' @>{\exp}>> \Hat{\bold{N}}'\\
@V{\pi_h}VV                          @VV{p_h}V \\
{\Hat{\frak{n}}}'/{\roman{C}^{h+1}(\Hat{\frak{n}}')}
@>{[\exp]}>> {\Hat{\bold{N}}'}/{\Hat{\bold{N}}'_{h+1}}
\endCD
$$
where $[\exp]$ denotes the exponential map on the quotient.
If $Z\in\roman{C}^{h}(\Hat{\frak{n}}')$, 
then $\pi_h(Z)$ belongs to the center
of the quotient Lie algebra 
${\Hat{\frak{n}}'}/{\roman{C}^{h+1}(\Hat{\frak{n}}')}$. Hence we obtain\,:
$$\matrix\format\r&\l\\
[\exp]\left(\pi_h(Z_1)\right)&=
[\exp](-\pi_h(X))\cdot [\exp](\pi_h(Z))\cdot [\exp](\pi_h(X-Z))\\
&=
[\exp](-\pi_h(X)))\cdot [\exp]\left(\pi_h(Z)+\pi_h(X-Z)\right)\\
&=
[\exp](-\pi_h(X)))\cdot [\exp](\pi_h(X)))=
\bold{1}_{{\Hat{\bold{N}}'}/{\Hat{\bold{N}}'_{h+1}}}\,.\endmatrix$$
Since $[\exp]$ is a diffeomorphism, we obtain that
$Z\in\roman{C}^{h+1}(\Hat{\frak{n}}')$.
\par
We show by recurrence that for every
$Z\in\roman{C}^{(m-i)}(\Hat{\frak{n}}')$
there is  some $X\in\frak{n}$ such that $\exp(-X)\cdot\exp(Z)\in
\bold{Q}_0$. This is trivially true when $m=0$, as $Z=0$ in
this case. If $Z\in\roman{C}^{(m-i)}(\Hat{\frak{n}}')$ for some $i>0$,
and $X\in\frak{n}'$ is such that $X-Z\in\frak{q}_0$,
then $\exp(-X)\cdot\exp(Z)\cdot\exp(X-Z)=\exp(Z_1)$ for some
$Z_1\in\roman{C}^{(m-i+1)}(\Hat{\frak{n}}')$. By the recursive assumption,
there is $X_1\in\frak{n}'$ such that
$\exp(-X_1)\exp(Z_1)\in\bold{Q}_0$. 
Then $g=\exp(X_1)\cdot\exp(X)\in\bold{N}'$
and $g^{-1}\cdot\exp(Z)\in\bold{Q}_0$. For $i=m$ we obtain our contention.
\enddimo
From Theorem {\TeoremaID} we obtain\,:
\thm{Let $M=M(\frak{g},\frak{q})$ and $M'=M(\frak{g},\frak{q}')$
be parabolic $CR$ manifolds.
If $\frak{q}\subset\frak{q}'\subset\frak{q}+\bar{\frak{q}}$,
then the $\bold{G}$-equivariant fibration
$M@>>>M'$ is $CR$ and has a totally complex
simply connected fiber.}
\dimo
We already noted in \S{\paragB} that the $CR$ algebra associated to
the fiber $F$ of
the fibration $M@>>>M'$, and hence $F$ itself,
is totally complex when 
$\frak{q}\subset\frak{q}'\subset\frak{q}+\bar{\frak{q}}$.
By Theorem {\TeoremaID} the connected components of the fiber are
the product of a Euclidean complex nilmanifold and of a manifold
$M(\frak{g}'_0,\Hat{\frak{g}}'_0\cap\frak{q})$, for a totally complex parabolic
$CR$ algebra
$(\frak{g}'_0,\Hat{\frak{g}}'_0\cap\frak{q})$.
This $M(\frak{g}'_0,\Hat{\frak{g}}'_0\cap\frak{q})$ is an open orbit of
a connected
real form $\bold{G}'_0$ of a connected complex Lie group $\Hat{\bold{G}}_0'$
with Lie algebra $\Hat{\frak{g}}'_0$, and thus is simply connected
by [W, Theorem 5.4].\enddimo
\cor{Let $(\frak{g},\frak{q}')$ be the weakly nondegenerate reduction
of the effective parabolic $CR$ algebra $(\frak{g},\frak{q})$. 
Then 
\form{f:M=M(\frak{g},\frak{q}) @>>> M'=M(\frak{g},\frak{q}')}
is a $\bold{G}$-equivariant $CR$ fibration
with complex  simply connected fibers.\qed}
We give here a simple general criterion that ensures the 
existence and the connectedness
of the fiber of some $\bold{G}$-equivariant fibration. 
\prop{Let $(\frak{g},\frak{q})$, $(\frak{g},\frak{q}')$ 
be two effective parabolic $CR$ algebras. Assume that
$\frak{g}_+=\frak{q}\cap\frak{g}\subset\frak{g}'_+=\frak{q}'\cap\frak{g}$
and that $\frak{g}_+$ contains a Cartan subalgebra $\frak{h}$ 
that is maximally noncompact among the Cartan subalgebras of
$\frak{g}$ that are contained in $\frak{g}'_+$.
Then the germ of local $\bold{G}$-equivariant submersion 
$(M(\frak{g},\frak{q}),\o)@>>>(M(\frak{g},\frak{q}'),\o')$,
defined by the projection $\frak{g}/\frak{g}_+@>>>\frak{g}/\frak{g}'$,
extends to a 
$\bold{G}$-equivariant fibration 
$\pi:M(\frak{g},\frak{q})@>>>M(\frak{g},\frak{q}')$ with
connected fibers.}
\edef\ProposizioneIF{\number\q.\number\x}\dimo
Let $\bold{H}=\bold{Z}_{\bold{G}}(\frak{h})=
\{h\in\bold{G}\,|\, \roman{Ad}_{\frak{g}}(h)(H)=H\,,\;
\forall H\in
\frak{h}\}$ be the Cartan subgroup of
$\bold{G}$ corresponding to $\frak{h}$. We have
$\roman{Ad}_{\Hat{\frak{g}}}(h)(\frak{q})=\frak{q}$ and
$\roman{Ad}_{\Hat{\frak{g}}}(h)(\frak{q}')=\frak{q'}$ for all
$h\in \bold{H}$. Hence $\bold{H}\subset\bold{G}_+\cap\bold{G}'_+$.
Since $\frak{h}$ is maximally
noncompact in $\frak{g}_+'$ and a fortiori in $\frak{g}_+$,
by [Kn Prop.7.90], all connected components of $\bold{G}'_+$ and
$\bold{G}_+$ intersect $\bold{H}$. The connected component of
the identity $\bold{G}_+^0$ of $\bold{G}_+$ is
contained in the connected component
$[{\bold{G}'_+}]^0$ of the identity in $\bold{G}_+'$. Since
$\bold{G}_+$ is generated by $\bold{G}_+^0$ and $\bold{H}$,
and likewise 
$\bold{G}_+'$ is generated by $[{\bold{G}_+'}]^0$ and $\bold{H}$,
we obtain at the same time that
$\bold{G}_+\subset\bold{G}'_+$ and that 
the fiber $\bold{G}'_+/\bold{G}_+$ is connected.\enddimo
\par
Using Proposition {\ProposizioneIF}, we can prove\,:
\thm{Let $(\frak{g},\frak{q})$ and $(\frak{g},\frak{q}')$ be two
effective parabolic $CR$ algebras. If 
$$\frak{q}\cap\bar{\frak{q}}=\frak{q}'\cap\bar{\frak{q}}\,'\,,$$
then the $CR$ manifolds
$M=M(\frak{g},\frak{q})$ and $M'=M(\frak{g},\frak{q}')$ 
are diffeomorphic, by a $\bold{G}$-equivariant diffeomorphism.}
\edef\TEoremaIG{\number\q.\number\x}\dimo
Let $\frak{h}$ be a Cartan subalgebra of $\frak{g}$, contained
in $\frak{g}_+=\frak{q}\cap\frak{g}=\frak{g}'_+=\frak{q}'\cap\frak{g}$
and maximally noncompact as a Cartan subalgebra of $\frak{g}_+$.\par
Let $A,A'\in\frak{h}_{\Bbb{R}}$ be such that
$\frak{q}=\frak{q}_A$, $\frak{q}'=\frak{q}_{A'}$. We can assume that
$\frak{q}\neq\frak{q}'$, so that $A$ and $A'$ are linearly independent.
Then we set $A_t=A+t(A'-A)$, for $0\leq t\leq 1$, so that
$A_0=A$ and $A_1=A'$. Let
us take a partition $t_0=0<t_1<\cdots<t_{m-1}<t_m=1$ such that
the rank of $\roman{ad}_{\Hat{\frak{g}}}(A_t)$ is constant
for
$t$ in the open intervals $t_{j-1}<t<t_{j}$, $1\leq j\leq m$,
so that $\frak{q}_{A_t}=\frak{q}_{A_{t'}}$ for $t_{j-1}<t,t'<t_j$.
Let $M_j=M(\frak{g},\frak{q}_{A_{t_j}})$, for $0\leq j\leq m$,
and $N_j=M(\frak{g},\frak{q}_{A_{(t_{j-1}+t_j)/2}})$, for
$1\leq j\leq m$. Since\,:
$\frak{q}_{A_{(t_{j-1}+t_j)/2}}\subset\frak{q}_{A_{t_{j-1}}}\cap
\frak{q}_{A_{t_{j}}}$, there are $\bold{G}$-equivariant maps\,:
$$\CD M_{j-1} @<{f_j}<< N_j @>{F_j}>> M_j\endCD
$$
for all $1\leq j\leq m$. By Proposition {\ProposizioneIF}, all these
maps, being covering maps with connected fibers, are diffeomorphisms.
Thus\,: 
$$(F_m\circ f_{m}^{-1})\circ (F_{m-1}\circ f_{m-1}^{-1})\circ\cdots
\circ (F_1\circ f_1^{-1}):M\longrightarrow M'$$
is a $\bold{G}$-equivariant diffeomorphism.
\enddimo


\se{Fit Weyl chambers and $CR$ geometry of $M(\frak{g},\frak{q})$}
Let $(\frak{g},\frak{q})$ be an effective parabolic $CR$ algebra and
$(\vartheta,\frak{h})$ a Cartan pair adapted to
$(\frak{g},\frak{q})$.
We keep the notation of the preceding sections, for roots,
parabolic sets, Cartan decomposition, etc.
In particular, we denote by $\sigma:\frak{h}^*_{\Bbb{R}}\ni
\alpha @>>>\bar\alpha\in\frak{h}^*_{\Bbb{R}}$ the adjoint map of
the restriction to 
$\frak{h}_{\Bbb{R}}=\frak{h}^-\oplus i\,\frak{h}^+$
of the conjugation in $\Hat{\frak{g}}$ defined by the real form
$\frak{g}$.  We say that a root
$\alpha$ is
{\it real} if $\bar\alpha=\alpha$,
{\it imaginary} if $\bar\alpha=-\alpha$,
{\it complex} if $\bar\alpha\neq\pm\alpha$
and denote by $\rootre$, $\rootim$ and $\rootcp$ the sets of real,
imaginary and complex roots, respectively. When $\alpha$ is imaginary,
the eigenspace $\Hat{\frak{g}}^{\alpha}$ is contained either in
$\Hat{\frak{k}}=\Bbb{C}\otimes_{\Bbb{R}}\frak{k}$, 
or in $\Hat{\frak{p}}=\Bbb{C}\otimes_{\Bbb{R}}\frak{p}$.
In the first case we say that $\alpha$ is {\it compact}, in the
second that $\alpha$ is {\it noncompact}. Thus $\rootim$ is the disjoint
union of the set $\rootimc$ of
the compact and of the set
$\rootimp$ of the noncompact imaginary roots\,:
$\rootim=\rootimc\cup\rootimp$.
\subhead \number\q.1. S-fit and V-fit Weyl chambers \endsubhead\par\noindent
The conjugation $\sigma$ defines an involution in $\frak{h}^*_{\Bbb{R}}$
that belongs to
$\alaff$. Vice versa, 
every involution $\sigma$ in $\alaff$ can be obtained 
from a conjugation with
respect to a real form $\frak{g}$ of $\Hat{\frak{g}}$. 
Note that, in general, $\sigma$ does not uniquely determine the isomorphism
class of $\frak{g}$. 
Let us describe the structure of an arbitrary involution 
$\sigma$ in $\alaff$\,:
\thm{Let $\sigma$ be an involution in $\alaff$.
Then there exist: a set 
of
pairwise strongly
orthogonal roots $\alpha_1\,,\;\hdots\, ,\;\alpha_m$ 
in $\Cal{R}$, with $\sigma(\alpha_j)=-\alpha_j$
for $j=1,\hdots,m$, a Weyl chamber $C\in\frak{C}(\Cal{R})$,
and an involution 
$\jmath\in\alaffine{C}$, with $\jmath(\alpha_i)=\alpha_i$, and
hence  
commuting with $s_{\alpha_i}$, for all $i=1,\hdots,m$,
such that\,:
\form{\sigma\,=\, \jmath\circ s_{\alpha_1}\circ \cdots\circ s_{\alpha_m}\,;}
\form{\alpha\in\Cal{R}^+(C) \Longrightarrow
\cases
\text{either $\sigma(\alpha)=-\alpha$}\\
\text{or $\sigma(\alpha)\in\Cal{R}^+(C)$}
\endcases}
}\edef\formulcamera{\number\q.\number\t} 
\edef\teoremaCI{\number\q.\number\x}\dimo
Let $F^{-}(\sigma)=\{\alpha\in\frak{h}^*_{\Bbb{R}}\,
|\,\sigma(\alpha)=-\alpha\}$.
Take a maximal subset $\alpha_1,\hdots,\alpha_m$ of 
pairwise orthogonal roots
in $F^-(\sigma)\cap\Cal{R}$ and consider 
$\jmath=\sigma\circ s_{\alpha_1}\circ\cdots\circ s_{\alpha_m}$.
We have $\jmath(\alpha_i)=\sigma(-\alpha_i)=\alpha_i$ 
for all $i=1,\hdots,m$. We claim
that $\jmath(\alpha)\neq -\alpha$ for all $\alpha\in\Cal{R}$.
Indeed, if there was $\alpha\in\Cal{R}$ with $\jmath(\alpha)=-\alpha$,
from
$(\alpha|\alpha_i)=(\jmath(\alpha)|\jmath(\alpha_i))=-(\alpha|\alpha_i)$
we obtain that $(\alpha|\alpha_i)=0$ for all $i=1,\hdots,m$. Hence
$s_{\alpha_i}(\alpha)=\alpha$ for all $\alpha$ and therefore
$\sigma(\alpha)=\jmath(\alpha)=-\alpha$, contradicting the
fact that $\alpha_1,\hdots,\alpha_m$ was a maximal
system of pairwise orthogonal roots in $\Cal{R}\cap F^-(\sigma)$.
\par
\smallskip
To obtain that $\alpha_1,\hdots,\alpha_m$ are  strongly
orthogonal it suffices to choose the sequence $\alpha_1,\hdots,\alpha_m$
with a maximal sum $\sum_{i=1}^m{\|\alpha_i\|^2}$.
Indeed, if $\alpha_j$ and $\alpha_h$ are orthogonal, but not strongly
orthogonal, then both $\alpha_j+\alpha_h$ and $\alpha_j-\alpha_h$
are roots. Setting $\alpha'_i=\alpha_i$ for $i\neq j,h$, and
$\alpha'_j=\alpha_j+\alpha_h$\,, $\alpha'_h=\alpha_j-\alpha_h$,
we obtain a new sequence $\alpha'_1,\hdots,\alpha'_m$ of pairwise
orthogonal roots in $F^-(\sigma)\cap\Cal{R}$. It is contained in a 
maximal one and the inequality
$\sum_{i=1}^m{
\|\alpha_i'\|^2}=\left(\sum_{i=1}^m{\|\alpha_i\|^2}\right)
+\|\alpha_j\|^2+\|\alpha_h\|^2>\sum_{i=1}^m{\|\alpha_i\|^2}$,
contradicts the maximality of $\sum_{i=1}^m{\|\alpha_i\|^2}$.
\par\smallskip
We claim that there exists a Weyl chamber $C$ such that\,:
$$\jmath(\Cal{R}^+(C))=\Cal{R}^+(C)\, . \tag $*$ $$
Indeed ($*$) is equivalent to $\jmath(\Cal{B}(C))\subset\Cal{R}^+(C)$.
For a Weyl chamber $C$, 
denote by $n_C$ the number of the elements
in $\Cal{R}^+(C)\cap\jmath(\Cal{R}^+(C))$. Fix
$C$ with $n_C$ maximum.
 If $C$ does not satisfy
($*$), take $\alpha\in\Cal B(C)$ with 
$\jmath(\alpha)\notin\Cal{R}^+(C)$ and consider the
chamber $C'=s_{\alpha}(C)$.
From $\Cal{R}^+(C')=\left(\Cal{R}^+(C)\setminus\{\alpha\}\right)
\cup\{-\alpha\}$ and $\jmath(-\alpha)\in\Cal{R}^+(C)\setminus\{\alpha\}
\subset\Cal{R}^+(C')$, we obtain
$n_{C'}=n_C+1$, contradicting our choice of $C$. Hence 
$C$ satisfies ($*$) and therefore also
$(\formulcamera)$.
This completes the proof. \enddimo
Using Theorem {\teoremaCI}, we obtain the formula\,:
\form{\sigma(\beta)=\jmath(\beta)-
\sum_{j=1}^m{\left(\beta|\alpha_j^{\lor}\right)\,
\alpha_j}\,,\quad\text{with}\quad\alpha_j^{\lor}=
2\alpha_j/\|\alpha_j\|^2\,,
\quad\forall\beta\in\frak{h}^*_{\Bbb{R}}\,.}

Likewise, we have the following\,:
\thm{Let $\sigma$ be an involution in $\alaff$. Then there
exists a set of pairwise strongly orthogonal roots 
$\delta_1,\hdots,\delta_m\in\Cal{R}$, 
with $\sigma(\delta_j)=\delta_j$ 
for $j=1,\hdots,m$, a Weyl chamber $C\in\frak{C}(\Cal{R})$
and
an involution
$\varpi$, that commutes with $s_{\delta_j}$, satisfies
$\varpi(\delta_j)=-\delta_j$
for all $j=1,\hdots,m$, 
and transforms $C$ into $C^{\roman{opp}}$, such that\,:
\form{\sigma=\varpi\circ s_{\delta_1}\circ\cdots 
\circ s_{\delta_m}\,,}
\form{\alpha\in\Cal{R}^+(C)\Longrightarrow
\cases
\text{either}\; \sigma(\alpha)=\alpha\\
\text{or}\; \sigma(\alpha)\in\Cal{R}^-(C)\,.
\endcases}
}\edef\teoremaCII{\number\q.\number\x}\edef\FormulaEE{\number\q.\number\t}
\dimo
We take $\sigma'=s_0\circ\sigma$, where $s_0$ is the symmetry with
respect to the origin of $\frak{h}^*_{\Bbb{R}}$.
By the preceding Theorem, $\sigma'=\jmath\circ s_{\delta_1}\circ\cdots
\circ s_{\delta_m}$, where $\jmath\in\alaffine{C}$ for some
$C\in\frak{C}(\Cal{R})$, and $\delta_1,\hdots,\delta_m$ is
a maximal system of strongly orthogonal
roots  in
$F^-(\sigma')\cap\Cal{R}=\{\alpha\in\Cal{R}\,|\,
\sigma(\alpha)=\alpha\}$, with $\jmath(\delta_j)=\delta_j$.
The statement follows by taking $\varpi=s_0\circ\jmath$.
\enddimo 
\medskip
With $\varpi$ and $\delta_1,\,\hdots\,,\delta_m$ as in 
Theorem {\teoremaCII}, we obtain the formula\,:
\form{\sigma(\beta)=\varpi(\beta)
+\sum_{j=1}^m{\left(\beta|\delta_j^{\lor}\right)\,
\delta_j}\,,\quad\text{with}\quad
\delta_j^{\lor}=2\,\delta_j/\|\delta_j\|^2\,,\quad
\forall\beta\in\frak{h}^*_{\Bbb{R}}\,.}

A Weyl chamber $C\in\frak{C}(\Cal{R})$ that satisfies condition
$(\formulcamera)$ (resp. $(\FormulaEE)$) is said to be
\footnote{If we choose $\frak{h}$ maximally noncompact,
then in an S-adapted Weyl chamber $C$ the conjugation can be
described by a Satake diagram; if instead we take $\frak{h}$ with a
maximal compact part, in a V-adapted Weyl chamber the conjugation is
described by a Vogan diagram (see e.g. [Ar], [Kn]).}
 {\it S-adapted}
(resp. {\it V-adapted}) to the conjugation $\sigma$.\par
\smallskip
Our description of the minimal orbits in [AMN] in terms of
{\it cross-marked Satake diagrams}
used the fact
that, for a parabolic minimal $(\frak{g},\frak{q})$, there are
a maximally noncompact
Cartan subalgebra $\frak{h}$ of $\frak{g}$,
adapted to $(\frak{g},\frak{q})$, and an S-adapted Weyl chamber $C$ in
$\frak{C}(\Cal{R},\Cal{Q})$.
For general $CR$ algebras $(\frak{g},\frak{q})$, 
there could be no adapted Cartan subalgebras $\frak{h}$
that are either maximally compact or maximally noncompact in $\frak{g}$\,.
This is a major drawback in the classification of the orbits of
$\bold{G}$ in $\frak{M}$
(see e.g. the references in [BL]), but, while discussing fundamentality,
weak nondegeneracy and some topological properties,
it turns out that the choice of 
$\frak{h}$ is not as crucial as that of 
special Weyl chambers $C$ in 
$\frak{C}(\Cal{R},\Cal{Q})$. 
In general 
$\frak{C}(\Cal{R},\Cal{Q})$ may not contain any 
Weyl chamber that is either S- or  V-adapted to $\sigma$.
In the following lemmas we describe chambers in
$\frak{C}(\Cal{R},\Cal{Q})$ that are as close as possible to being
S- or V-adapted.
\lem{Let $(\vartheta,\frak{h})$ be a Cartan pair adapted
to $(\frak{g},\frak{q})$. Then there exists a Weyl chamber
$C\in\frak{C}(\Cal{R},\Cal{Q})$ that satisfies the equivalent
conditions\,:
\roster
\item"$(i)$" If $\alpha\succ_C 0$ and $\bar\alpha\prec_C 0$, then
both $\alpha$ and $-\bar\alpha$ belong to $\Cal{Q}^n$\,.
\item"$(ii)$" $\bar\alpha\succ_C 0$ for all 
$\alpha\in\Cal{B}(C)\setminus\left(\Phi_C\cup\rootim\right)$\,.
\endroster
Assume that $C$ satisfies the equivalent conditions $(i)$ and $(ii)$. 
Then\,:
\roster
\item"$(iii)$" If moreover $\frak{h}$
is maximally noncompact among the Cartan subalgebras of $\frak{g}$
contained in $\frak{g}_+=\frak{q}\cap\frak{g}$, then 
$\Cal{B}(C)\cap\rootimp\subset\Phi_C$.
\endroster} 
\dimo Choose $C\in\frak{C}(\Cal{R},\Cal{Q})$ with
a maximal
$\Cal{R}^+(C)\cap\sigma\left(\Cal{R}^+(C)\right)$. 
Then $(ii)$ is satisfied. Indeed, if there was
$\alpha\in\Cal{B}(C)\setminus\left(\Phi_C\cup\rootim\right)$
with $\bar\alpha\prec_C 0$, we would take $C'=s_{\alpha}(C)$. Then
$\Cal{R}^+(C')=\left(\Cal{R}^+(C)\setminus\{\alpha\}\right)\cup\{-\alpha\}
\subset\Cal{Q}$, so that $C'\in\frak{C}(\Cal{R},\Cal{Q})$,
and
$\Cal{R}^+(C')\cap \sigma\left(\Cal{R}^+(C')\right)
\supsetneqq
\Cal{R}^+(C)\cap\sigma\left(\Cal{R}^+(C)\right)$, yielding a contradiction.
Clearly $(i)\Rightarrow (ii)$. 
Vice versa, if $\alpha=\sum_{\beta\in\Cal{B}(C)}{k^{\beta}_{\alpha}\beta}
\in\Cal{R}^+(C)$ and $\bar\alpha\prec_{C}0$\,,
then either $\alpha\in\rootim$, or else there is some 
$\beta\in\supp{C}(\alpha)\cap\rootcp$ with $\bar\beta\prec_{C}0$\,;
by $(ii)$, we have $\beta\in\Phi_{C}$ and hence $\alpha\in\Cal{Q}^n$.
The same argument, applied to $-\bar\alpha$, shows that also
$-\bar\alpha\in\Cal{Q}^n$. This completes the proof of the equivalence
$(i)\Leftrightarrow(ii)$.
\par\smallskip
Finally, if $\alpha\in\left(\Cal{B}(C)\cap\rootimp\right)\setminus\Phi_C$,
both $\alpha$ and $\bar\alpha=-\alpha$ belong to $\Cal{Q}$.
Let $\Gamma=\{(X_{\alpha},H_{\alpha})_{\alpha\in\Cal{R}}\}$ be a
Chevalley system, as in [B2].
Then $T_{\alpha}=X_{\alpha}-X_{-\alpha}\in\frak{p}\cap\frak{g}_+$
(for this construction cf. [S1], [S2]) 
is a semisimple
element of $\frak{g}$
that commutes with all elements of $\frak{h}^-$ and of
$\frak{j}^+=\{H\in\frak{h}^+\,|\, \alpha(H)=0\}$. 
Hence $\frak{j}=\frak{h}^-\oplus\Bbb{R}\, T_{\alpha}\oplus\frak{j}^+$
is a Cartan subalgebra of $\frak{g}$, contained in $\frak{g}^+$,
with $\frak{j}^-=\frak{h}^-\oplus\Bbb{R}\,T_{\alpha}
\supsetneqq\frak{h}^-$. Thus, if
$\frak{h}^-$ is maximal, we have 
$\Cal{B}(C_0)\cap\rootimp\subset\Phi_C$.
\enddimo\edef\lemmaia{\number\q.\number\x}
An alternative construction of a
Weyl chamber $C$ satisfying $(i)$ and $(ii)$ of Lemma {\lemmaia}
is
the following 
(which is a particular application of the general construction that
will be described at the beginning of \S{6}).
Fix a Weyl chamber $C_0\in\frak{C}(\Cal{R})$ that is S-adapted to $\sigma$
(this means that $\sigma\left(\Cal{R}^+(C_0)\setminus\rootim
\right)\subset\Cal{R}^+(C_0)$), and consider the Borel subalgebra\,:
$$\frak{b}_{C_0}=\Hat{\frak{h}}\oplus\bigoplus_{\alpha\in\Cal{R}^+(C_0)}{
\Hat{\frak{g}}^{\alpha}}\subset\Hat{\frak{g}}\,.$$
Then  $\frak{b}=\frak{q}^n\oplus\left(\frak{q}^r\cap\frak{b}_{C_0}\right)$ 
is a Borel subalgebra of $\Hat{\frak{g}}$, corresponding to a
Weyl chamber $C\in\frak{C}(\Cal{R},\Cal{Q})$ that
satisfies $(i)$ and $(ii)$ of Lemma {\lemmaia}.
\par\smallskip
A Weyl chamber $C\in\frak{C}(\Cal{R},\Cal{Q})$ that satisfies the 
equivalent 
conditions $(i)$ and $(ii)$ 
of Lemma {\lemmaia} 
is called {\it S-fit} to $(\frak{g},\frak{q})$.
\medskip
With arguments similar to those employed to prove Lemma {\lemmaia},
we obtain\,:
\lem{Let $(\vartheta,\frak{h})$ be a Cartan pair adapted
to $(\frak{g},\frak{q})$. Then there exists a Weyl chamber
$C\in\frak{C}(\Cal{R},\Cal{Q})$ that satisfies the equivalent
conditions\,:
\roster
\item"$(iv)$" If $\alpha\succ_C 0$ and $\bar\alpha\succ_C 0$, then
both $\alpha$ and $\bar\alpha$ belong to $\Cal{Q}^n$\,.
\item"$(v)$" $\bar\alpha\prec_C 0$ for all 
$\alpha\in\Cal{B}(C)\setminus\left(\Phi_C\cup\rootre\right)$\,.
\endroster
Assume that $C$ satisfies the equivalent conditions $(iv)$ and $(v)$. 
Then\,:
\roster
\item"$(vi)$" If moreover $\frak{h}$
is maximally compact among the Cartan subalgebras of $\frak{g}$
contained in $\frak{g}_+=\frak{q}\cap\frak{g}$, then 
$\Cal{B}(C)\cap\rootre\subset\Phi_C$.\qed
\endroster} 
\edef\lemmaIB{\number\q.\number\x}

A Weyl chamber $C\in\frak{C}(\Cal{Q},\Cal{R})$ that satisfies
the equivalent conditions $(iv)$ and $(v)$ of Lemma {\lemmaIB}
is called {\it V-fit} to $(\frak{g},\frak{q})$.\par
V-fit Weyl chambers can be obtained as follows\,: if 
$\frak{b}_{C_0}=\Hat{\frak{h}}
\oplus\sum_{\alpha\in\Cal{R}^+(C_0)}{\Hat{\frak{g}}^{\alpha}}$ 
is the Borel subalgebra associated to a Weyl chamber
$C_0$ that is V-adapted to the conjugation $\sigma$, then 
$\frak{q}^n\oplus\left(\frak{q}^r\cap\frak{b}_{C_0}\right)$ 
is the Borel subalgebra
associated to a Weyl chamber 
$C\in\frak{C}(\Cal{R},\Cal{Q})$ that is
V-fit to $(\frak{g},\frak{q})$.
\subhead \number\q.2. Fundamental parabolic $CR$ algebras \endsubhead\par\noindent
\thm{Let $(\frak{g},\frak{q})$ be an effective parabolic $CR$ algebra.
Fix a Cartan subalgebra $\frak{h}$ of $\frak{g}$ adapted to
$(\frak{g},\frak{q})$ and fix an S-fit Weyl chamber 
$C\in\frak{C}(\Cal{R},\Cal{Q})$ for $(\frak{g},\frak{q})$. Then
$(\frak{g},\frak{q})$ is fundamental if and only if\,:
\form{\forall\alpha_0\in\Phi_C\cap\bar{\Cal{Q}}\,^n
\;\cases
\text{either}\;\exists\beta\in\Cal{B}(C)\setminus\Phi_C\;\text{with}\;
\alpha_0\in\supp{C}(\bar\beta)\,,\\
\text{or}\;\exists\beta\in\Phi_C\;\text{with}\;
\bar\beta\in\Cal{R}^-(C)\;\text{and}\; \alpha_0\in\supp{C}(\bar\beta)\,.
\endcases}}\edef\FormulaEI{\number\q.\number\t}
\edef\TeorEF{\number\q.\number\x}
\dimo The parabolic $CR$ algebra $(\frak{g},\frak{q})$ is fundamental
if, and only if, there is no complex parabolic subalgebra $\frak{q}'$
of $\Hat{\frak{g}}$ with 
$\frak{q}+\bar{\frak{q}}\subset\frak{q}'\subsetneq\Hat{\frak{g}}$.
All complex parabolic $\frak{q}'$ that contain $\frak{q}$ are of the
form $\frak{q}'=\frak{q}_{\Psi_C}$ for some set of simple roots
$\Psi_C\subset\Phi_C$. We can limit ourselves to consider the cases where
$\Psi_C=\{\alpha_0\}$ for some $\alpha_0\in\Phi_C$. Thus we obtain\,:
\form{(\frak{g},\frak{q}_{\Phi_C})\;\text{is not fundamental}
\Longleftrightarrow 
\cases\exists\alpha_0\in\Phi_C\;\text{such that}\\
\{\beta\in\Cal{R}\,|\, \beta\succeq_C\alpha_0\}
\subset
\bar{\Cal{Q}}\,^n\,.\endcases}
Indeed, $\frak{q}+\bar{\frak{q}}\subset \frak{q}'$ if and only if
$\frak{q}'\,^n\subset\frak{q}^n\cap\bar{\frak{q}}\,^n$.
Thus it suffices to further restrict our consideration
to simple roots $\alpha_0\in\Phi_C\cap\bar{\Cal{Q}}^{\,n}$
and check whether $\Cal{F}=\{\beta\in\Cal{R}\,|\,
\beta\succeq_C\alpha_0\}$ is contained or not in
$\bar{\Cal{Q}}^{\,n}$. \par
First we show that condition $(\FormulaEI)$ is sufficient.
Let $\alpha_0\in\Phi_C\cap\bar{\Cal{Q}}^n$.
If $\beta\in\Cal{B}(C)\setminus\Phi_C$ and $\alpha_0\in\supp{C}(\bar\beta)$,
then 
$\bar\beta\succeq_C\alpha_0$
by the assumption that $C$ is S-fit
to $(\frak{g},\frak{q})$\,; 
hence $\bar\beta\in\Cal{F}\setminus\bar{\Cal{Q}}^n$.
Likewise, if $\beta\in\Phi_C$, $\bar\beta\prec_C 0$ and
$\alpha_0\in\supp{C}(\bar\beta)$, then 
$-\bar\beta\in\Cal{F}\setminus\bar{\Cal{Q}}\,^n$.
This completes the proof of sufficiency.\par
To prove that $(\FormulaEI)$ is also necessary,  
fix again $\alpha_0\in\Phi_C\cap\bar{\Cal{Q}}^{\,n}$. If $(\frak{g},\frak{q})$
is fundamental, then $\Cal{F}\not\subset\bar{\Cal{Q}}^{\, n}$, and hence
there is a root $\alpha$ with $\alpha\succeq_C\alpha_0$
and
$\bar\alpha\notin{\Cal{Q}}^n$. If $\bar\alpha\succ_C 0$,
then $\supp{C}(\bar\alpha)\cap\Phi_C=\emptyset$. Since
$\alpha_0\in\supp{C}({\alpha})\subset\bigcup_{\beta\in\supp{C}(\bar\alpha)}{
\supp{C}(\bar\beta)}$, there is at least a $\beta\in\Cal{B}(C)\setminus
\Phi_C$ with $\alpha_0\in\supp{C}(\bar\beta)$. Assume now that
$\alpha_0$ does not belong to $\supp{C}(\bar\beta)$ for any
$\beta\in\Cal{B}\setminus\Phi_C$. Then, for a root $\alpha\succ_C\alpha_0$
that does not belong to $\bar{\Cal{Q}}^{\,n}$, we have
$\bar\alpha\prec_C 0$.
From $\alpha_0\in\bigcup_{\beta\in\supp{C}(
\bar\alpha)}{\supp{C}(\bar\beta)}$ we obtain that $\alpha_0\in\supp{C}
(\bar\beta)$ for some $\beta\in\Cal{B}(C)$ with $\bar\beta\prec_C 0$,
and hence in $\Phi_C$. The
proof is complete.\enddimo
Theorem {\TeorEF} provides a criterion, only
involving the conjugation of the simple roots of an S-fit Weyl chamber,
for an effective parabolic 
$CR$ algebra to be totally real. We found convenient to formulate the
criterion for an arbitrary Weyl chamber $C\in\frak{C}(\Cal{R},\Cal{Q})$.
\lem{Let $(\frak{g},\frak{q})$ be an effective parabolic $CR$ algebra,
and $C\in\frak{C}(\Cal{R},\Cal{Q})$. A necessary and sufficient
condition in order that $(\frak{g},\frak{q})$ be totally real is that\,:
\form{\left\{
\matrix
 \alpha\in\Cal{B}(C)\setminus\Phi_C&\;\Longrightarrow\;&
\roman{supp}_C(\bar\alpha)\cap\Phi_C=\emptyset\\\vspace{1\jot}
\alpha\in\Phi_C&\;\Longrightarrow \;& \bar\alpha\in\Cal{R}^+(C)\,.
\endmatrix\right.}}\edef\formulaIF{\number\q.\number\t}
\dimo 
The case where $\Cal{Q}=\Cal{R}$ is trivial. Assume that $\Cal{Q}\neq\Cal{R}$.
The first condition in (\number\q.\number\t) implies that 
$\Cal{B}(C)\setminus\Phi_C \subset\bar\Cal{Q}^{\, r}$.
Hence $\Cal{Q}^r=\bar\Cal{Q}^{\, r}$. In particular,
if $\alpha\in\Cal{Q}^n$, then $\bar\alpha\in\Cal{R}\setminus\Cal{Q}^r$.
Then, since $\Cal{Q}^n=\Cal{R}^+(C)\setminus\Cal{Q}^r$, 
the second condition implies
that $\bar\alpha\in\Cal{Q}^n$ for all $\alpha\in\Phi_C$.
Hence $\Phi_C\subset\bar\Cal{Q}^n$. Therefore
$\Cal{B}(C)\subset\overline{\Cal{Q}}$, so that we also have
$\overline{\Cal{Q}}^{\, n}={\Cal{Q}}^n$ and hence $\Cal{Q}=
\overline{\Cal{Q}}$. The condition is obviously also necessary.
\enddimo
\medskip\edef\lemmaid{\number\q.\number\x}

Lemma {\lemmaid} prompts a recursive method to construct the totally
real basis of the canonical $\frak{g}$-equivariant fibration 
$(\frak{g},\frak{q})@>>>
(\frak{g},\frak{q}')$, with totally real basis and 
fundamental fiber.  
\par
After taking any
$C\in\frak{C}(\Cal{R},\Cal{Q})$, we define recursively\,:
\form{\left\{\matrix\format\l&\c&\;\l\\
\Upsilon_C^{(0)}&=&\{\alpha\in\Phi_C\,|\, \bar\alpha\in\Cal{R}^+(C)\}\\
\vspace{1\jot}
\Upsilon_C^{(1)}&=& \Upsilon_C^{(0)}\setminus\bigcup_{\alpha\in\Cal{B}(C)
\setminus
\Upsilon_C^{(0)}}\supp{C}(\bar\alpha)\\
\vspace{1\jot}
\Upsilon_C^{(h+1)}&=& \Upsilon_C^{(h)}\setminus\bigcup_{\alpha\in\Cal{B}(C)
\setminus
\Upsilon_C^{(h)}}\supp{C}(\bar\alpha)\quad\text{for}\quad h\geq 1
\\
\vspace{2\jot}
\Upsilon_C&=&\bigcap_{h\geq 0}\Upsilon_C^{(h)}\quad\text{(finite
intersection)}
\,.
\endmatrix\right.}\edef\FFrmEJ{\number\q.\number\t}
One easily verifies, using the previous results, that\,:
\prop{The natural $\frak{g}$-equivariant fibration 
$(\frak{g},\frak{q}_{\Phi_C})@>>>(\frak{g},\frak{q}_{\Upsilon_C})$
is the fundamental reduction of $(\frak{g},\frak{q}_{\Phi_C})$.
In particular, $(\frak{g},\frak{q}_{\Phi_C})$ is
fundamental if and only if $\Upsilon_C=\emptyset$. \qed}
\bigskip
\subhead \number\q.3. Weakly nondegenerate parabolic $CR$ algebras \endsubhead\par\noindent
We turn now to weak nondegeneracy for an effective parabolic $CR$
algebra $(\frak{g},\frak{q})$. We shall see 
that this property 
can be better examined in terms of 
V-fit Weyl chambers.
\lem{Let $(\frak{g},\frak{q})$ be an effective
parabolic $CR$ algebra, $\frak{h}$ a Cartan
subalgebra 
of $\frak{g}$ 
adapted to $(\frak{g},\frak{q})$, and $C\in\frak{C}(\Cal{R},\Cal{Q})$
a V-fit Weyl chamber for $(\frak{g},\frak{q})$.
Let $\Phi_C=\Cal{B}(C)\cap\Cal{Q}^n$
and $\frak{q}'=\frak{q}_{\Psi_C}\supset\frak{q}$ 
(cf. $(\Fformcorr)$ for the notation),
with
$\Psi_C\subset\Phi_C$. Then the $\frak{g}$-equivariant
homomorphism of $CR$ algebras
$(\frak{g},\frak{q})@>>>(\frak{g},\frak{q}')$ is
a $CR$ fibration with a totally complex fiber if and only if
\form{\bar\alpha\prec_C 0 \quad\forall\alpha\in\Phi_C\setminus
\Psi_C\,.}}\edef\formulazza{\number\q.\number\t}
\edef\lemmazzo{\number\q.\number\x}
\dimo
The $\frak{g}$-equivariant $CR$ homomorphism 
$(\frak{g},\frak{q})@>>>(\frak{g},\frak{q}')$ is
a $CR$ fibration 
with a totally complex fiber 
if, and only if, $\frak{q}\subset\frak{q}'\subset
\frak{q}+\bar{\frak{q}}$. Thus, we need to show that these inclusions
are equivalent to $(\formulazza)$ when $C\in\frak{C}(\Cal{R},\Cal{Q})$
satisfies $(iv)$ and $(v)$ of Lemma {\lemmaIB}.\par
Clearly, it suffices to consider the case where the difference
$\Phi_C\setminus\Psi_C$ consists of a single simple root $\alpha_0$.
We shall assume in the following that
$\Phi_C\setminus\Psi_C=\{\alpha_0\}$.\par\smallskip
First we prove that, if
$\bar\alpha_0\prec_C 0$, then
$\Cal{Q}'=\Cal{Q}_{\Psi_C}\subset\Cal{Q}\cup\bar{\Cal{Q}}$\,. 
Assume by contradiction that there is
$\beta\in\Cal{Q}'\setminus\left(\Cal{Q}\cup\bar{\Cal{Q}}\right)$.
Then $-\beta,-\bar\beta\in\Cal{Q}^n\subset\Cal{R}^+(C)$.
Moreover $\beta\in{\Cal{Q}'}^r$, because
${\Cal{Q}'}^n\subset\Cal{Q}^n\subset
\left(\Cal{Q}\cup\bar{\Cal{Q}}\right) $\,.
Thus
$\supp{C}({\beta})\cap\Psi_C=\emptyset$, and hence
$\supp{C}(\beta)\cap\Phi_C=\{\alpha_0\}$,
because $-\beta\in\Cal{Q}^n$.
As $C$ is V-fit,
$\bar\alpha\prec_C 0$ for all 
$\alpha\in\supp{C}(\beta)\setminus\rootre$. Since $\beta$ is not real,
this implies that $\bar\beta\succ_C 0$, giving a contradiction.
Hence $\Cal{Q}'\subset\left(\Cal{Q}\cup\overline{\Cal{Q}}\right)$,
and
this proves that $\frak{q}\subset\frak{q}'\subset\frak{q}+\bar{\frak{q}}$
when $\bar\alpha_0\in\Cal{R}^-(C)$.
\smallskip
Vice versa, assume that
$\frak{q}\subset\frak{q}'\subset\frak{q}+\bar\frak{q}$. 
In particular, $-\alpha_0\in\bar{\Cal{Q}}$. 
If $\bar\alpha_0\in{\Cal{Q}}^r$, then $\bar\alpha_0$ and
$\alpha_0=\Bar{\Bar{\alpha}}_0$ belong to opposite cones
$\Cal{R}^{\pm}(C)$ and thus
$\bar\alpha_0\prec_C 0$. When $\bar\alpha_0\notin\Cal{Q}^r$,
we have
$-\bar\alpha_0\in\Cal{Q}^n\subset\Cal{R}^+(C)$,
and thus
$\bar\alpha_0\prec_C 0$.
The proof is complete.
\enddimo
\thm{Let $(\frak{g},\frak{q})$ be an effective
parabolic $CR$ algebra, and let $\frak{h}$ be a Cartan
subalgebra of $\frak{g}$ adapted to
$(\frak{g},\frak{q})$. Let $C\in\frak{C}(\Cal{R},\Cal{Q})$ be
V-fit to $(\frak{g},\frak{q})$. 
Then $(\frak{g},\frak{q})$ is weakly nondegenerate if and only if\,:
\form{\bar\alpha\succ_C 0
\qquad\forall\alpha\in\Phi_C\,.}}
\edef\formulaid{\number\q.\number\t}
\edef\TEOEI{\number\q.\number\x}
\dimo
Fix a Weyl chamber $C\in\frak{C}(\Cal{R},\Cal{Q})$ for which
$(iv)$ and $(v)$ of Lemma {\lemmaIB} are valid.
By Lemma {\lemmazzo}, 
the necessary and sufficient condition for
$(\frak{g},\frak{q})$
to be weakly degenerate is that there exists $\alpha_0\in\Phi_C$
contradicting $(\formulaid)$.\enddimo 
\cor{Let $(\frak{g},\frak{q})$ be a weakly nondegenerate parabolic
$CR$ algebra. Let $\Cal{Q}$ be the parabolic set associated to $\frak{q}$
in $\Cal{R}=\Cal{R}(\Hat{\frak{g}},\Hat{\frak{h}})$, for an admissible
Cartan subalgebra $\frak{h}$ of $(\frak{g},\frak{q})$.
Then $(\frak{g},\frak{q})$ is totally real if and only if
$\bar{\Cal{Q}}^{\,n}\cap\Cal{Q}^r=\emptyset$.}\edef\corIF{\number\q.\number\x}
\dimo Choose a V-fit Weyl chamber
$C\in\frak{C}(\Cal{R},\Cal{Q})$ for $(\frak{g},\frak{q})$. 
The second condition in $(\formulaIF)$ being automatically satisfied
because $C$ is V-fit and $(\frak{g},\frak{q})$ weakly nondegenerate,
we observe that the first line in 
$(\formulaIF)$ is equivalent to the condition that
$\bar{\Cal{Q}}^{\,n}\cap\Cal{Q}^r=\emptyset$.\enddimo
\bigskip
\subhead \number\q.4. A criterion of minimality \endsubhead\par\noindent
In [AMN] we only used maximally noncompact adapted Cartan subalgebras.
Here
we give a
characterization of
effective parabolic minimal $CR$ algebras $(\frak{g},\frak{q})$ 
in terms of
the set of roots $\Cal{Q}$
associated to $\frak{q}$
by {\it any} choice of an
adapted Cartan pair $(\vartheta,\frak{h})$.
\par\medskip
\prop{A necessary and sufficient condition in order that 
a weakly nondegenerate
effective parabolic $CR$ algebra $(\frak{g},\frak{q})$
be parabolic minimal\footnote{cf. [AMN]. Recall that
the condition is equivalent to the fact that the corresponding
$CR$ manifold $M(\frak{g},\frak{q})$ 
is the minimal orbit.}, is that 
\form{\Cal{Q}^n\cap\bar{\Cal{Q}}^{\,-n}\subset\Cal{R}_{\bullet}\,.}}
\edef\formuih{\number\q.\number\t}\edef\PROPie{\number\q.\number\x}\dimo
A necessary and sufficient condition for $(\frak{g},\frak{q})$ to be
minimal is that \par\noindent
$\frak{g}_++\frak{k}=\frak{g}$. By complexification
this condition can be rewritten as\,:
\form{\Hat{\frak{g}}=\frak{q}\cap\bar{\frak{q}}+\Hat{\frak{k}}\,.}
\edef\formii{\number\q.\number\t}\edef\lemmaJE{\number\q.\number\x}
Since $\Hat{\frak{h}}\subset\frak{q}\cap\bar{\frak{q}}$, it suffices
to show that $(\formuih)$ is equivalent to\,:
\form{\Hat{\frak{g}}^\alpha\subset\frak{q}\cap\bar{\frak{q}}+\Hat{\frak{k}}
\quad\text{
for all $\alpha\in\Cal{R}$}\,.}\edef\formij{\number\q.\number\t}
\smallskip
To make the equivalence more clear, we first prove\,:
\lem{For all $\alpha\notin\left(\Cal{Q}^n\cap\bar{\Cal{Q}}^{\,-n}\right)
\cup\left(
\Cal{Q}^{-n}\cap\bar{\Cal{Q}}^{\,n}\right)
=\left(\Cal{Q}^n\setminus\bar\Cal{Q}\right)
\cup\left(\bar{\Cal{Q}}^{\, n}\setminus\Cal{Q}\right)$ we have
$$\Hat{\frak{g}}^{\alpha}\subset\frak{q}\cap\bar{\frak{q}}+\Hat{\frak{k}}
\,.$$}\edef\lemmif{\number\q.\number\x} 
\dimo 
We get\,:
$\Cal{R}\setminus\left[\left(\Cal{Q}^n\setminus\bar\Cal{Q}\right)
\cup\left(\bar{\Cal{Q}}^{\, n}\setminus\Cal{Q}\right)\right]=
\Cal{Q}^r\cup\bar{\Cal{Q}}^{\, r}\cup \left(\Cal{Q}^n\cap\bar{\Cal{Q}}\right)
\cup \left(\bar{\Cal{Q}}^{\, n}\cap{\Cal{Q}}\right)\,.$\par
Clearly $\Hat{\frak{g}}^{\alpha}\subset\frak{q}\cap\bar{\frak{q}}
\subset\frak{q}\cap\bar{\frak{q}}+\Hat{\frak{k}}$ when
$\alpha\in\left(\Cal{Q}^n\cap\bar{\Cal{Q}}\right)
\cup \left(\bar{\Cal{Q}}^{\, n}\cap{\Cal{Q}}\right)$.\par\smallskip
Since the statement is invariant if we interchange $\frak{q}$ and
$\bar{\frak{q}}$, to complete the proof of the lemma it suffices
to show that $\Hat{\frak{g}}^{\alpha}\subset
\frak{q}\cap\bar{\frak{q}}+\Hat{\frak{k}}$ for all
$\alpha\in\Cal{Q}^r$. \par   
Let $C\in\frak{C}(\Cal{R},\Cal{Q})$ be V-fit,
so that 
$\sigma\left(\left(\Cal{R}^+(C)\cap\Cal{Q}^r\right)\setminus\rootre\right)
\subset\Cal{R}^-(C)$. We also use the notation
$\Hat{\frak{k}}^{\alpha,-\bar\alpha}=\Hat{\frak{k}}\cap
\left(\Hat{\frak{g}}^{\alpha}+\Hat{\frak{g}}^{-\bar\alpha}\right)$.
This is a one-dimensional complex subspace of $\Hat{\frak{g}}$ when
$\alpha\in\Cal{R}\setminus\rootimp$.
\par \smallskip
For $\alpha\in\Cal{Q}^r\cap\left(\rootre\cup\rootim\right)$, we have
$\alpha\in\Cal{Q}^r\cap{\bar{\Cal{Q}}}^{\, r}
\subset\Cal{Q}\cap{\bar{\Cal{Q}}}$ and hence 
$\Hat{\frak{g}}^{\alpha}\subset\frak{q}\cap\bar{\frak{q}}$.
\par 
If $\alpha\in\Cal{Q}^r\cap\Cal{R}^+(C)\cap\rootcp$, 
then 
$-\bar\alpha\in\Cal{R}^+(C)\subset\Cal{Q}$. Hence 
$-\alpha,-\bar\alpha\in\Cal{Q}\cap\bar{\Cal{Q}}$, so that\,:\quad
$\Hat{\frak{g}}^{-\alpha}\subset\frak{q}\cap\bar{\frak{q}}\subset
\frak{q}\cap\bar{\frak{q}}+\Hat{\frak{k}}\,,\quad
\Hat{\frak{g}}^{\alpha}\subset\Hat{\frak{g}}^{-\bar\alpha}+
\Hat{\frak{k}}^{(\alpha,-\bar\alpha)}\subset\frak{q}\cap\bar{\frak{q}}
+\Hat{\frak{k}}\,.$
\enddimo
We conclude now the proof of Proposition {\PROPie}. 
If condition $(\formuih)$ is satisfied, then
$\left(\Cal{Q}^n\cap\bar{\Cal{Q}}^{\,-n}\right)
\cup\left(\Cal{Q}^{-n}\cap\bar{\Cal{Q}}^{\,n}\right)
\subset\Cal{R}_\bullet$, hence
$\frak{g}^{\alpha}\subset\Hat{\frak{k}}$ for all
$\alpha\in\left(\Cal{Q}^n\cap\bar{\Cal{Q}}^{\,-n}\right)
\cup\left(\Cal{Q}^{-n}\cap\bar{\Cal{Q}}^{\,n}\right)$.
In view of Lemma {\lemmif}, we obtain $(\formij)$ and hence 
$(\formii)$.
\smallskip
Vice versa, assume that there is $\alpha\in\Cal{Q}^n\cap
\bar{\Cal{Q}}^{\,-n}\setminus\rootimc$. If
$\alpha\in\rootimp$, then 
$\Hat{\frak{g}}^{\alpha}\in\Hat{\frak{p}}$, and cannot be
contained in $\frak{q}\cap\bar{\frak{q}}+\Hat{\frak{k}}$.
Otherwise, $\alpha\in\rootcp$ and 
$\alpha,-\bar\alpha\in\Cal{Q}^n\cap
\bar{\Cal{Q}}^{\,-n}$. This implies that
$\Hat{\frak{g}}^{\alpha}\oplus\Hat{\frak{g}}^{\bar\alpha}
\oplus\Hat{\frak{g}}^{-\alpha}\oplus\Hat{\frak{g}}^{-\bar\alpha}
=\Hat{\frak{k}}^{(\alpha,-\bar\alpha)}\oplus
\Hat{\frak{k}}^{(-\alpha,\bar\alpha)}\oplus
\Hat{\frak{p}}^{(\alpha,-\bar\alpha)}\oplus
\Hat{\frak{p}}^{(\alpha,-\bar\alpha)}$
has intersection $\{0\}$ with $\frak{q}\cap\bar{\frak{q}}$
and therefore $(\formii)$ cannot possibly hold true. \enddimo
\subhead \number\q.5. Cross-marked diagrams and examples \endsubhead\par\noindent
In the examples that will follow, here and in the next sections,
to describe specific parabolic $CR$ algebras, 
we shall utilize {\it cross-marked diagrams}. They are Dynkin diagrams,
where the simple roots in $\Cal{B}(C)$ for
a Weyl chamber $C\in\frak{C}(\Cal{R},\Cal{Q})$, are indicated by\,:
$$\align 
\bianca\; &\quad\text{ if the root is real;}\\
\nera\; &\quad\text{ if the root is compact imaginary;}\\
{\ssize\circledast} &\quad\text{ if the root is noncompact imaginary;}\\
{\ssize \oplus}&\quad\text{ if the root is complex and its conjugate belongs to
$\Cal{R}^+(C)$;}\\
{\ssize\ominus}&
\quad\text{ if the root is complex and its conjugate belongs to
$\Cal{R}^-(C)$}
\endalign
$$ 
and we cross-mark the roots in $\Phi_C$. Some extra information about
the action of $\sigma$ on the simple roots in $\Cal{B}(C)$ is 
provided by some arrows and dotted arrows joining pairs of simple roots
that have the same, or opposite, restriction to $\frak{h}^-\subset
\frak{h}_{\Bbb{R}}$, or that are the edges of segments whose
nodes are the support of real or imaginary roots. 
However, as we shall see, the most important information is
carried by the {\it colors} of the nodes.

\esemp{
Consider 
the $CR$ manifold $M$ consisting of $3$-planes
$\ell_3$ of $\Bbb{C}^6$ with $\roman{dim}_{\Bbb{C}}(\ell_3\cap
\bar\ell_3)=1$. This is an orbit of $\bold{SL}(6,\Bbb{R})$ in the
Grassmannian of $3$-planes of $\Bbb{C}^6$. Let
$\epsilon_1,\hdots,\epsilon_6$ be the canonical basis of $\Bbb{C}^6$.
It is convenient to represent the Lie algebra $\frak{sl}(6,\Bbb{R})$
in the basis 
$$e_1=\epsilon_1+i\epsilon_6\,,\;
e_2=\epsilon_2+i\epsilon_5\,,\;
e_3=\epsilon_3\,,\;
e_4=\epsilon_4\,,\;
e_5=\epsilon_2-i\epsilon_5\,,\;
e_6=\epsilon_1-i\epsilon_6\,.
$$
Then $M$ is the orbit of the $3$ plane generated by $e_1,e_2,e_3$. 
\par
We can take $\frak{h}_{\Bbb{R}}$ to be the set of 
real diagonal matrices.
The parabolic $\frak{q}$ is $\frak{q}_A$ for 
$A=\roman{diag}(1,1,1,-1,-1,-1)$.
The Weyl chamber $C$ corresponding to the canonical basis
$\alpha_i=e_i-e_{i+1}$ ($i=1,\hdots,5$) belongs to
$\frak{C}(\Cal{R},\Cal{Q})$ and is V-fit.
We have indeed
$$
\left\{\smallmatrix\\
\bar\alpha_1&=&e_6-e_5=&-\alpha_5\\
\bar\alpha_2&=&e_5-e_3=&-(\alpha_3&+\alpha_4)\\
\bar\alpha_3&=&e_3-e_4=&\alpha_3\\
\bar\alpha_4&=&e_4-e_2=&-(\alpha_2&+\alpha_3)\\
\bar\alpha_5&=&e_2-e_1=&-\alpha_1\endsmallmatrix \right.
$$
so that the associated diagram is\,:
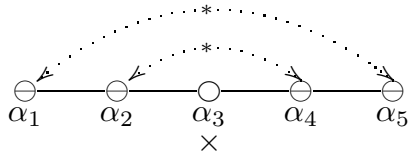
$$
\xymatrix @M=0pt @R=2pt @!C=5pt{
\ominus\ar@{-}[r]
\ar@/^2pc/@{<.>}[rrrr]!U|{\ast}
&\ominus\ar@{-}[r]
\ar@/^1pc/@{<.>}[rr]!U|{\ast}&\bbianca
\ar@{-}[r]
&\ominus
\ar@{-}[r]&\ominus\\
\alpha_1&\alpha_2&\alpha_{3}&\alpha_4&\alpha_5\\
&&\times}
$$
We have $\Phi_C=\{\alpha_3\}$, $\sigma\left([\Cal{R}^+(C)\cap\Cal{Q}^r]
\setminus\rootre\right)\subset\Cal{R}^-(C)$, $\bar{\Phi}_C=\Phi_C\subset
\Cal{R}^+(C)$. Hence, by Lemma {\lemmaIB},
$(\frak{g},\frak{q})$ is weakly nondegenerate.
Our $M$ is a $CR$ manifold of hypersurface type $(8,1)$.
Its Levi form has two positive, two negative and four zero eigenvalues.
Hence $M$ is fundamental and weakly, but not strictly nondegenerate.
\par
We obtain an S-fit chamber by describing $\frak{sl}(6,\Bbb{R})$ in the
basis\,:
$$e_1=\epsilon_1\,,\;e_2=\epsilon_2+i\epsilon_5\,,\;
e_3=\epsilon_3+i\epsilon_4\,,\; e_4=\epsilon_3-i\epsilon_4\,,\;
e_5=\epsilon_2-i\epsilon_5\,,\;e_6=\epsilon_6\,.$$
One verifies that, with the basis of simple roots $\alpha_i=e_i-e_{i+1}$ ($1\leq
i\leq 5$) the corresponding diagram is\,:
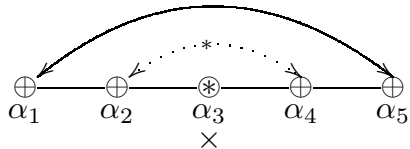
$$
\xymatrix @M=0pt @R=2pt @!C=5pt{
\oplus\ar@{-}[r]
\ar@/^2pc/@{<->}[rrrr]!U
&\oplus\ar@{-}[r]
\ar@/^1pc/@{<.>}[rr]!U|{\ast}&\grigia
\ar@{-}[r]
&\oplus
\ar@{-}[r]&\oplus\\
\alpha_1&\alpha_2&\alpha_{3}&\alpha_4&\alpha_5\\
&&\times}
$$
We have $\bar\alpha_1=\alpha_1+\alpha_2+\alpha_3+\alpha_4\succ_C\alpha_3$,
in accordance with the fact that $(\frak{g},\frak{q})$ is fundamental.}
\par
\esemp{Consider the $CR$ manifold $M^*$ consisting of all flags
$\ell_1\subset\ell_2\subset\ell_4\subset\ell_5\subset\Bbb{C}^6$
(the subscript is the dimension of the linear subspace), with
$$
\roman{dim}_{\Bbb{C}}(\ell_1\cap\bar\ell_1)=0\,,\quad
\ell_2=\ell_1+\bar{\ell}_1\,,\quad
\roman{dim}_{\Bbb{C}}(\ell_4\cap\bar\ell_4)=3\,,\quad
\roman{dim}_{\Bbb{C}}(\ell_5\cap\bar\ell_5)=3\, .
$$
We note that $M^*$ is not connected. Thus an orbit of $\bold{SL}(6,\Bbb{R})$
in $M^*$ will be
a connected component $M$ of $M^*$. We can better describe such an
$M$ in terms of the choice of a suitable basis of $\Bbb{C}^6$.
Denoting by $\epsilon_1,\hdots,\epsilon_6$ 
the canonical basis of 
$\Bbb{C}^6$, we introduce the basis\,:
$$e_1=\epsilon_1+i\epsilon_2,\;
e_2=\epsilon_1-i\epsilon_2,\;
e_3=\epsilon_3,\;
e_4=\epsilon_4+i\epsilon_6,\;
e_5=\epsilon_5,\;
e_6=\epsilon_4-i\epsilon_6\,.$$ 
Our $M$ is the orbit, under the action of $\bold{SL}(6,\Bbb{R})$, 
of the flag
$$\langle e_1\rangle \subset
\langle e_1,e_2\rangle 
\subset
\langle e_1,e_2,e_3,e_4\rangle\subset
\langle e_1,e_2,e_3,e_4,e_5\rangle\,.$$
We can take the Cartan subalgebra $\frak{h}$ of $\frak{sl}(6,\Bbb{R})$
that is described, in the basis $e_1,\hdots,e_6$, by\,:
$$\frak{h}=\{\roman{diag}(\lambda,\bar\lambda,a,
\mu,b,\bar\mu)\,|\, \lambda,\mu\in\Bbb{C},\;a,b\in\Bbb{R}\,,\;
2\roman{Re}(\lambda+\mu)+a+b=0\}$$
and consider the corresponding root system
$\Cal{R}$ of $\frak{sl}(6,\Bbb{C})$ with respect to
$\Hat{\frak{h}}$. We identify $\frak{h}_{\Bbb{R}}$ with the space
of real diagonal matrices in $\frak{sl}(6,\Bbb{C})$ and
the $e_h$'s to the evaluation of the $h$-th diagonal entry of
$H\in\frak{h}_{\Bbb{R}}$. Let $(\frak{g},\frak{q})$ be the
effective parabolic $CR$ algebra associated to $M$ and
$\Cal{Q}$ the parabolic set of $\frak{q}$. 
Then the Weyl chamber $C\in\frak{C}(\Cal{R})$
with basis $\Cal{B}(C)=\{\alpha_i=e_i-e_{i+1}\,|\, 1\leq i\leq 5\}$
belongs to $\frak{C}(\Cal{R},\Cal{Q})$ and
$\Phi_C=\{\alpha_1,\alpha_2,\alpha_4,\alpha_5\}$.
We have\,:
$$\left\{\smallmatrix
\bar\alpha_1&=&\!\!-\alpha_1\;\\
\bar\alpha_2&=&\alpha_1+\alpha_2\\
\bar\alpha_3&=&\alpha_3+\alpha_4&\!\!+\alpha_5\\
\bar\alpha_4&=&\!\!-\alpha_5\;\\
\bar\alpha_5&=&\!\!-\alpha_4\,.\,\endsmallmatrix\right.
$$
Thus ($\FFrmEJ)$ yields in this case\,:  
$$\cases
\Upsilon_C^{(0)}=\{\alpha_2\}\\
\Upsilon_C^{(1)}=\Upsilon_C^{(0)}=\Upsilon_C\,.\endcases
$$
Thus the basis of the fundamental reduction is the Grassmannian of
$2$-planes $\ell_2$ in $\Bbb{C}^6$ with $\bar{\ell}_2={\ell}_2$.
We give below the diagrams for $(\frak{g},\frak{q})$,
for its fundamental reduction $(\frak{g},\frak{q}')$,
and for the fiber $(\frak{g},\frak{q}'')$, (recall that the basis of
the fundamental reduction of a parabolic is still
parabolic).
$$(\frak{g},\frak{q})\qquad
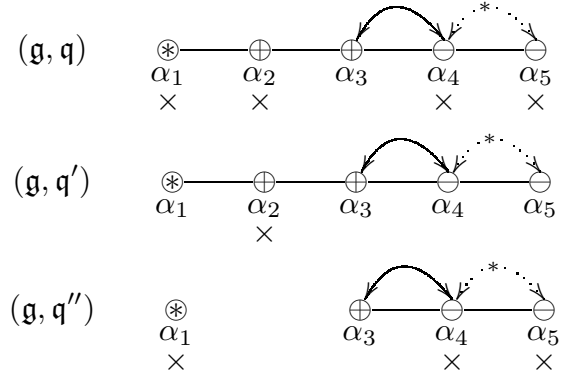
\xymatrix @M=0pt @R=2pt @!C=5pt{
\grigia\ar@{-}[r]
&\oplus\ar@{-}[r]&\oplus\ar@{-}[r]
\ar@/^1pc/@{<->}[r]!U&\ominus\ar@{-}[r]
\ar@/^1pc/@{<.>}[r]!U|{\ast}&\ominus\\
\alpha_1&\alpha_2&\alpha_{3}&\alpha_4&\alpha_5\\
\times&\times&&\times&\times}
$$
$$(\frak{g},\frak{q}')\qquad
\xymatrix @M=0pt @R=2pt @!C=5pt{
\grigia\ar@{-}[r]
&\oplus\ar@{-}[r]&\oplus\ar@{-}[r]\ar@/^1pc/@{<->}[r]!U&\ominus\ar@{-}[r]
\ar@/^1pc/@{<.>}[r]!U|{\ast}&
\ominus\\
\alpha_1&\alpha_2&\alpha_{3}&\alpha_4&\alpha_5\\
&\times}
$$
$$(\frak{g},\frak{q}'')\qquad
\xymatrix @M=0pt @R=2pt @!C=5pt{
\grigia
&\qquad&\oplus\ar@{-}[r]\ar@/^1pc/@{<->}[r]!U
&\ominus\ar@{-}[r]\ar@/^1pc/@{<.>}[r]!U|{\ast}
&\ominus\\
\alpha_1&&\alpha_{3}&\alpha_4&\alpha_5\\
\times&&&\times&\times}
$$
The fiber $M''$ is the product of a complex disk (a connected component
of the set of $\Bbb{CP}^1\setminus\Bbb{RP}^1\simeq S^2\setminus S^1$)
and a connected
$CR$ manifold $N$, consisting of flags in $\Bbb{C}^6/\Bbb{C}^2$:
these
can be identified to pairs $\ell_2\subset\ell_3\subset
\Bbb{C}^4$ such that $\roman{dim}_{\Bbb{C}}(\ell_2\cap\bar\ell_2)=1$,
$\ell_3\not\supset\ell_2+\bar{\ell}_2$
and  $\roman{dim}_{\Bbb{C}}(\ell_3\cap\bar\ell_3)=2$.
It is convenient to utilize the basis $e_4,e_3,e_5,e_6$ of
$\Bbb{C}^4\simeq\langle e_3,e_4,e_5,e_6\rangle\subset\Bbb{C}^6$.
Then $\beta_1=e_4-e_3,\beta_2=e_3-e_5,\beta_3=e_5-e_6$ is the
basis related to a Weyl chamber $C_N$ in which the diagram associated
to the parabolic $CR$ algebra of $N$ is\,:
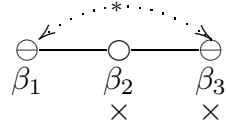
$$\xymatrix @M=0pt @R=2pt @!C=5pt{
\ominus\ar@{-}[r]\ar@/^1pc/@{<.>}[rr]!U|{\ast}
&\bbianca\ar@{-}[r]&\ominus\\
\beta_1&\beta_{2}&\beta_3\\
&\times&\times}
$$
Since $C_N$ is V-fit,
we see from this diagram that
$N$ is weakly degenerate. The basis $N'$ of its
weakly nondegenerate reduction,
consists of planes $\ell_2$ of $\Bbb{C}^4$
with $\roman{dim}_{\Bbb{C}}(\ell_2\cap\bar{\ell}_2)=1$, and
corresponds to the diagram\,:
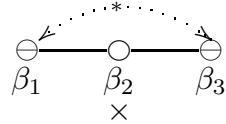
$$\xymatrix @M=0pt @R=2pt @!C=5pt{
\ominus\ar@{-}[r]\ar@/^1pc/@{<.>}[rr]!U|{\ast}
&\bbianca\ar@{-}[r]&\ominus\\
\beta_1&\beta_{2}&\beta_3\\
&\times}
$$
The parabolic $CR$ manifold $N'$ is of hypersurface type $(3,1)$, with
a degenerate Levi form of signature $(1,-1,0)$. 
Thus it is $1$-pseudoconcave (see e.g. [HN])
and weakly, but not strictly, nondegenerate.
The fiber $F$ of the $\bold{SL}(4,\Bbb{R})$-equivariant
weakly nondegenerate reduction
$N@>>>N'$, that lies above
a given $2$-plane $\ell_2$ with $\roman{dim}_{\Bbb{C}}(\ell_2+
\bar{\ell}_2)=3$, 
is isomorphic to the pencil of $3$-planes in $\Bbb{C}^4$, that
contain $\ell_2$ and are distinct from $(\ell_2+\bar{\ell}_2)$.
Thus 
$F\simeq \Bbb{CP}^1\setminus\{\text{a point}\}\simeq\Bbb{C}$.
Note that the $CR$ algebra that is naturally associated to the
fiber fails in this case to be  parabolic. In fact, the $CR$ algebra
$(\frak{g}^\sharp,\frak{q}^\sharp)$ of the fiber is given by\,:
$$\frak{g}^\sharp=\left.\left\{
\pmatrix
\lambda&0&\bar{\zeta}&0\\
z&h&s&\bar{z}\\
0&0&k&0\\
0&0&\zeta&\bar{\lambda}
\endpmatrix\right|{\matrix\format\l\\
\lambda,z,\zeta\in\Bbb{C}\\
h,k,s\in\Bbb{R}\\
h+k+2\roman{Re}\lambda=0
\endmatrix}\right\}\,,
$$ 
$$
\frak{q}^\sharp
=\left.\left\{
\pmatrix
\lambda_1&0&{\zeta_2}&0\\
0&\lambda_2&\theta&{z}_2\\
0&0&\lambda_3&0\\
0&0&\zeta_1&{\lambda}_4
\endpmatrix\right|{\matrix\format\l\\
\lambda_i,z_i,\zeta_i,\theta\in\Bbb{C}\\\vspace{2\jot}
\sum_{i=1}^4{\lambda_i}=0
\endmatrix}\right\}\,,
$$
where both
$\frak{g}^\sharp$ and its complexification
$\Hat{\frak{g}}^{\sharp}$ are nilpotent.}


\se{Canonical fibrations over a parabolic $CR$ manifold}
Given an effective parabolic $CR$ algebra $(\frak{g},\frak{q})$,
we constructed in the previous sections new parabolic
complex subalgebras $\frak{q}'$ of $\Hat{\frak{g}}$ with
$\frak{q}\subset\frak{q}'$, to obtain smooth fibrations 
$M(\frak{g},\frak{q})@>>>M(\frak{g},\frak{q}')$ , namely
with a weakly nondegenerate basis and totally complex fibers,
and with totally real basis and fundamental
fibers. Here, and in the next section,  we consider
smooth fibrations $M(\frak{g},\frak{q}')@>>>M(\frak{g},\frak{q})$,
obtained by
choosing special parabolic $\frak{q}'\subset\Hat{\frak{g}}$ with
$\frak{q}'\subset\frak{q}$, and that will be useful 
to find suitable Weyl chambers in $\frak{C}(\Cal{R},\Cal{Q})$ 
and to investigate
the topology of the general $M(\frak{g},\frak{q})$.\par
We keep the notation of the preceding sections. In particular,
we fix a Cartan pair $(\vartheta,\frak{h})$, assuming that
it is adapted to all the parabolic $CR$ algebras that we shall consider.
\par We have\,:
\prop{Let $(\vartheta,\frak{h})$ be a Cartan pair adapted to
the effective parabolic $CR$ algebras
$(\frak{g},\frak{q})$ and $(\frak{g},\frak{e})$. 
Then\,:
\form{\frak{l}=\frak{q}^n\oplus\left(\frak{q}^r\cap\frak{e}\right)}
is a parabolic complex Lie subalgebra of $\Hat{\frak{g}}$ with\,:
\form{\left\{
\matrix\format\r&\l\\
\frak{h}\;&\subset\frak{l}\\
\frak{l}^n&=\frak{q}^n\oplus\left(\frak{q}^r\cap\frak{e}^n\right)
\supset\frak{q}^n\\
\frak{l}^r&=\frak{q}^r\cap\frak{e}^r\subset\frak{q}^r\\
\frak{l}\;\,&=\frak{l}^n\oplus\frak{l}^r\subset\frak{q}\,.
\endmatrix\right.}}\edef\formulaiD{\number\q.\number\t}
\advance\t by-1
\edef\formulaiC{\number\q.\number\t}
\advance\t by1 \edef\propoIB{\number\q.\number\x}
\dimo Let $\Cal{Q}$, $\Cal{E}$ be the parabolic sets in $\Cal{R}$
corresponding to the complex parabolic Lie subalgebras $\frak{q}$,
$\frak{e}$, respectively. To prove that $\frak{l}$ is parabolic, we
need to prove that
\form{\Cal{L}=\Cal{Q}^n\cup\left(\Cal{Q}^r\cap\Cal{E}\right)
=\Cal{Q}^n\cup\left(\Cal{Q}^r\cap\Cal{E}^r\right)\cup
\left(\Cal{Q}^r\cap\Cal{E}^n\right)}
is a parabolic subset of $\Cal{R}$.\par
Let
$A,B\in\frak{h}_{\Bbb{R}}$ be such that
$\Cal{Q}=\Cal{Q}_A$, $\Cal{E}=\Cal{Q}_B$ and fix $\epsilon>0$ so small
that $\epsilon\,|\alpha(B)|<\alpha(A)$ for all $\alpha\in\Cal{Q}^n$.
Then we claim that 
$$\Cal{L}=\{\alpha\in\Cal{R}\,|\, \alpha(A+\epsilon B)\geq 0\}\,.$$
Indeed, when $\alpha\in\Cal{Q}^n$, then
$\alpha(A+\epsilon B)\geq \alpha(A)-\epsilon\,|\alpha(B)|>0$;
when $\alpha\notin\Cal{Q}$, then $-\alpha\in\Cal{Q}^n$ and hence
$\alpha(A+\epsilon B)<0$; finally for $\alpha\in\Cal{Q}^r$, we have
$\alpha(A+\epsilon B)=\epsilon\,\alpha(B)$ and hence
$\alpha\in\Cal{L}$ if and only if $\alpha\in\Cal{E}$.
\par
The proof of $(\formulaiD)$ is straightforward.\enddimo
\medskip
Vice versa, when $\frak{l}$ is a complex parabolic subalgebra 
of $\Hat{\frak{g}}$ with
$\frak{h}\subset\frak{l}\subset\frak{q}$, then
$\frak{l}=\frak{q}^n\oplus\left(\frak{q}^r\cap\frak{l}\right)$,
so that $(\formulaiC)$ gives a way to construct all complex parabolic 
subalgebras of $\Hat{\frak{g}}$ with
$\frak{h}\subset\frak{l}\subset\frak{q}$.\par\medskip

\prop{Let $(\vartheta,\frak{h})$ be a Cartan pair adapted to the
effective parabolic $CR$ algebra
$(\frak{g},\frak{q})$\,.
Consider the complex parabolic Lie subalgebra 
\form{\frak{w}= \frak{q}^n\oplus \left(\frak{q}^r\cap\vartheta(\frak{q})
\right)}
of $\Hat{\frak{g}}$. The
$\sigma$-invariant reductive subalgebra 
$\frak{q}^r\cap\bar{\frak{q}}^{\, r}$ is a complement 
in $\frak{w}$
of its nilradical
$\frak{w}^n$. We have\,: 
\form{\cases
\frak{w}^r=\frak{q}^r\cap\bar{\frak{q}}^{\, r}\,,\quad
\frak{w}^n=\frak{q}^n\oplus\left(\frak{q}^r\cap\vartheta(\frak{q}^n)
\right)\supset\frak{q}^n\,,\quad
\frak{w}=\frak{w}^n\oplus\frak{w}^r\subset\frak{q}\,,\\\vspace{2\jot}
\frak{w}^n\cap\bar\frak{w}^{\, n}=\frak{q}^n\cap\bar{\frak{q}}^{\, 
n}\,,\quad \frak{w}^r=\bar{\frak{w}}^{\, r}\,,\quad
\frak{w}\cap\bar\frak{w}=\frak{q}^r\cap\bar{\frak{q}}^{\, r}
\,\oplus\,\frak{q}^n\cap\bar{\frak{q}}^{\, n}\,,\endcases}
and $\frak{w}$ 
is the smallest parabolic subalgebra
of $\Hat{\frak{g}}$ that satisfies the conditions\,:
\form{\frak{q}^r\cap\bar{\frak{q}}^r\subset
\frak{w}\subset\frak{q}
\quad
\text{and}\quad 
\frak{w}\,+\,\bar{\frak{w}}\,=\,\frak{q}\,+\,\bar{\frak{q}}\,.}}
\edef\formulaJC{\number\q.\number\t}
\advance\t by-1
\edef\formulaJB{\number\q.\number\t}
\advance\t by-1

\advance\t by2 \edef\proposizioneJZ{\number\q.\number\x}
\dimo By Proposition {\propoIB}, $\frak{w}$ is complex parabolic in
$\Hat{\frak{g}}$.
Indeed $\vartheta(\frak{q})$ is complex parabolic in
$\Hat{\frak{g}}$ and contains $\frak{h}$.
The parabolic set associated to $\vartheta(\frak{q})$ is
$\vartheta(\Cal{Q})=\{\alpha\,|\,-\bar\alpha\in\Cal{Q}\}=\bar{\Cal{Q}}^r\cup
\bar{\Cal{Q}}^{\,-n}$.
Hence the parabolic set corresponding to $\frak{w}$ is\,:
\form{\Cal{W}=\Cal{Q}^n\cup\left(\Cal{Q}^r\cap\vartheta(\Cal{Q})\right)=
\left(\Cal{Q}^r\cap\bar{\Cal{Q}}^{\, r}\right)\cup
\Cal{Q}^n\cup\left(\Cal{Q}^r\cap\bar{\Cal{Q}}^{\,-n}\right)\,.}
We obtain 
$(\formulaJB)$ by using Proposition {\propoIB}.\par
We have 
$\Cal{W}^r=\bar{\Cal{W}}^{\,r}=\Cal{Q}^r\cap\bar{\Cal{Q}}^{\,r}$, 
and
$\Cal{W}\cup\bar{\Cal{W}}=\Cal{W}^r\cup\Cal{W}^n
\cup\bar{\Cal{W}}^{\, n}=
\Cal{Q}\cup\bar{\Cal{Q}}$. The right hand side of this
equality 
can be written as a disjoint union\,:
$$\Cal{Q}\cup\bar{\Cal{Q}}=\left(\Cal{Q}^r\cap\bar{\Cal{Q}}^{\,r}\right)
\cup\left(\Cal{Q}^n\cup\bar{\Cal{Q}}^{\,n}\right)
\cup\left(\Cal{Q}^r\setminus
\bar{\Cal{Q}}\right)\cup \left(\bar{\Cal{Q}}^{\,r}\setminus
{\Cal{Q}}\right)\,.$$
In particular, $\Cal{Q}^r\setminus{\bar{\Cal{Q}}}\subset\Cal{W}$
for the parabolic set $\Cal{W}$ of any complex parabolic $\frak{w}$
that satisfies
$(\formulaJC)$,
and this shows that the
$\frak{w}$ we constructed 
is the smallest complex parabolic subalgebra of $\Hat{\frak{g}}$
that satisfies
$(\formulaJC)$.\enddimo
Since $\frak{w}\subset\frak{q}$, we have a 
$\frak{g}$-equivariant $CR$ fibration
$(\frak{g},\frak{w})@>>>(\frak{g},\frak{q})$.
We call $(\frak{g},\frak{w})$ the {\it canonical $CR$-lift} of
$(\frak{g},\frak{q})$. \par
\thm{The canonical $CR$ lift 
$(\frak{g},\frak{w})@>>>(\frak{g},\frak{q})$ is a 
$\frak{g}$-equivariant $CR$ fibration
(in particular a $CR$ submersion) with totally
complex fibers. When $(\frak{g},\frak{q})$ is
weakly nondegenerate, then 
$(\frak{g},\frak{w})@>>>(\frak{g},\frak{q})$ is the weakly
nondegenerate reduction 
of $(\frak{g},\frak{w})$ (cf. {\rm [MN4]}).}
\dimo The statement is a consequence of Lemma {\lemmazzo} and of the
next lemma.\enddimo
\lem{Let $(\vartheta,\frak{h})$ be a Cartan pair adapted to the
effective parabolic $CR$ algebra $(\frak{g},\frak{q})$.
Let $(\frak{g},\frak{w})$ be its canonical $CR$-lift. 
Then a Weyl chamber $C\in\frak{C}(\Cal{R},\Cal{Q})$ is V-fit for
$(\frak{g},\frak{q})$ if and only if $C\in\frak{C}(\Cal{R},\Cal{W})$
and is V-fit for $(\frak{g},\frak{w})$.}

\dimo Assume that $C\in\frak{C}(\Cal{R},\Cal{Q})$ is V-fit for
$(\frak{g},\frak{q})$. We want to show that $\Cal{R}^+(C)\subset\Cal{W}$.
Since $\Cal{Q}^n\cup\left(\Cal{Q}^r\cap
\bar{\Cal{Q}}^{\,r}\right)\subset\Cal{W}$, it suffices to prove that
$\alpha\in\Cal{W}$ for 
$\alpha\in\left(\Cal{R}^+(C)\cap\Cal{Q}^r\right)
\setminus\bar{\Cal{Q}}^{\,r}$.
Since $C$ is V-fit for $(\frak{g},\frak{q})$, 
for such a root $\alpha$ we have 
$\bar\alpha\prec_C0$. Hence $\vartheta(\alpha)=-\bar\alpha\in\Cal{R}^+(C)
\subset\Cal{Q}$, i.e. $\alpha\in\Cal{Q}^r\cap\vartheta(\Cal{Q})\subset
\Cal{W}$. A chamber $C\in\frak{C}(\Cal{R},\Cal{W})$ that is V-fit
for $(\frak{g},\frak{q})$ is also V-fit for $(\frak{g},\frak{w})$,
because $\frak{w}\subset\frak{q}$.
\par
Let now
$C\in\frak{C}(\Cal{R},\Cal{W})
\subset\frak{C}(\Cal{R},\Cal{Q})$ be V-fit for $(\frak{g},\frak{w})$.
Note that $\Cal{R}^+(C)\cap\Cal{Q}^r=\left(\Cal{R}^+(C)\cap\Cal{W}^r\right)
\cup\left(\Cal{W}^n\setminus\Cal{Q}^n\right)$, and
$\Cal{W}^n\setminus\Cal{Q}^n=\Cal{Q}^r\cap\bar{\Cal{Q}}^{\,-n}$.
Since $\bar\alpha\prec_C0$ for all $\alpha\in\bar{\Cal{Q}}^{\,-n}$,
we obtain that $\sigma(\Cal{R}^+(C)\cap\Cal{Q}^r\setminus\rootre)\subset
\Cal{R}^-(C)$ and therefore $C$ is also V-fit for $(\frak{g},\frak{q})$.
\enddimo
\edef\LemmaHC{\number\q.\number\x}
\prop{Let $(\vartheta,\frak{h})$ be a Cartan pair adapted to the
effective parabolic $CR$ algebra $(\frak{g},\frak{q})$, and
let $(\frak{g},\frak{w})$ be its canonical lift.
Let $C\in\frak{C}(\Cal{R},\Cal{Q})$ be V-fit
for $(\frak{g},\frak{q})$, and $\Phi_C=\Cal{B}(C)\cap\Cal{Q}^n$,
$\tilde{\Phi}_C=\Cal{B}(C)\cap\Cal{W}^n$. Then\,:
\form{
\tilde{\Phi}_C=\Phi_C\cup\{\alpha\in\Cal{B}_C\,|\,\supp{C}(\bar\alpha)
\cap\Phi_C\neq\emptyset\}\,.}}
\dimo By 
Lemma {\LemmaHC},
$\frak{w}=\frak{q}_{\tilde{\Phi}_C}$ for some set of simple roots
$\tilde{\Phi}_C$ with $\Phi_C\subset\tilde{\Phi}_C\subset\Cal{B}(C)$.
Let $\Cal{Q}=\Cal{Q}_A=\{\alpha\in\Cal{R}\,|\,\alpha(A)\geq 0\}$, with
$A\in\frak{h}_{\Bbb{R}}$.
Fix a real $\epsilon>0$ with
$\epsilon\,|\bar\alpha(A)|<\alpha(A)$ for all $\alpha\in\Cal{Q}^n$.
Then\,:
\form{\Cal{W}=\{\alpha\in\Cal{R}\,|\,
\alpha(A-\epsilon\bar{A})\geq 0\,\}\,.}
Indeed, $\Cal{W}\subset\Cal{Q}$ because $\alpha(A-\epsilon\bar{A})<0$
when $\alpha\in\Cal{R}$ and $\alpha(A)<0$; moreover
$\Cal{Q}^n\subset\Cal{W}^n$, and a root $\alpha\in\Cal{Q}^r$
belongs to $\Cal{W}$ if and only if $\bar\alpha(A)\leq 0$, i.e. 
if and only if $\vartheta(\alpha)\in\Cal{Q}$. \par
Thus we have\,:
$$\align
\tilde{\Phi}_C&=\Cal{B}(C)\cap\Cal{W}^n=\{\alpha\in\Cal{B}(C)\,|\,
\alpha(A-\epsilon\bar{A})>0\}\\
&=\Phi_C\cup\{\alpha\in\Cal{B}(C)\,|\, \alpha(A)=0\,,\;
\bar\alpha(A)<0\}\\
&=\Phi_C\cup\{\alpha\in\Cal{B}(C)\,|\, \supp{C}(\bar\alpha)\cap\Phi_C\neq
\emptyset\,\}\,,
\endalign
$$
because $\tilde{\Phi}_C\setminus\Phi_C\subset\rootcp$ 
and, for a complex $\alpha$ in $\Cal{B}(C)\setminus\Phi_C$ we
have $\bar\alpha\in\Cal{R}^-(C)$\,:
hence $\bar\alpha(A)<0$
whenever $\bar\alpha(A)\neq 0$, i.e.
$\supp{C}(\bar\alpha)\cap\Phi_C\neq\emptyset$.
\enddimo\edef\PropHD{\number\q.\number\x}
We can slightly improve the criterion of weak non-degeneracy 
of Theorem {\TEOEI}, by using Weyl chambers adapted to the canonical
$CR$ lift.
We have
indeed\,:
\prop{Let $(\vartheta,\frak{h})$ be a Cartan pair adapted to the
effective parabolic $CR$ algebra $(\frak{g},\frak{q})$,
and let $(\frak{g},\frak{w})$ be
the canonical $CR$-lift of $(\frak{g},\frak{q})$.
If $C\in\frak{C}(\Cal{R},\Cal{W})$\,, then\,:
\form{(\frak{g},\frak{q})\quad\text{is weakly non-degenerate if and
only if}\quad \bar\alpha\succ_C0\quad\forall\alpha\in\Phi_C\,.}}
\dimo{} $(\Rightarrow)$\quad We argue by contradiction.
Assume that $\bar\alpha_0\prec_C0$ for some $\alpha_0\in\Phi_C$. 
We want to prove that $\Cal{Q}'\cup\bar{\Cal{Q}}\,'=\Cal{Q}\cup\bar{\Cal{Q}}$
for the parabolic set 
\form{\Cal{Q}'=\Cal{Q}_{\Psi_C}
=\Cal{Q}\cup\left\{\beta\in\Cal{R}\,|\,\supp{C}(\beta)\cap\Phi_C
\subset\{\alpha_0\}\right\}\,,}\edef\formulaKG{\number\q.\number\t}
corresponding to 
$\Psi_C=\Phi_C\setminus\{\alpha_0\}\subset\Cal{B}(C)$.
It suffices to verify that $\bar\beta\in\Cal{Q}$ if
$\beta\prec_C0$ and $\supp{C}(\beta)\cap\Phi_C=\{\alpha_0\}$.
Since
$C\in\frak{C}(\Cal{R},\Cal{W})$, 
each $\alpha\in\Cal{B}(C)$
either belongs to $\Cal{W}^n=
\Cal{Q}^n\cup\left(\Cal{Q}^r\cap\bar{\Cal{Q}}^{\,-n}\right)$, or to
$\Cal{W}^r=\Cal{Q}^r\cap\bar{\Cal{Q}}^{\,r}$. 
Thus $\Cal{B}(C)\setminus\Phi_C\subset\bar{\Cal{Q}}^{\,r}
\cup\bar{\Cal{Q}}^{\,-n}$, and, being $\bar\alpha_0\prec_C0$, we get
$-\bar\alpha\in\Cal{Q}$ for all $\alpha\in\supp{C}(\beta)$,
yielding
$\beta\in\bar{\Cal{Q}}$.
\smallskip
\noindent $(\Leftarrow)$ 
Assume that $(\frak{g},\frak{q})$ is weakly degenerate. Then, 
for some $\alpha_0\in\Phi_C$, $(\formulaKG)$ defines a parabolic set
$\Cal{Q}'$ with
$\Cal{Q}'\cup\bar{\Cal{Q}}\,'=\Cal{Q}\cup\bar{\Cal{Q}}$. 
In particular, $-\alpha_0\in\bar{\Cal{Q}}$. If 
$-\alpha_0\in\bar{\Cal{Q}}^n$, then $-\bar\alpha_0\in\Cal{Q}^n\subset
\Cal{R}^+(C)$ and $\bar\alpha_0\prec_C0$. Otherwise,
$\bar\alpha_0\in{\Cal{Q}}^{\,r}\cap\bar{\Cal{Q}}^n$ and, because
${\Cal{Q}}^{\,r}\cap\bar{\Cal{Q}}^n\subset\Cal{W}^{-n}$,
the condition that
$C\in\frak{C}(\Cal{R},\Cal{W})$ implies that $\bar\alpha_0\prec_C 0$.
\enddimo\edef\PRPFF{\number\q.\number\x}
Proposition {\PRPFF} gives a way to construct the weakly non-degenerate
reduction of a parabolic $CR$ algebra\,:
\cor{Let $(\frak{g},\frak{w})$ be the canonical $CR$-lift of the
parabolic $CR$ algebra $(\frak{g},\frak{q})$ and $\frak{h}$ 
a Cartan subalgebra of $\frak{g}$ adapted to $(\frak{g},\frak{q})$.
If $C\in\frak{C}(\Cal{R},\Cal{W})\subset\frak{C}(\Cal{R},\Cal{Q})$,
$\Phi_C=\Cal{B}(C)\cap\Cal{Q}^n$, and 
\form{\Psi_C=\{\alpha\in\Phi(C)\,|\, \bar\alpha\in\Cal{R}^+(C)\}\,,}
then the  $\frak{g}$-equivariant $CR$ fibration
$(\frak{g},\frak{q})@>>>(\frak{g},\frak{q}_{\Psi_C})$ 
is the weakly nondegenerate
reduction of $(\frak{g},\frak{q})$.\qed}
\prop{
Let $(\vartheta,\frak{h})$ be a Cartan pair adapted to the
effective parabolic $CR$ algebra $(\frak{g},\frak{q})$, and 
let $(\frak{g},\frak{w})$ be the canonical lift of
$(\frak{g},\frak{q})$.
Then the natural $\bold{G}$-equivariant projection
$M(\frak{g},\frak{w})@>>>M(\frak{g},\frak{q})$ is 
a $CR$ fibration with totally complex connected fibers.}

\dimo Our $CR$ manifolds are described as homogeneous spaces by the
quotients
$M(\frak{g},\frak{q})=\bold{G}/\bold{G}_+$
with $\bold{G}_+=\bold{N}_{\bold{G}}(\frak{q})$ and
$M(\frak{g},\frak{w})=\bold{G}/\bold{W}_+$, where
$\bold{W}_+=\bold{N}_{\bold{G}}(\frak{w})$. \par
We want to prove that
every connected component of $\bold{G}_+$ contains 
an element of $\bold{W}_+$. \par
Let $\frak{g}_+=\frak{g}\cap\frak{q}$ be the Lie algebra
of $\bold{G}_+$ and
consider its decomposition in $(\formulaAA)$ of
Proposition {\proposizioneAE}\,: we have
$\frak{g}_+=\frak{n}\oplus\frak{g}_0$, where $\frak{n}$ is the
ideal of the nilpotent elements of the radical of $\frak{g}_+$
and $\frak{g}_0=\frak{q}^r\cap\bar{\frak{q}}\,^r\cap\frak{g}$ 
is a $\vartheta$-invariant reductive complement
of $\frak{n}$ in $\frak{g}_+$.\par
Being algebraic, $\bold{G}_+$ has a {\it Chevalley decomposition}
(see [Ch, Chap.5, Sect.4]) into the semidirect product 
$\bold{N}\rtimes \bold{G}_0$ of
the analytic subgroup with Lie algebra $\frak{n}$ and 
of a closed Lie subgroup $\bold{G}_0$ with Lie algebra
$\frak{g}_0$.\par
Let $g\in\bold{G}_+$ and denote by $\Gamma_g$ the
connected component of $g$ in $\bold{G}_+$. Since $\bold{N}$ is
connected, we can as well take from the start $g$ in
$\bold{G}_0$, so that in particular
$\roman{Ad}_{\frak{g}}(g)(\frak{g}_0)=\frak{g}_0$.
\par
Since $\frak{g}_0$ is a real form of $\frak{q}^r\cap\bar\frak{q}\,^r$,
by complexification we obtain
that $\roman{Ad}_{\hat{\frak{g}}}(g)(\frak{q}^r\cap\bar\frak{q}\,^r)=
\frak{q}^r\cap\bar\frak{q}\,^r$. Thus
$\roman{Ad}_{\hat{\frak{g}}}(g)(\frak{w})$ is a parabolic
complex subalgebra of $\Hat{\frak{g}}$ with
$\frak{q}^r\cap\bar\frak{q}\,^r\subset
\roman{Ad}_{\hat{\frak{g}}}(g)(\frak{w})\subset\frak{q}$.
Since $\roman{dim}_{\Bbb{C}}(\roman{Ad}_{\hat{\frak{g}}}(g)(\frak{w}))
=\roman{dim}_{\Bbb{C}}(\frak{w})$, it follows 
from the characterization of $\frak{w}$ in  Proposition {\proposizioneJZ},
that
$\roman{Ad}_{\hat{\frak{g}}}(g)(\frak{w})
=\frak{w}$,
and hence $g\in\bold{W}_+$. \enddimo
\par
Next  we describe a construction that is similar to
the one discussed above. We keep the notation 
introduced therein.
\prop{Let $(\vartheta,\frak{h})$ be a Cartan pair adapted to the
effective parabolic $CR$ algebra $(\frak{g},\frak{q})$, and let\,:
\form{\frak{v}=\frak{q}^n\oplus\left(\frak{q}^r\cap\bar{\frak{q}}\right)\,.}
Then $\frak{v}$ is a parabolic subalgebra of $\Hat{\frak{g}}$, such that\,:
\form{\left\{\matrix\format\l\\
\frak{v}^r=\bar{\frak{v}}^{\, r}=\frak{q}^r\cap\bar\frak{q}^{\,r}\,,
\quad
\frak{v}^n=\frak{q}^n\oplus\left(\frak{q}^r\cap\bar{\frak{q}}^{\,n}
\right)
\supset\frak{q}^n\,,
\quad
\frak{v}=\frak{v}^r\oplus\frak{v}^n\subset\frak{q}\\\vspace{2\jot}
\frak{v}^n\cap\bar{\frak{v}}^n=\frak{q}^n\cap\bar\frak{q}^{\,n}\,,
\quad
\bar\frak{v}^{\,r}=\frak{v}^r\,,\quad
\frak{v}\cap\bar{\frak{v}}=\frak{q}\cap\bar{\frak{q}}\,.
\endmatrix\right.}
It is uniquely determined by the condition of being the
smallest complex parabolic subalgebras $\frak{q}'\subset\Hat{\frak{g}}$ 
with\,:
\form{\frak{q}\cap\bar{\frak{q}}\subset\frak{q}'\subset\frak{q}\,.}}
We note that the construction of $\frak{v}$ is independent from the choice 
of the Cartan pair~$(\vartheta,\frak{h})$.\par
\dimo \edef\formulaKC{\number\q.\number\t}
\advance\t by-1
\edef\formulaKB{\number\q.\number\t}\advance\t by1
\edef\proposKA{\number\q.\number\x}
All Cartan subalgebras $\frak{h}$ adapted to
$(\frak{g},\frak{q})$ are also contained in $\bar{\frak{q}}$.
We apply Proposition {\propoIB}. The parabolic set
corresponding to $\frak{v}$ is\,:
\form{\Cal{V}=\Cal{Q}^n\cup\left(\Cal{Q}^r\cap\bar{\Cal{Q}}\right)
=\left(\Cal{Q}^r\cap\bar{\Cal{Q}}^{\,r}\right)
\cup\Cal{Q}^n\cup\left(\Cal{Q}^r\cap\bar{\Cal{Q}}^{\,n}\right)
\,.}
Then it is easy to verify $(\formulaKB)$ by using
Proposition {\propoIB}.\par
All complex parabolic $\frak{q}'\subset\Hat{\frak{g}}$
that satisfy $(\formulaKC)$, also satisfy
${\frak{q}'}\,^n\supset\frak{q}^n$,
because $\frak{q}'\subset\frak{q}$,
and $\frak{q}^r\cap\bar{\frak{q}}\subset\frak{q}\cap\bar{\frak{q}}$\,;
hence $\frak{v}\subset\frak{q}'$ if $\frak{q}'$ satisfies
$(\formulaKC)$. \enddimo
The parabolic subalgebra $(\frak{g},\frak{v})$ 
defined in Proposition {\proposKA} is called {\it the weakest
$CR$ model} of $(\frak{g},\frak{q})$.\par
By Theorem {\TEoremaIG}, since $\frak{q}\cap\bar{\frak{q}}=\frak{v}\cap
\bar{\frak{v}}$, we have\,:
\thm{Let $(\frak{g},\frak{q})$ be an effective parabolic $CR$ algebra
and $(\frak{q},\frak{v})$ its weakest $CR$ model. Then the
holomorphic projection 
$\Hat{M}(\Hat{\frak{g}},\frak{v})@>>>\Hat{M}(\Hat{\frak{g}},\frak{q})$
restricts to a smooth
diffeomorphism $M(\frak{g},\frak{v})@>>>M(\frak{g},\frak{q})$.\qed}
\lem{Let $(\frak{g},\frak{q})$ be an effective parabolic $CR$ algebra
and
$\frak{h}$ a Cartan subalgebra of $\frak{g}$ adapted to
$(\frak{g},\frak{q})$. Let $\Cal{Q}=\Cal{Q}_A=
\{\alpha\in\Cal{R}\,|\,\alpha(A)\geq 0\}$, with $A\in\frak{h}_{\Bbb{R}}$.
Then the parabolic set 
$\Cal{V}$ of the parabolic complex $\frak{v}\subset\Hat{\frak{g}}$ of
the weakest $CR$ model $(\frak{g},\frak{v})$ of $(\frak{g},\frak{q})$
is given by\,:
\form{\Cal{V}=\{\alpha\in\Cal{R}\,|\, \alpha(A+\epsilon\bar{A})\geq 0\}}
where $\epsilon$ is any positive real number with
$\epsilon|\bar\alpha(A)|<\alpha({A})$ for all $\alpha\in\Cal{Q}^n$.}
\dimo Since $\alpha(A+\epsilon\bar{A})<0$ when
$\alpha\notin\Cal{Q}$, we have $\Cal{Q}^n\subset\Cal{V}^n$ and hence
$\Cal{V}\subset\Cal{Q}$.
Moreover $\alpha\in\Cal{Q}^r$ belongs
to $\Cal{V}$ if and only if $\bar{\alpha}(A)\geq 0$, i.e.
if and only if $\alpha\in\bar{\Cal{Q}}$.\enddimo
\edef\LemmaKM{\number\q.\number\x}
\prop{Let $(\vartheta,\frak{h})$ be a Cartan pair adapted to the
effective parabolic $CR$ algebra $(\frak{g},\frak{q})$, and
let $(\frak{g},\frak{w})$ be the canonical lift of
$(\frak{g},\frak{q})$. Then
$(\frak{g},\frak{w})$ coincides with its
weakest $CR$ model.}
\dimo We use the notation of Lemma {\LemmaKM}. By Proposition
{\PropHD}, 
$\frak{w}=\frak{q}_B$ with $B=A-\epsilon\bar{A}\in\frak{h}_{\Bbb{R}}$
for $0<\epsilon<\epsilon_0$. Then, by Lemma {\LemmaKM}, 
 the weakest $CR$ model of
$(\frak{g},\frak{w})$ is $(\frak{g},\frak{v}')$ where
$\frak{v}'=\frak{q}_C$ with $C=A-\epsilon'
\bar{A}\in\frak{h}_{\Bbb{R}}$, with
$\epsilon'=\frac{\epsilon-\delta}{1-\epsilon\delta}$ for 
$0<\delta<\delta_0$ sufficiently small, and hence
$\frak{v}'=\frak{w}$.
\enddimo
We have\,: 
\lem{Let $(\frak{g},\frak{v})$ be the 
weakest $CR$ model of $(\frak{g},\frak{q})$ and
$\frak{h}$ a $CR$ algebra of $\frak{g}$ adapted to 
$(\frak{g},\frak{q})$ (and hence also to $(\frak{g},\frak{v})$).
If
$\Cal{Q}$, $\Cal{V}$ are the parabolic sets corresponding to
$\frak{q}$ and $\frak{v}$ in $\Cal{R}=\Cal{R}(\Hat{\frak{g}},\Hat{\frak{h}})$,
then\,: 
\form{\Cal{Q}^n\cap\bar{\Cal{Q}}^{\,-n}=
 \Cal{V}^n\cap\bar{\Cal{V}}^{\,-n}\,.}}
\dimo We have\,: \edef\formulaKD{\number\q.\number\t}
$\Cal{V}^n=\Cal{Q}^n\cup\left(\Cal{Q}^r\cap\bar{\Cal{Q}}^{\,n}\right)$ and
hence $\bar{\Cal{V}}^{\,-n}=\bar{\Cal{Q}}^{\,-n}\cup
\left(\bar{\Cal{Q}}^{\,r}\cap{\Cal{Q}}^{\,-n}\right)$.
Since\,:\par\centerline{
$\left(\Cal{Q}^r\cap\bar{\Cal{Q}}^{\,n}\right)\cap
\left(\bar{\Cal{Q}}^{\,r}\cap{\Cal{Q}}^{\,-n}\right)=\emptyset$,\quad
$\left(\Cal{Q}^r\cap\bar{\Cal{Q}}^{\,n}\right)\cap
\bar{\Cal{Q}}^{\,-n}=\emptyset$, \quad
$\Cal{Q}^n\cap\left(\bar{\Cal{Q}}^{\,r}\cap{\Cal{Q}}^{\,-n}\right)=\emptyset$,}
\par\noindent we obtain $(\formulaKD)$.
\enddimo
\prop{The effective parabolic $CR$ algebra $(\frak{g},\frak{q})$
is minimal if and only if its $CR$ weakest $CR$ model
$(\frak{g},\frak{v})$ is minimal.}
\dimo 
The statement follows by $(\formulaKD)$ and 
the characterization of parabolic minimal
given in Lemma {\lemmaJE}.\enddimo
\lem{Let $(\frak{g},\frak{q})$ be an effective parabolic $CR$
algebra and $(\frak{g},\frak{v})$ its weakest $CR$ model.
Then a Weyl chamber $C\in\Cal{C}(\Cal{R},\Cal{Q})$ is S-fit for 
$(\frak{g},\frak{q})$ if and only if
$C\in\Cal{C}(\Cal{R},\Cal{V})$ and is S-fit for $(\frak{g},\frak{v})$.}
\dimo
We have $\frak{C}(\Cal{R},\Cal{V})\subset\frak{C}(\Cal{R},\Cal{Q})$
because $\Cal{V}\subset\Cal{Q}$. Since $\Cal{V}^n=\Cal{Q}^n\cup
\left(\Cal{Q}^r\cap\bar{\Cal{Q}}^{\,n}\right)$, for
$C\in\frak{C}(\Cal{R},\Cal{V})$ and $\alpha\in\Cal{V}^n\setminus\Cal{Q}^n$,
we get $\bar\alpha\in\Cal{Q}^n\subset\Cal{V}^n\subset\Cal{R}^+(C)$.
If moreover $C\in\frak{C}(\Cal{R},\Cal{V})$ is S-fit
for $(\frak{g},\frak{v})$, then $\bar\alpha\succ_C0$ also for
$\alpha\in
\left(\Cal{V}^r\cap\Cal{R}^+(C)\right)\setminus
\rootim$; since $\Cal{Q}^r\cap\Cal{R}^+(C)=
\left(\Cal{V}^r\cap\Cal{R}^+(C)\right)\cup\left(\Cal{V}^n\setminus\Cal{Q}^n
\right)$, it follows that $C$ is 
S-fit also for $(\frak{g},\frak{q})$. \par
To complete the proof, it suffices to show that
an S-fit Weyl chamber for $(\frak{g},\frak{q})$ is
admissible for $(\frak{g},\frak{v})$.
Let $C$ be S-fit for $(\frak{g},\frak{q})$. 
The elements $\alpha$ of $\Cal{V}^n\setminus\Cal{Q}^n$
belong to $\Cal{Q}^r\cap\bar{\Cal{Q}}^{\,n}$\,: this implies that
$\bar\alpha\succ_C0$ and therefore that $\alpha\succ_C0$, by the
assumption that $C$ is S-fit for $(\frak{g},\frak{q})$. 
Hence $\Cal{V}^n\subset\Cal{R}^+(C)$ and therefore
$C\in\frak{C}(\Cal{R},\Cal{V})$ when $C\in\frak{C}(\Cal{R},\Cal{Q})$
is S-fit for $(\frak{g},\frak{q})$.
\enddimo
\edef\LemmaGG{\number\q.\number\x}
\prop{Let $(\frak{g},\frak{q})$ be a parabolic $CR$ algebra,
with $\frak{q}=\frak{q}_{\Phi_C}$ for an S-fit Weyl chamber
$C\in\frak{C}(\Cal{R},\Cal{Q})$.
Then its weakest $CR$ model is
$(\frak{g},\frak{q}_{\Phi_C^{\sharp}})$, for
\form{\Phi_C^{\sharp}=\Phi_C\cup\{\beta\in\Cal{B}(C)\cap\rootcp\,|\, 
\supp{C}(\bar\beta)\cap\Phi_C\neq\emptyset\}\,.}}
\dimo Let $C\in\frak{C}(\Cal{R},\Cal{Q})$ be S-fit for
$(\frak{g},\frak{q})$. If $Cal{Q}=\Cal{Q}_A$ for $A\in\frak{h}_{\Bbb{R}}$,
then 
$\Phi_C=\{\beta\in\Cal{B}(C)\,|\,\beta(A)>0\}$ and, by
Lemma {\LemmaKM} and Lemma {\LemmaGG},
$C\in\frak{C}(\Cal{Q},\Cal{V})$ and thus
$\frak{v}=\frak{q}_{\Phi_C^{\sharp}}$ for
$\Phi_C^{\sharp}=
\{\beta\in\Cal{B}(C)\,|\,\beta(A)+\epsilon\bar\beta(A)>0\}$
with $\epsilon>0$ and sufficiently small,
yielding the
characterization in the statement of the Proposition.\enddimo

From Corollary {\corIF} we obtain\,:
\cor{Let $(\frak{g},\frak{v})$ be the weakest $CR$ model of
a parabolic $CR$ algebra $(\frak{g},\frak{q})$. Then\,:
\roster\item"$(i)$" $(\frak{g},\frak{v})$ is either
totally real, or weakly degenerate.
\item"$(ii)$" If $(\frak{g},\frak{q})$ is weakly non-degenerate,
then $\frak{v}=\frak{q}$ if and only if
$(\frak{g},\frak{q})$ is totally real.\qed
\endroster}\edef\corLG{\number\q.\number\x}
By using Corollary {\corLG}, starting from a parabolic $CR$ algebra
$(\frak{g},\frak{q})$, 
we can construct a chain of parabolic $CR$ algebras 
$(\frak{g},\frak{q}_h)$, $(\frak{g},\frak{v}_h)$  and
$\frak{g}$-equivariant 
$CR$ homomorphisms\,:
\form{\CD
(\frak{g},\frak{q})\\
@VVV \\
(\frak{g},\frak{q}_1) @<<< (\frak{g},\frak{v}_1)\\
@. @VVV \\
@. (\frak{g},\frak{q}_2)@<<<(\frak{g},\frak{v}_2) \\
@. @. @VVV \\
@. @. (\frak{g},\frak{q}_3)@<<<(\frak{g},\frak{v}_3) \\
@. @. @.       \vdots 
\endCD}\edef\FRMLFX{\number\q.\number\t}
where each vertical arrow is a weakly nondegenerate reduction and
each horizontal arrow is a lifting to the weakest $CR$ model.\par
If we denote by $\frak{P}(\Hat{\frak{g}})$ the set of all parabolic
complex Lie subalgebras of $\Hat{\frak{g}}$, and set
$\frak{v}_{0}=\frak{g}_{-1}=\frak{g}$, we have, for all integers $h>0$\,:
\form{\cases
\frak{q}_h=\text{the largest $\frak{a}\in\frak{P}(\Hat{\frak{g}})$ such
that}\quad \frak{v}_{h-1}\subset\frak{a}\subset\frak{v}_{h-1}+
\bar{\frak{v}}_{h-1}\,,\\
\frak{v}_h=\text{the smallest $\frak{a}\in\frak{P}(\Hat{\frak{g}})$ such
that}\quad \frak{q}_{h}\cap\bar\frak{q}_h\subset\frak{a}\subset\frak{q}_{h}
\, .
\endcases}
This characterization shows that the construction in $(\FRMLFX)$
is uniquely determined and independent of the choices of the adapted
Cartan pairs in $(\frak{g},\frak{q}_h)$ and $(\frak{g},\frak{v}_h)$.
\par
We know that  the $\bold{G}$-equivariant maps 
$M(\frak{g},\frak{v}_h)@>>>M(\frak{g},\frak{q}_h)$ (for $h>0$)
are smooth diffeomorphisms, while the $\bold{G}$-equivariant maps 
$M(\frak{g},\frak{v}_h)@>>>M(\frak{g},\frak{q}_{h+1})$ (for $h\geq 0$)
are $CR$ fibrations. In particular, there is a smallest integer $m\geq -1$
such that 
$\roman{dim}_{\Bbb{R}}M(\frak{g},\frak{v}_h)>
\roman{dim}_{\Bbb{R}}M(\frak{g},\frak{v}_{h+1})$ for $h<m+1$, and
$M(\frak{g},\frak{v}_h)=M(\frak{g},\frak{v}_{m+1})
=M(\frak{g},\frak{q}_m)$ for all $h>m+1$.
Hence, by the characterization of the fibers in
Theorem {\TeoremaID}, we have\,:
\prop{With the notation above\,: let $(\frak{g},\frak{q}_h)$,
$(\frak{g},\frak{v}_h)$ be the sequence of weakly nondegenerate 
parabolic $CR$ algebra and of their weakest $CR$ models defined above.
Then there exists a smallest integer $m\geq 0$ such that
$\frak{q}_h=\frak{v}_h=\frak{q}_m$ for all $h>m$.
Moreover, $(\frak{g},\frak{q}_m)$ is totally real
and $m$ is the smallest nonnegative integer for which
$(\frak{g},\frak{q}_m)$ is totally real.\par
By composition, we obtain
a $\bold{G}$-equivariant  fibration 
$M(\frak{g},\frak{q}) @>>> M(\frak{g},\frak{q}_m)$ where the basis
$M(\frak{g},\frak{q}_m)$ is a totally real parabolic $CR$ manifold
and each connected component of the fiber is
a Cartesian product of Euclidean complex nilmanifolds and of
simply connected
totally complex parabolic $CR$ manifolds.
\qed}
\edef\ProposizioneLJ{\number\q.\number\x}
\esemp{
Let $\epsilon_1,\hdots,\epsilon_6$ be the canonical basis of
$\Bbb{R}^6\subset\Bbb{C}^6$. Let $\bold{G}=\bold{SL}(6,\Bbb{R})$ consist
of the matrices of $\Hat{\bold{G}}=\bold{SL}(6,\Bbb{C})$ which have real
entries in the canonical basis.
We consider 
in $\Bbb{C}^6$ the basis\,:
$$e_1=\epsilon_1+i\epsilon_4,\;
e_2=\epsilon_2+i\epsilon_5,\;
e_3=\epsilon_3+i\epsilon_6,\;
e_4=\epsilon_1-i\epsilon_4,\;
e_5=\epsilon_2-i\epsilon_5,\;
e_6=\epsilon_3-i\epsilon_6$$
and we want to investigate the $\bold{G}$-orbit $M$ of the flag
$$\langle e_1\rangle\subset
\langle e_1, e_2\rangle\subset
\langle e_1, e_2,e_3\rangle\subset
\langle e_1, e_2,e_3,e_4\rangle\subset
\langle e_1, e_2,e_3,e_4,e_5\rangle\,.
$$
Let $(\frak{g},\frak{q})$, with $\frak{g}=\frak{sl}(6,\Bbb{R})$,
be the associated parabolic $CR$ algebra.
We consider the Cartan subalgebra $\frak{h}$ of $\frak{sl}(6,\Bbb{R})$
that is represented, in the basis $e_1,\hdots,e_6$, by the diagonal 
matrices $\roman{diag}(\lambda_1,\lambda_2,\lambda_3,
\bar\lambda_1,\bar\lambda_2,\bar\lambda_3)$ with
$\roman{Re}(\lambda_1+\lambda_2+\lambda_3)=0$. Then
$\frak{h}_{\Bbb{R}}$ consists of the $6\times 6$ traceless real diagonal
matrices. \par
The cross-marked diagram associated to $(\frak{g},\frak{q})$ in the
adapted Weyl chamber $C$ with simple roots
$\Cal{B}(C)=\{\alpha_i=e_1-e_{i+1}\,,\; 1\leq i\leq 5\}$, is the following\,:
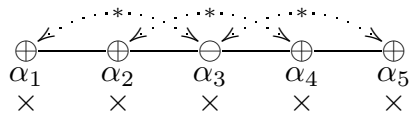
$$
\xymatrix @M=0pt @R=2pt @!C=5pt{
\oplus\ar@{-}[r]\ar@/^1pc/@{<.>}[rr]|{\ast}
&\oplus\ar@{-}[r]\ar@/^1pc/@{<.>}[rr]|{\ast}&\ominus\ar@{-}[r]
\ar@/^1pc/@{<.>}[rr]|{\ast}&\oplus
\ar@{-}[r]&\oplus\\
\alpha_1&\alpha_2&\alpha_{3}&\alpha_4&\alpha_5\\
\times&\times&\times&\times&\times}
$$
Since $\frak{q}$ is a complex Borel subalgebra of $\Hat{\frak{g}}$,
the chamber $C$ is S-fit and then $(\frak{g},\frak{q})$ is fundamental
by Theorem {\TeorEF},
because $\bar\alpha_3=-(\alpha_1+\alpha_2+\alpha_3+\alpha_4+\alpha_5)$.
The chamber $C$ is also V-fit, and we obtain the weakly non-degenerate
basis $(\frak{g},\frak{q}_1)$ 
by dropping the cross under the simple root $\alpha_3$
with 
$\bar\alpha_3\prec_C0$. The diagram associated to
$(\frak{g},\frak{q}_1)$ is\,:
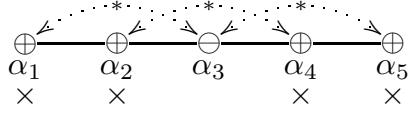
$$
\xymatrix @M=0pt @R=2pt @!C=5pt{
\oplus\ar@{-}[r]\ar@/^1pc/@{<.>}[rr]|{\ast}
&\oplus\ar@{-}[r]\ar@/^1pc/@{<.>}[rr]|{\ast}
&\ominus\ar@{-}[r]\ar@/^1pc/@{<.>}[rr]|{\ast}
&\oplus
\ar@{-}[r]&\oplus\\
\alpha_1&\alpha_2&\alpha_{3}&\alpha_4&\alpha_5\\
\times&\times&&\times&\times}
$$
The parabolic $\frak{q}_1$ is defined by the element
$A_1=\roman{diag}(2,1,0,0,-1,-2)\in\frak{h}_{\Bbb{R}}$. 
To compute its weakest $CR$ model $(\frak{g},\frak{v}_1)$,
we observe that $\bar{A}_1=\roman{diag}(0,-1,-2,2,1,0)$, so that
$A_1+\epsilon\bar{A}_1=
\roman{diag}(2,1-\epsilon,-2\epsilon,2\epsilon,-1+\epsilon,-2)$.
To take a Weyl chamber adapted to 
$(\frak{g},\frak{v}_1)$, 
it is convenient to consider  the basis 
obtained from $e_1,\hdots,e_6$ by reordering its elements according
to the decreasing ordering of
the diagonal entries of $A_1+\epsilon\bar{A}_1$. We obtain the new
basis\,:
$e_1,e_2,e_4,e_3,e_5,e_6\,.$ 
With 
$\alpha'_1=e_1-e_2,\;
\alpha_2'=e_2-e_4,\;
\alpha_3'=e_4-e_3,\;
\alpha_4'=e_3-e_5,\;
\alpha_5'=e_5-e_6$
being the simple roots
a Weyl chamber $C'\in\frak{C}(\Cal{R},\Cal{V}_1)$,
we obtain the diagram\,:
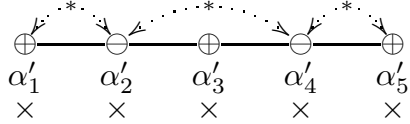
$$
\xymatrix @M=0pt @R=2pt @!C=5pt{
\oplus\ar@{-}[r]\ar@/^1pc/@{<.>}[r]|{\ast}
&\ominus\ar@{-}[r]\ar@/^1pc/@{<.>}[rr]|{\ast}
&\oplus\ar@{-}[r]
&\ominus\ar@/^1pc/@{<.>}[r]|{\ast}
\ar@{-}[r]&\oplus\\
\alpha_1'&\alpha_2'&\alpha_{3}'&\alpha_4'&\alpha_5'\\
\times&\times&\times&\times&\times}
$$
The weakly nondegenerate reduction $(\frak{g},\frak{q}_2)$ of
$(\frak{g},\frak{v}_1)$ has the diagram\,:
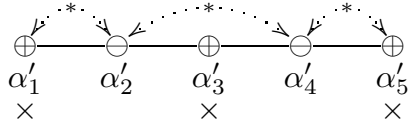
$$
\xymatrix @M=0pt @R=2pt @!C=5pt{
\oplus\ar@{-}[r]\ar@/^1pc/@{<.>}[r]|{\ast}
&\ominus\ar@{-}[r]\ar@/^1pc/@{<.>}[rr]|{\ast}
&\oplus\ar@{-}[r]
&\ominus\ar@/^1pc/@{<.>}[r]|{\ast}
\ar@{-}[r]&\oplus\\
\alpha_1'&\alpha_2'&\alpha_{3}'&\alpha_4'&\alpha_5'\\
\times&&\times&&\times}
$$
We have $\frak{q}_2=\frak{q}_{A_2}$, with\par\centerline{
$A_2=\roman{diag}(2,1,-1,1,-1,-2)$,\quad
$\bar{A}_2=\roman{diag}(1, -1, -2, 2, 1, -1)$,}\par\noindent  so that
$A_2+\epsilon\bar{A}_2=
\roman{diag}(2+\epsilon,1-\epsilon,-1-2\epsilon,1+2\epsilon,-1+\epsilon,
-2-\epsilon)$.
To describe $(\frak{g},\frak{v}_2)$, for $\frak{v}_2=\frak{q}_{A_2+
\epsilon\bar{A}_2}$, it is convenient to consider the Weyl chamber
$C''\in\frak{C}(\Cal{R},\Cal{V}_2)$ that corresponds to the simple
roots related 
to the ordered basis
$e_1,e_4,e_2,e_5,e_3,e_6$ of $\Bbb{C}^6$, for which 
the entries of $A_2+\epsilon\bar{A}_2$ are decreasing:
with $\alpha_1''=e_1-e_4,\; 
\alpha_2''=e_4-e_2,\; 
\alpha_3''=e_2-e_5,\; 
\alpha_4''=e_5-e_3,\; 
\alpha_5''=e_3-e_6\,$, 
we obtain for $(\frak{g},\frak{v}_2)$ the diagram\,:
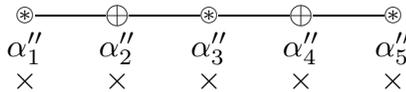
$$
\xymatrix @M=0pt @R=2pt @!C=5pt{
\ggrigia\ar@{-}[r]
&\oplus\ar@{-}[r]&\ggrigia\ar@{-}[r]
&\oplus
\ar@{-}[r]&\ggrigia\\
\alpha_1''&\alpha_2''&\alpha_{3}''&\alpha_4''&\alpha_5''\\
\times&\times&\times&\times&\times}
$$
Since $\frak{v}_2$ is Borel, $C''$ is V-fit for $(\frak{g},\frak{v}_2)$
and
the diagram of the
weakly non-degenerate reduction $(\frak{g},\frak{q}_3)$ of
$(\frak{g},\frak{v}_2)$ is obtained by dropping the crosses
under the $\alpha''_i$'s with $\bar\alpha''_i\prec_{C''}0$\,:
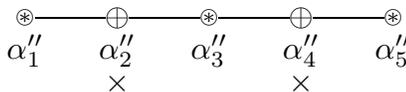
$$
\xymatrix @M=0pt @R=2pt @!C=5pt{
\ggrigia\ar@{-}[r]
&\oplus\ar@{-}[r]&\ggrigia\ar@{-}[r]
&\oplus
\ar@{-}[r]&\ggrigia\\
\alpha_1''&\alpha_2''&\alpha_{3}''&\alpha_4''&\alpha_5''\\
&\times&&\times&}
$$
The parabolic $CR$ algebra
$(\frak{g},\frak{q}_3)$ is totally real, as
$\frak{q}_3=\frak{q}_{A_3}$ with\par\centerline{
$A_3=\roman{diag}(1,0,-1,1,0,-1)=\bar{A}_3$.}\par\noindent Hence
$(\frak{g},\frak{q}_2)=(\frak{g},\frak{v}_h)=(\frak{g},\frak{q}_h)$
for all $h\geq 3$.
The map $M(\frak{g},\frak{q})@>>>M(\frak{g},\frak{q}_3)$ 
is given by $(\ell_1,\ell_2,\ell_3,\ell_4,\ell_5) @>>> (\ell_1+\bar\ell_1,
\ell_2+\bar\ell_2)$. }
\edef\ESMPKZ{\number\q.\number\y}


\se{The isotropy subgroups}
Let $\frak{g}$ be a real semisimple Lie algebra, 
$\Hat{\frak{g}}$ its complexification, $\vartheta$ a Cartan involution 
of $\frak{g}$. Let $\frak{h}$ be a $\vartheta$-invariant Cartan subalgebra
of $\frak{g}$ and
$\Cal{R}=\Cal{R}(\Hat{\frak{g}},\Hat{\frak{h}})$ the root system
of $\Hat{\frak{g}}$ with respect to the complexification $\Hat{\frak{h}}$
of $\frak{h}$.
Denote by $\Lambda(\Cal{R})$ the additive subgroup of $\frak{h}^*_{\Bbb{R}}$
generated by $\Cal{R}$ and by $\Pi(\Cal{R})$ the {\it lattice of weights},
consisting of all $\eta\in\frak{h}_{\Bbb{R}}^*$ for which
$(\eta|\alpha^{\vee})=2(\eta|\alpha)/\|\alpha\|^2\in\Bbb{Z}$.
We have $\Lambda(\Cal{R})\subset\Pi(\Cal{R})$.\par
Given a lattice (i.e. a free Abelian group) $\Cal{L}$, a 
{\it character} of $\Cal{L}$ is
a homomorphism $\chi:\Cal{L}@>>>\Bbb{C}^*$
of $\Cal{L}$ into the multiplicative group
$\Bbb{C}^*=\Bbb{C}\setminus\{0\}$ of non-zero complex numbers.
If $b_1,\hdots,b_{\ell}\in\Cal{L}$ is a basis of
$\Cal{L}$ over $\Bbb{Z}$, a character $\chi\in\Hom(\Cal{L},\Bbb{C}^*)$
is completely determined by its values 
$\lambda_i=\chi(b_i)$ (for $1\leq i\leq\ell$) on the basis, so that
$\Hom(\Cal{L},\Bbb{C}^*)\simeq[\Bbb{C}^*]^\ell$\,.\par
We keep the notation of the previous sections, in particular
$\Hat{\bold{G}}$ is a connected and simply connected
complex Lie group with Lie algebra $\Hat{\frak{g}}$ and $\bold{G}$
its analytic subgroup with Lie algebra $\frak{g}$. 
It is a covering group of any linear group with Lie algebra $\frak{g}$.
Let $\Hat{\bold{H}}$
be the Cartan subgroup of $\Hat{\bold{G}}$ corresponding to $\Hat{\frak{h}}$\,:
\form{\Hat{\bold{H}}=\bold{Z}_{\Hat{\bold{G}}}(\Hat{\frak{h}})\,=\,
\{z\in\Hat{\bold{G}}\,|\,\roman{Ad}_{\Hat{\frak{g}}}(z)(H)=H\,,\;
\forall H\in\Hat{\frak{h}}\,\}\,.}\edef\FRMXX{\number\q.\number\t}
All finite dimensional $\Bbb{C}$-linear representations of 
the complex semisimple Lie algebra $\Hat{\frak{g}}$ are differentials
of representations of the complex Lie group $\Hat{\bold{G}}$.
Each element $h$ of $\Hat{\bold{H}}$ defines a character 
$\chi_h\in\Hom(\Pi(\Cal{R}),\Bbb{C}^*)$, and vice versa, a character
$\chi\in\Hom(\Pi(\Cal{R}),\Bbb{C}^*)$ defines the element 
$h_{\chi}\in\Hat{\bold{H}}$. To explain this correspondence,
take first a faithful representation 
$\boldsymbol{\rho}:\Hat{\bold{G}}\hookrightarrow
\bold{SL}_{\Bbb{C}}(V)$, corresponding to 
$\rho:\Hat{\frak{g}}\hookrightarrow\frak{sl}_{\Bbb{C}}(V)$,
for a finite dimensional complex linear space
$V$\,, and define $\boldsymbol{\rho}(h_{\chi})(v)=\chi(\omega)\,v$ for 
$v\in V^{\omega}=\{v\in V\,|\, \rho(H)(v)=\omega(H)\,v\,,\;\forall H\in
\Hat{\frak{h}}\}$\,, for $\omega\in\Pi(\Cal{R})$. The Cartan subgroup
$\Hat{\bold{H}}$ is analytic and 
$\exp\colon\Hat{\frak{h}}@>>>\Hat{\bold{H}}$ is onto, so that 
the correspondence with
$\Hom(\Pi(\Cal{R}),\Bbb{C}^*)$ can also be described by
$\chi_{\exp(H)}(\omega)=\exp(\omega(H))$ for $H\in\Hat{\frak{h}}$.
With 
$\ell=\roman{dim}_{\Bbb{C}}(\Hat{\frak{h}})$=rank of $\Hat{\frak{g}}$
we have\,:
\form{\Hat{\bold{H}}=\bold{Z}_{\Hat{\bold{G}}}(\Hat{\frak{h}})=
\{h_\chi\,|\,\chi\in\Hom(\Pi(\Cal{R}),\Bbb{C}^*)\} \simeq
\Hom(\Pi(\Cal{R}),\Bbb{C}^*)\simeq\left[\Bbb{C}^*\right]^\ell\,.}
\par\smallskip
The Cartan subgroup $\bold{H}$ of $\bold{G}$ corresponding to $\frak{h}$
is the centralizer of $\frak{h}$ in $\bold{G}$\,:
\form{\bold{H}=\bold{Z}_{\bold{G}}(\frak{h})=\{h\in\bold{G}\,|\,
\roman{Ad}_{\frak{g}}(h)(H)=H\,,\;\forall H\in\frak{h}\,\}\,.}
For a lattice $\Cal{L}\subset\frak{h}_{\Bbb{R}}^*$, with
$\sigma(\Cal{L})=\Cal{L}$, we set\par\centerline{
$\Hom^{\sigma}(\Cal{L},\Bbb{C}^*)=\{\chi\in\Hom(\Cal{L},\Bbb{C}^*)\,|\,
\chi(\bar\eta)=\overline{\chi(\eta)}\,,\;\forall\eta\in \Cal{L}\}$.}\par
We obtain\,:
\lem{Let $\frak{g}$ be a real semisimple Lie algebra, $\frak{h}$
a Cartan subalgebra of $\frak{g}$. Then\,: 
\form{\bold{H}=\bold{Z}_{\bold{G}}(\frak{h})=\{h_{\chi}\,|\,
\chi\in\Hom^{\sigma}(\Pi(\Cal{R}),\Bbb{C}^*)\,\}\,.}}
\edef\LEMAIA{\number\q.\number\x}
\dimo
The action of $\sigma$ on $\Hom(\Pi(\Cal{R}),\Bbb{C}^*)$, given by 
$\sigma(\chi)(\eta)=\overline{\chi(\bar{\eta})}$ coincides, under 
the correspondence $(\FRMXX)$, with the action of $\sigma$ on $\Hat{\bold{H}}$.
Then the statement follows from the fact that $\bold{H}=
\Hat{\bold{H}}\cap\bold{G}$ 
and $\bold{H}=\Hat{\bold{H}}^{\sigma}$ because of $(\FRMLDB)$.
\enddimo

In view of the preceding Lemma, we now give a more explicit description of
$\Hom^{\sigma}(\Pi(\Cal{R}),\Bbb{C}^*)$.
Fix a Weyl chamber $C\in\frak{C}(\Cal{R})$ 
that is S-adapted to the conjugation
$\sigma$\,.
Set $\Cal{B}(C)=\{\beta_1,\hdots,\beta_a\}\cup\{\tau_1,\hdots,\tau_b\}$,
with $\{\beta_1,\hdots,\beta_a\}=\Cal{B}(C)\setminus\rootim$ and
$\{\tau_1,\hdots,\tau_b\}=\Cal{B}(C)\cap\rootim$. By Theorem {\teoremaCI},
the conjugation $\sigma$ is described in $\Cal{B}(C)$
by an involutive 
permutation $\jmath:\beta_i@>>>\jmath(\beta_i)=\beta_{i'}$ of
$\{\beta_1,\hdots,\beta_a\}$ and a matrix of nonnegative integers
$(k_{i,q})$\,,
$1\leq i\leq a\,,\;
1\leq q\leq b$\,, 
such that\,:
\form{\cases
\bar\beta_i=\beta_{i'}+\sum_{q=1}^{b}{k_{i,q}\tau_q}\,,\\
\jmath(\beta_i)=\beta_{i'}\,,\quad k_{i',p}=k_{i,p} 
\quad\text{for}\; 1\leq i\leq a\,,\; 1\leq p\leq b\,.
\endcases}
In
$\Pi(\Cal{R})$ we consider 
the basis $\Cal{B}^*(C)=
\{\omega_1,\hdots,\omega_a,\theta_1,\hdots,\theta_b\}$, 
adjoint of $\Cal{B}(C)$\,,
defined by\,:
\form{\left\{{\matrix\format\l\quad&\l\\
(\omega_i|\beta_j^{\vee})=\delta_{i,j} 
& \; (\omega_i|\tau_q^{\vee})=0 \\\vspace{2\jot}
(\theta_p|\beta_j^{\vee})=0 &\;
(\theta_p|\tau_q^{\vee})=\delta_{p,q} 
\endmatrix}\right.\qquad
\text{for}\; 1\leq i,j\leq a\,, \;\quad 1\leq p,q\leq b\,.}
The conjugation $\sigma$ is described in $\Cal{B}^*(C)$ by\,:
\form{\cases
\bar\omega_i=\omega_{i'} &\text{for}\; 1\leq i,i'\leq a\,,\;
\jmath(\beta_i)=\beta_{i'}\\
\bar\theta_p=-\theta_p+\sum_{j=1}^a{k'_{j,p}\omega_j} 
&\text{for}\; 1\leq p\leq b\,,
\endcases}
\par \smallskip
\noindent where $k'_{j,p}=k_{j,p} \|\tau_p\|^2/\|\beta_j\|^2$.
\par\smallskip
The characters $\chi\in\Hom^{\sigma}(\Pi(\Cal{R}),\Bbb{C}^*)$
are those satisfying\,:
\form{\cases
\chi(\omega_i)\in\Bbb{C}^*\,,\; \chi(\theta_p)\in\Bbb{C}^*\,,&\text{for}\;
1\leq i\leq a\,,\;1\leq p\leq b\,,\\
\overline{\chi(\omega_i)}\,=\,\chi(\omega_{i'})
&\text{for}\; 1\leq i,i'\leq a\,,\; \jmath(\beta_i)=\beta_{i'}\\
|\chi(\theta_p)|^2=\prod_{j=1}^a{\left[\chi(\omega_j)\right]^{k'_{j,p}}}
&\text{for}\; 1\leq p\leq b\,.
\endcases}
\par\smallskip
Each lattice 
$\Cal{L}$ in $\frak{h}_{\Bbb{R}}^*$, with
$\Lambda(\Cal{R})\subset\Cal{L}\subset\Pi(\Cal{R})$,
is the set of weights of all 
finite dimensional linear representations
of an essentially unique connected
complex semisimple Lie group $\Hat{\bold{G}}_{\Cal{L}}$.
For instance, we have
$\Hat{\bold{G}}_{\Pi(\Cal{R})}=\Hat{\bold{G}}$ 
(the simply connected complex Lie group with Lie algebra $\Hat{\frak{g}}$)
and
$\Hat{\bold{G}}_{\Lambda(\Cal{R})}=\Interc{g}$ (the group of inner
automorphisms of the complex semisimple Lie algebra $\Hat{\frak{g}}$)
(see e.g. [OV, Ch.3, Theorem 2.11]).
When moreover
$\sigma(\Cal{L})=\Cal{L}$, the analytic Lie subgroup of 
$\Hat{\bold{G}}_{\Cal{L}}$ with Lie algebra $\frak{g}$ is a real form
$\bold{G}_{\Cal{L}}$ of $\Hat{\bold{G}}_{\Cal{L}}$.
Vice versa, every connected linear semisimple Lie group with Lie algebra
$\frak{g}$ can be obtained in this way.
The Cartan
subgroup\,:   
\form{\bold{H}_{\Cal{L}}=\bold{Z}_{\bold{G}_{\Cal{L}}}(\frak{h})}
\noindent of $\bold{G}_{\Cal{L}}$, relative to $\frak{h}$,  is a 
real Lie subgroup of
the complex Cartan subgroup\,:  
\form{\Hat{\bold{H}}_{\Cal{L}}=
\bold{Z}_{\Hat{\bold{G}}_{\Cal{L}}}(\Hat{\frak{h}})
=\{h_{\chi}\,|\,\chi\in\Hom(\Cal{L},\Bbb{C}^*)\}} 
\noindent of
$\Hat{\bold{G}}_{\Cal{L}}$, relative to $\Hat{\frak{h}}$.
From Lemma {\LEMAIA} we obtain\,:
\cor{Let $\Cal{L}$ be a lattice in $\frak{h}_{\Bbb{R}}^*$ with
$\Lambda(\Cal{R})\subset\Cal{L}\subset\Pi(\Cal{R})$ and
$\sigma(\Cal{L})=\Cal{L}$. Denote by
$\chi@>>>\chi^\flat$ the restriction homomorphism
$\Hom(\Pi(\Cal{R}),\Bbb{C}^*) @>>> \Hom(\Cal{L},\Bbb{C}^*)$. Then\,:
\form{\bold{H}_{\Cal{L}}=\left.\left\{h_{\chi^\flat}\,\right|
\chi\in\Hom^\sigma(\Pi(\Cal{R}),\Bbb{C}^*)\right\}\,.}}
\dimo The covering map $\bold{G}@>>>\bold{G}_{\Cal{L}}$ transforms
the Cartan subgroup $\bold{H}$ of $\bold{G}$
into the Cartan subgroup $\bold{H}_{\Cal{L}}$
of $\bold{G}_{\Cal{L}}$ . \enddimo \edef\CORIB{\number\q.\number\x}

From Corollary {\CORIB} we obtain the exact sequence\,:
\form{
\bold{1}@>>>\Hom^{\sigma}(\Pi(\Cal{R})/\Cal{L},\Bbb{C}^*)
@>>>\Hom^{\sigma}(\Pi(\Cal{R}),\Bbb{C}^*)
@>>>\bold{H}_{\Cal{L}}@>>>\bold{1}\,.
}
\edef\FFrmC{\number\q.\number\t}
We can utilize $(\FFrmC)$ to compute the group $\pi_0(\bold{H}_{\Cal{L}})$ 
of the connected components 
of the Cartan subgroup $\bold{H}_{\Cal{L}}$. We have indeed\,:
\thm{Keep the previous notation.
Let $C\in\frak{C}(\Cal{R})$ be
S-adapted to the
conjugation $\sigma$. Denote by $\omega_1,\hdots,\omega_a$ the weights
in $\Cal{B}^*(C)$ that vanish on $\Cal{B}(C)~\cap~\rootim$, and by 
$\{\theta_1,\hdots,\theta_b\}$ those 
vanishing on $\Cal{B}(C)\setminus\rootim$. By reordering, we
assume that
$\{\omega_1,\hdots,\omega_c\}$ is the set of weights in
$\Cal{B}^*(C)$ with 
$\bar\omega_i=\omega_i$. 
Define the non negative integers 
$k'_{i,p}$, for $1\leq i\leq a$, $1\leq p\leq b$, by
$\bar\theta_p=-\theta_p+\sum_{i=1}^a{k'_{i,p}\omega_i}$.
Consider the subgroups
of the free Abelian group $\Bbb{Z}_2^c$\,:
\form{\matrix\format\r&\l\\
\bold{A}&=\{(\eta_1,\hdots,\eta_c)\in\Bbb{Z}_2^c
\,|\, \sum_{i=1}^c{k'_{i,p}\eta_i}\equiv 0\mod{2}\,,\;
\forall\, 1\leq p\leq b\}\,,\\\vspace{2\jot}
\bold{A}_{\Cal{L}}&=\{
(\eta_1,\hdots,\eta_c)\in\bold{A}
\,|\, \sum_{i=1}^c{k_i\eta_i}\equiv 0\mod{2}\,,\;\text{if}\;
\sum_{i=1}^c{k_i\omega_i}\in\Cal{L}
\}\,.\endmatrix}
Then\,:
\form{\pi_0(\bold{H}_{\Cal{L}})\cong \bold{A}/
\bold{A}_{\Cal{L}}\,.}}
\dimo The exact sequence $(\FFrmC)$
yields the exact sequence for the groups of 
the connected components\,:
$$\pi_0(\Hom^{\sigma}(\Pi(\Cal{R})/\Cal{L},\Bbb{C}^*))
@>>>
\pi_0(\Hom^{\sigma}(\Pi(\Cal{R}),\Bbb{C}^*))
@>>> \pi_0(\bold{H}_{\Cal{L}})@>>>\bold{1}\,.$$ 
The statement follows because
$\pi_0(\Hom^{\sigma}(\Pi(\Cal{R}),\Bbb{C}^*))\cong \bold{A}$,
and, in this isomorphism,
the image of $\pi_0(\Hom^{\sigma}(\Pi(\Cal{R})/\Cal{L},\Bbb{C}^*))$
in $\pi_0(\Hom^{\sigma}(\Pi(\Cal{R}),\Bbb{C}^*))$ is identified with
$\bold{A}_{\Cal{L}}$.
\enddimo
\esemp{Consider the $\roman{B}_2$ system of roots
$\Cal{R}=\{\pm e_1\pm e_2\}\cup\{\pm e_1,\pm e_2\}\subset\Bbb{R}^2$,
with $\sigma$ defined by $\sigma(e_1)=e_2$ and $\sigma(e_2)=e_1$.
In the Weyl chamber $C$ with simple roots $\{e_1-e_2,e_2\}$, that
is S-adapted to $\sigma$, it can be represented by the diagram\,:
$$
\xymatrix @M=0pt @R=2pt @!C=5pt{
{\circledast}\ar@{=}[rr]&>&\bbianca\\
\tau&&\beta}
$$
We have, for the adjoint basis, $\theta=e_1$ and $\omega=(e_1+e_2)/2$.
We obtain $\bar\theta=-\theta+2\omega$ and $\bar\omega=\omega$.
Then $k'=2$, $m=1$, and $\bold{A}=\Bbb{Z}_2$. The quotient 
$\Pi(\Cal{R})/\Lambda(\Cal{R})$ is a group of order $2$, with the generator
$[\omega]$ and $2[\omega]\equiv 0$. Hence
$\pi_0(\bold{H})\simeq\Bbb{Z}_2$, while  
$\pi_0(\bold{H}_{\Lambda(\Cal{R})})\cong\Bbb{Z}_2/\Bbb{Z}_2\cong\bold{1}$.}
\par
Let $\bold{G}_{\Cal{L}}$ be any connected real linear Lie group with
Lie algebra $\frak{g}$, and $\frak{h}$ a Cartan subalgebra of $\frak{g}$.
The Cartan subgroup
$\bold{H}_{\Cal{L}}$ 
is a normal subgroup of the normalizer
\form{\bold{N}_{{\Cal{L}}}=
\bold{N}_{\bold{G}_{\Cal{L}}}({\frak{h}})=
\{g\in\bold{G}_{\Cal{L}}\,|\, \roman{Ad}(g)(\frak{h})=\frak{h}\}}
of $\frak{h}$ in $\bold{G}_{\Cal{L}}$\,.
The quotient
\form{\bold{W}_{\Cal{L}}=
\bold{N}_{\Cal{L}}/\bold{H}_{\Cal{L}}}
is the {\it analytic Weyl group} corresponding to $\frak{h}$ and
$\Cal{L}$. \par
We know (see e.g.
[Wa, Ch.1,\S{4},p.115]) that actually
the {\it analytic Weyl group} only depends, modulo natural
isomorphisms, upon the real Lie algebra $\frak{g}$ and its
Cartan subalgebra $\frak{h}$;  to stress this fact, we shall 
write
$\bold{W}(\frak{g},\frak{h})$ instead of 
$\bold{W}_{\Cal{L}}$\,.
\par\smallskip
Let us consider now an effective 
parabolic $CR$ algebra $(\frak{g},\frak{q})$. 
\par
We keep the notation introduced in previous sections, in particular we
fix a Cartan pair $(\vartheta,\frak{h})$ adapted to
$(\frak{g},\frak{q})$.\par
We have the decomposition 
\form{\frak{g}_+=\frak{q}\cap\frak{g}=
\frak{n}\oplus\frak{g}_0=\frak{r}\oplus\frak{s}_0=
\frak{n}\oplus\frak{z}_0
\oplus\frak{s}_0}\edef\FormulaIN{\number\q.\number\t}
where\,:\par
$\frak{r}=\frak{n}\oplus\frak{z}_0$ is the radical of $\frak{g}_+$\,,\par
$\frak{n}$ is the ideal of the $\roman{ad}_{\frak{g}}$-nilpotent
elements of the radical of $\frak{g}_+$, \par
$\frak{z}_0$ is a maximal Abelian subalgebra of 
$\roman{ad}$-semisimple elements
of $\frak{r}$\,,\par
$\frak{g}_0=\frak{z_0}\oplus\frak{s}_0$ is
the $\vartheta$-invariant reductive complement of $\frak{n}$ in
$\frak{g}_+$, \par
$\frak{s_0}=[\frak{g}_0,\frak{g}_0]$ is the semisimple ideal 
and $\frak{z}_0$ is the center of $\frak{g}_0$. \par\smallskip
The Cartan subalgebra $\frak{h}$ of $\frak{g}$, contained in $\frak{g}_+$,
decomposes into the direct sum\,:
\form{\frak{h}=\frak{z}_0\oplus\frak{h}_0}\edef\FormulaIO{\number\q.\number\t}
where $\frak{h}_0=\frak{h}\cap\frak{s}_0$ is a Cartan subalgebra of
$\frak{s}_0$. The subalgebra $\frak{z}_0$ is characterized by\,:
\form{\frak{z}_0\,=\,\{H\in\frak{h}\,|\,\alpha(H)=0\,,\;\forall\alpha\in
\Cal{Q}^r\cap\bar{\Cal{Q}}^{\,r}\,\}\,.}\edef\FRMIZ{\number\q.\number\t}
Indeed its complexification $\Hat{\frak{z}}_0$ is the center of
the reductive complex Lie algebra $\frak{q}^r\cap\bar{\frak{q}}^{\,r}$.
\par
The isotropy subgroup 
\form{\bold{G}_+=\bold{N}_{\bold{G}}(\frak{q})=\{g\in\bold{G}\,|\,
\roman{Ad}(g)(\frak{q})=\frak{q}\}}\edef\FormulaIP{\number\q.\number\t}
is a closed real-algebraic
subgroup of $\bold{G}$, with Lie algebra
$\frak{g}_+$. Thus we have the
{Chevalley decomposition}\,:
\form{{\bold{G}}_+=
\bold{G}_0\rtimes\bold{N}}\edef\FormulaIQ{\number\q.\number\t}
where $\bold{G}_0$ is a closed, real-algebraic, reductive subgroup
of $\bold{G}$, 
with Lie algebra $\frak{g}_0$, and $\bold{N}$ is
a unipotent closed connected and simply connected subgroup of $\bold{G}$,
diffeomorphic to its Lie algebra $\frak{n}$.\par\smallskip
We also define $\bold{S}_0$ to be the analytic 
semisimple Lie subgroup of $\bold{G}$
with Lie algebra $\frak{s}_0$.
\prop{The group $\bold{G}_0$ is the normalizer in $\bold{G}$ of
$\frak{g}_0$ and the centralizer in $\bold{G}$ of $\frak{z}_0$\,:
\form{\matrix\format\r&\l\\
\bold{G}_0&
=\bold{N}_{\bold{G}}(\frak{g}_0)=
\{g\in\bold{G}\,|\,\roman{Ad}(g)(\frak{g}_0)=\frak{g}_0\}
\\\vspace{2\jot}
&=\bold{Z}_{\bold{G}}(\frak{z}_0)\,=\,
\{g\in\bold{G}\,|\,\roman{Ad}(g)(H)=H\,,\;\forall H\in\frak{z}_0\,\}\,.
\endmatrix}} 
\dimo \edef\TRMIH{\number\q.\number\x}

While proving this Proposition we can assume that 
$\frak{h}$ is maximally noncompact
among the Cartan subalgebras of $\frak{g}$
that are contained in $\frak{g}_+$. Indeed, all admissible Cartan subalgebras
of $(\frak{g},\frak{q})$ are conjugate to the direct sum of $\frak{z}_0$
and a Cartan subalgebra $\frak{h}_0$ of $\frak{s}_0$. Thus our assumption
means that we have chosen a maximally noncompact Cartan subalgebra
$\frak{h}_0$ of $\frak{s}_0$.
\par
We first prove that
\form{\bold{G}_0=\{g\in\bold{G}_+\,|\,
\roman{Ad}(g)(\frak{g}_0)=\frak{g}_0\}\,.}

If $g\in\bold{G}_+$,
then $\roman{Ad}(g)(\frak{g}_0)$ is a maximal reductive
subalgebra of $\frak{g}_+$.
By [OV, Ch.I,\S{6.5}] there is an element $n\in\bold{N}$ such that
$\frak{g}_0=\roman{Ad}(g\,n^{-1})(\frak{g}_0)$. This yields a decomposition
$g=g_0\,n$, with $g_0\in\bold{N}_{\bold{G}_+}\!(\frak{g}_0)$ and
$n\in\bold{N}$. Since $\frak{z}_0\subset\frak{g}_0$, and
$[\frak{z}_0,\frak{n}]=\frak{n}$,
we have $\bold{N}\cap\bold{N}_{\bold{G}_+}\!(\frak{g}_0)\,=\,\{1\}$.
This yields the decomposition $\bold{G}_+=\bold{N}_{\bold{G}_+}\!(\frak{g}_0)
\rtimes\bold{N}$. Since we also have $(\FormulaIQ)$ and
$\frak{g}_0$ is the Lie algebra of both $\bold{G}_0$ and
$\bold{N}_{\bold{G}_+}\!(\frak{g}_0)$, we obtain that 
$\bold{G}_0=\bold{N}_{\bold{G}_+}\!(\frak{g}_0)$\,.
\par
Next we show that $\bold{Z}_{\bold{G}}(\frak{z}_0)=
\bold{N}_{\bold{G}}(\frak{g}_0)$.
If $g_0\in\bold{N}_{\bold{G}}(\frak{g}_0)$, 
then $\roman{Ad}(g_0)(\frak{z}_0)=\frak{z}_0$.
Moreover, $\roman{Ad}(g_0)(\frak{h}_0)$ is another
maximally noncompact Cartan subalgebra of $\frak{s}_0$.
Therefore
there is an inner automorphism of $\frak{s}_0$, and thus also an
element $g_{1}\in\bold{S}_0$
such that 
$g_{1}\circ g_0(\frak{h})=\frak{h}$ 
(see e.g. [S1], [S2], [OV, Ch.4,\S{4.7}]).
Let $g=g_1\circ g_0$. Since $\bold{S}_0\subset\bold{Z}_{\bold{G}}(\frak{z}_0)$,
it suffices to show that $g\in \bold{Z}_{\bold{G}}(\frak{z}_0)$.
We have
$\roman{Ad}(g)(\frak{z}_0)=\frak{z}_0$ and
$\roman{Ad}(g)(\frak{h}_0)=\frak{h}_0$. In particular,
$g\in\bold{N}_{\bold{G}}(\frak{h})$. 
Since $\roman{Ad}(g)$ 
commutes with the conjugation $\sigma$ in
$\Hat{\frak{g}}$,
and $\roman{Ad}(g)
(\frak{q}^n)=\frak{q}^n$, 
$\roman{Ad}(g)(\frak{q}^r)=\frak{q}^r$, we also have\,:
\form{\roman{Ad}(g)(\frak{q}^r\cap\bar{\frak{q}}^{\, r})
=\frak{q}^r\cap\bar{\frak{q}}^{\, r}
\quad\text{and}\quad \roman{Ad}(g)\left(\frak{q}^n
+\left[\bar{\frak{q}}\cap\frak{q}^r
\right]\right)=\frak{q}^n
+\left[\bar{\frak{q}}\cap\frak{q}^r
\right]\,.
}\edef\FFormulaIX{\number\q.\number\t}
The analytic Weyl group 
$\anweyl=\bold{N}_{\bold{G}}(\frak{h})/\bold{H}$
can be identified to a subgroup of the Weyl group $\alweyl=
\bold{N}_{\Hat{\bold{G}}}(\Hat{\frak{h}})/\Hat{\bold{H}}$.
The element $s_g\in\alweyl$ defined by $g$ satisfies\,:
$\roman{Ad}(g)(\Hat{\frak{g}}^{\alpha})\subset\Hat{\frak{g}}^{s_g(\alpha)}$ 
for all $\alpha\in\Cal{R}$. Because of $(\FFormulaIX)$, $g$ normalizes
the complex parabolic subalgebra $\frak{v}=\frak{q}^n+\left(\frak{q}^r
\cap\bar\frak{q}\right)$, with 
$\frak{v}^r=\frak{q}^r\cap\bar{\frak{q}}^{\,r}$. Hence
$s_{g}$ is the composition of symmetries $s_{\alpha}$ with
$\alpha\in\Cal{Q}^r\cap\bar{\Cal{Q}}^{\, r}$. If $\alpha\in\Cal{Q}^r\cap
\bar{\Cal{Q}}^{\, r}$ and $H\in\frak{z}_0$ we obtain
for $X_{\beta}\in\Hat{\frak{g}}^{\beta}$\,:
$$\align
[\roman{Ad}(g)(H),&\,\roman{Ad}(g)(X_{\beta})]=
\roman{Ad}(g)([H,X_{\beta}])\\
&=\beta(H)\, \roman{Ad}(g)(X_{\beta})
=s_{\alpha}(\beta)(H) \roman{Ad}(g)(X_{\beta})
\endalign
$$
because $s_{\alpha}(\beta)(H)=(\beta-(\beta|\alpha^{\vee})\alpha)(H)=\beta(H)$,
since $\alpha(H)=0$ by $(\FRMIZ)$. This shows that 
$[\roman{Ad}(g)(H),X_{\beta}]=\beta(H)X_{\beta}$
for all $\beta\in\Cal{R}$ and hence that $\roman{Ad}(g)(H)=H$\,.
Moreover $\bold{Z}_{\bold{G}}(\frak{z}_0)\subset\bold{N}_{\bold{G}}(\frak{g}_0)$
because $\frak{g}_0$ is the centralizer of $\frak{z}_0$ in $\frak{g}$,
and hence the equality follows.\par
Finally we observe that 
$\bold{Z}_{\bold{G}}(\frak{z}_0)\subset\bold{Q}=
\bold{N}_{\Hat{\bold{G}}}(\frak{q})$.
Indeed, if $g_0\in\bold{Z}_{\bold{G}}(\frak{z}_0)$, 
we can find $g_1\in\bold{S}_0$
such that $g=g_1g_0$ normalizes our Cartan subalgebra $\frak{h}$, still
centralizing $\frak{z}_0$. Being
$\Cal{Q}=\{\alpha\in\Cal{R}\,|\, \alpha(A)\geq 0\}$ for an element
$A\in\frak{h}_{\Bbb{R}}\cap\Hat{\frak{z}_0}$, clearly $\roman{Ad}(g)(\frak{q})=
\frak{q}$. Hence also $\roman{Ad}(g_0)(\frak{q})=\frak{q}$. This shows
that $\bold{Z}_{\bold{G}}(\frak{z}_0)\subset\bold{G}_+$, completing the proof
of the Theorem.
\enddimo
\par\smallskip
Next we prove that our choice of $\bold{G}=\bold{G}_{\Pi(\Cal{R})}$ brings 
that also $\bold{S}_0$ is the semisimple real group associated to the
full weights lattice of $\Hat{\frak{s}}_0$. Indeed we have\,:
\thm{The complexification $\Hat{\bold{S}}_0$ of
the linear Lie group $\bold{S}_0$ is simply connected.}
\edef\TRMID{\number\q.\number\x}
First we prove the following\,:
\lem{Let $E$ be a linear hyperplane of $\frak{h}^*_{\Bbb{R}}$ and let
$\Cal{R}'=E\cap\Cal{R}$. Then every basis of simple roots $\{\alpha_1,
\hdots,\alpha_m\}$ of $\Cal{R}'$ is a subset of a basis of simple
roots $\{\alpha_1\,,\hdots\,,\alpha_m\,,\alpha_{m+1}\,,
\hdots\,,\alpha_{\ell}\}$ of
$\Cal{R}$.}\edef\LMMIJ{\number\q.\number\x}
\dimo Let $A\in\frak{h}_{\Bbb{R}}$ be such that 
$E=\{\eta\in\frak{h}^*_{\Bbb{R}}\,|\, \eta(A)=0\}$. Next we take
a regular element $B\in\frak{h}_{\Bbb{R}}$ such that
$\{\alpha_1,\hdots,\alpha_m\}$ is the system of simple roots in
${\Cal{R}'}^+=\{\alpha\in\Cal{R}'\,|\,\alpha(B)>~0\}$. For small $\epsilon>0$
the element $A+\epsilon B$ is regular and $\alpha_1,\hdots,\alpha_m$
are simple in $\Cal{R}^+=\{\alpha\in\Cal{R}\,|\,\alpha(A+\epsilon B)>0\}$.
Indeed, we take $\epsilon>0$ such that 
$\epsilon\,|\alpha(B)|<|\alpha(A)|$ when $\alpha(A)\neq 0$.
Assume by contradiction that, for some $1\leq i\leq m$, 
the root $\alpha_i$ is not simple in $\Cal{R}^+$\,, i.e. that
we have
$\alpha_i=\beta+\gamma$, with $\beta(A+\epsilon\, B)>0$,
$\gamma(A+\epsilon\,B)>0$, $\beta(B)\geq\gamma(B)$. Since
$\alpha_i$ is simple in ${\Cal{R}'}^+$\,,
we have $\gamma(B)<0<\alpha_i(B)<\beta(B)$. Hence $\gamma(A)>0$.
Moreover $\beta(A)\geq 0$,
because otherwise $\beta(A+\epsilon\,B)<0$ by our choice
of $\epsilon$. Thus we obtain\,:
$$\alpha_i(A)=\beta(A)+\gamma(A)\geq\gamma(A)>0\,,$$
contradicting $\alpha_i\in E$. \enddimo\edef\LMMIF{\number\q.\number\x}
\demo{Proof of Theorem {\TRMID}}
We keep the notation introduced in the previous discussion.
While applying Lemma {\LMMIF} to our situation, we observe that,
if our parabolic $\Cal{Q}$ is $\{\alpha\in\Cal{R}\,|\,
\alpha(A)\geq 0\}$ for some $A\in\frak{h}_{\Bbb{R}}$, then
the set of roots $\alpha$ with $\Hat{\frak{g}}^{\,\alpha}\subset\Hat{\frak{s}}_0$
is $\Cal{R}'=\Cal{Q}^r\cap\bar{\Cal{Q}}^{\,r}=
\{\alpha\in\Cal{R}\,|\, \alpha(A+\epsilon\bar{A})=0\}$, where $\epsilon$ is 
any positive real number with
$\epsilon\,|\bar\alpha(A)|<\alpha(A)$ for
$\alpha\in\Cal{Q}^n$. 
This set $\Cal{R}'$ is naturally isomorphic to the root system 
$\Cal{R}(\Hat{\frak s}_0,\Hat{\frak{h}}_0)$
of the semisimple
complex Lie algebra $\Hat{\frak{s}}_0$ with respect to its
Cartan subalgebra $\Hat{\frak{h}}_0$\,.
We identify $\Pi(\Cal{R}')$ to the set of elements $\omega$ in 
$\langle\Cal{R}'\rangle_{\Bbb{R}}\subset\frak{h}_{\Bbb{R}}^*$
for which $(\omega|\alpha^{\vee})\in\Bbb{Z}$ for all 
$\alpha\in\Cal{R}'$. We have a natural projection
$\varpi:\Pi(\Cal{R})@>>>\Pi(\Cal{R}')$, that is defined by
$(\varpi(\eta)|\alpha^{\vee})=(\eta|\alpha^{\vee})$ for all
$\alpha\in\Cal{R}'$, and coincides with the orthogonal projection
onto $\langle\Cal{R}'\rangle_{\Bbb{R}}$.\par
 The lattice $\Cal{L}=\varpi(\Pi(\Cal{R}))$
satisfies $\Lambda(\Cal{R}')\subset\Cal{L}\subset\Pi(\Cal{R}')$
and is the set of weights of the finite dimensional linear representations
of $\bold{S}_0$.
\par
According to Lemma {\LMMIJ} we can fix
$C\in\frak{C}(\Cal{R},\Cal{Q})$, with 
$\Cal{B}(C)=\{\alpha_1,\hdots,\alpha_{\ell}\}$, in such a way that
$\Cal{B}'=\{\alpha_1,\hdots,\alpha_m\}$, with $m\leq\ell$, is a basis of
$\Cal{R}'$. Let $\Cal{B}^*(C)=\{\mu_1,\hdots,\mu_\ell\}\subset\Pi(\Cal{R})$ 
and ${\Cal{B}'}^*=\{\nu_1,\hdots,\nu_m\}\subset\Pi(\Cal{R}')
\subset\langle\Cal{R}'\rangle_{\Bbb{R}}$ be the
corresponding adjoint basis.
Then $\varpi(\mu_i)=\nu_i$
for $i=1,\hdots,m$. This shows that actually $\Cal{L}=\Pi(\Cal{R}')$,
proving our statement.
\enddimo
Now we give a more accurate description of the group $\pi_0(\bold{G}_+)$ 
of the connected components of $\bold{G}_+$.
\thm{Let $(\frak{g},\frak{q})$ be an effective parabolic $CR$ algebra.
Keep the notation of $(\FormulaIN)$,
$(\FormulaIO)$, $(\FormulaIP)$ and $(\FormulaIQ)$. 
Then 
$\bold{N}_{\bold{S}_0}(\frak{h}_0)$ is a 
closed normal Lie subgroup of 
$\bold{N}_{\bold{G}}(\frak{h})\cap\bold{Z}_{\bold{G}}(\frak{z}_0)$.
The natural inclusion
$\bold{N}_{\bold{G}}(\frak{h})\cap\bold{Z}_{\bold{G}}(\frak{z}_0)
\hookrightarrow
{\bold{G}}_+$
passes to the quotient to define a one-to-one 
correspondence of connected components\,:
\form{
\pi_0\left(\bold{N}_{\bold{G}}(\frak{h})\cap\bold{Z}_{\bold{G}}(\frak{z}_0)
/\bold{N}_{\bold{S}_0}(\frak{h})\right)
\;\longleftrightarrow \pi_0\left({\bold{G}}_+\right)\,.}}
\dimo\edef\FFormulaIV{\number\q.\number\t}
The inclusion $\bold{G}_0\hookrightarrow\bold{G}_+$ defines an isomorphism
$\pi_0(\bold{G}_0) @>{\sim}>>\pi_0(\bold{G}_+)$.
\par We showed in the proof of Proposition {\TRMIH} that
every connected component of $\bold{G}_0$ contains an element of
$\bold{N}_{\bold{G}}(\frak{h})\cap\bold{G}_0$
and that 
$\bold{N}_{\bold{G}}(\frak{h})\cap\bold{G}_0=\bold{N}_{\bold{G}}(\frak{h})
\cap \bold{Z}_{\bold{G}}(\frak{z}_0)$.
\par\smallskip
Since $\frak{s}_0$ is an ideal of 
$\frak{g}_0=\frak{s}_0\oplus\frak{z}_0$, the analytic subgroup
$\bold{S}_0$ is
normal in $\bold{G}_0$ and we already noticed that is
contained in
$\bold{Z}_{\bold{G}}(\frak{z}_0)$. Moreover,
$\bold{N}_{\bold{S}_0}(\frak{h})=
\bold{N}_{\bold{G}}(\frak{h})\cap\bold{S}_0\,.$ \par
Let $\bold{T}_0$ be the analytic Lie subgroup of $\bold{G}_0$
with Lie algebra $\frak{z}_0$. Then $\bold{S}_0\bowtie
\bold{T}_0$ is the connected component of the identity in
$\bold{G}_0$. We have\,: 
$$\bold{T}_0\subset\bold{H}
\subset \bold{N}_{\bold{G}}(\frak{h})\cap\bold{Z}_{\bold{G}}(\frak{z}_0)\,
.$$
Hence we obtain isomorphisms\,:
\form{\CD
\dsize\frac{\bold{N}_{\bold{G}}(\frak{h})\cap\bold{Z}_{\bold{G}}(
\frak{z}_0)}{\bold{N}_{\bold{S}_0}(\frak{h})\bowtie \bold{T}_0}
@>{\dsize\cong}>> \dsize\frac{\bold{G}_0}{\bold{S}_0\bowtie\bold{T}_0}
@>{\dsize\cong}>> \pi_0(\bold{G}_+)\,,\endCD}
yielding the isomorphism in $(\FFormulaIV)$.\enddimo
\edef\TeoremaIB{\number\q.\number\x}
We denote by $\bold{H}_0$ the Cartan subalgebra of $\bold{S}_0$
corresponding to $\frak{h}_0$. Since 
$\bold{S}_0\subset\bold{Z}_{\bold{G}}(\frak{z}_0)$, we have\,:
\form{\bold{H}_0=\{g\in\bold{S}_0\,|\,\roman{Ad}(g)(H)=H\,,\;
\forall H\in\frak{h}_0\,\}=\bold{H}\cap\bold{S}_0\,.}
\edef\FRMLIBH{\number\q.\number\t}
By $(\FRMLIBH)$, 
 we obtain for the analytic Weyl group $\bold{W}(\frak{s}_0,\frak{h}_0)
=\bold{N}_{\bold{S}_0}(\frak{h}_0)/\bold{H}_0$\,:
\form{\bold{W}(\frak{s}_0,\frak{h}_0)\cong 
\bold{N}_{\bold{S}_0}(\frak{h})
/\bold{Z}_{\bold{S}_0}(\frak{h})=
\bold{N}_{\bold{S}_0}(\frak{h})
/\left(\bold{H}\cap\bold{S}_0\right)
\,.}
Thus we can identify $\bold{W}(\frak{s}_0,\frak{h}_0)$ with a
subgroup of the centralizer 
$\bold{Z}_{\bold{W}(\frak{g},\frak{h})}(\frak{z}_0)$
of $\frak{z}_0$ in $\bold{W}(\frak{g},\frak{h})$.
Since $\bold{S}_0$ is normal in $\bold{G}_0$, it turns out
that $\bold{W}(\frak{s}_0,\frak{h}_0)$ is normal in
$\bold{Z}_{\bold{W}(\frak{g},\frak{h})}(\frak{z}_0)$.
We have\,:
\thm{Keep the notation and assumptions of Theorem {\TeoremaIB}.
We have a short exact sequence of groups and homomorphisms\,:
\form{
\bold{1}@>>>\frac{\bold{H}\;\,}{\bold{H}_0}
@>>>
\frac{\bold{N}_{\bold{G}}(\frak{h})\cap\bold{Z}_{\bold{G}}(\frak{z}_0)}{
\bold{N}_{\bold{S}_0}(\frak{h}_0)}
@>>>\frac{\bold{Z}_{\bold{W}(\frak{g},\frak{h})}(\frak{z}_0)}{
\bold{W}(\frak{s}_0,\frak{h}_0)}
@>>>\bold{1}\,.}
We obtain an exact sequence of finite groups\,:
\form{
\bold{1}@>>>\pi_0\left(
\frac{\bold{H}\;\,}{\bold{H}_0}
\right) @>>>\pi_0\left(
\bold{G}_+
\right)
@>>>{\frac{\bold{Z}_{\bold{W}(\frak{g},\frak{h})}(\frak{z}_0)}{
\bold{W}(\frak{s}_0,\frak{h}_0)}}
@>>>\bold{1}\,.}
When $\frak{h}$ is maximally compact in $\frak{g}_+$,
then $\bold{H}_0$ 
is connected and $(\number\q.\number\t)$ yields the
exact sequence\,:
\form{
\bold{1}@>>>\pi_0\left(
\bold{H}\right) @>>>\pi_0\left(
\bold{G}_+\right)
@>>> \frac{\bold{Z}_{\bold{W}(\frak{g},\frak{h})}(\frak{z}_0)}{
\bold{W}(\frak{s}_0,\frak{h}_0)}
@>>>\bold{1}\,.}
When $\frak{h}$ is maximally noncompact in $\frak{g}_+$, then
$\bold{Z}_{\bold{W}(\frak{g},\frak{h})}(\frak{z}_0)\cong
\bold{W}(\frak{s}_0,\frak{h}_0)$
and \advance\t by-1
$(\number\q.\number\t)$ yields the isomorphism\,:\advance\t by1
\form{\pi_0\left( 
\frac{\bold{H}\;\,}{\bold{H}_0}
\right) \cong \pi_0\left(
\bold{G}_+\right)\,.}}
\edef\FormulaKD{\number\q.\number\t}
\advance\t by-1
\edef\FormulaKC{\number\q.\number\t}
\advance\t by-1
\edef\FormulaKB{\number\q.\number\t}
\advance\t by-1
\edef\FormulaKA{\number\q.\number\t}
\advance\t by3 
\dimo \edef\TRMAZ{\number\q.\number\x}
Since $\bold{Z}_{\bold{W}(\frak{g},\frak{h})}(\frak{z}_0)$
is a finite group, $(\FormulaKB)$ is a consequence of
the exact homotopy sequence of a fiber bundle, of
$(\FormulaKA)$
and of Theorem {\TeoremaIB}. Thus a proof of $(\FormulaKA)$ also 
provides a proof of $(\FormulaKB)$.
\par
We observe that $\bold{H}$ is a normal subgroup
of $\bold{N}_{\bold{G}}(\frak{h})\cap\bold{Z}_{\bold{G}}(\frak{z}_0)$
and that the quotient
$\left(\bold{N}_{\bold{G}}(\frak{h})\cap\bold{Z}_{\bold{G}}(\frak{z}_0)
\right)/\bold{H}$ is naturally isomorphic to the centralizer
$\bold{Z}_{\bold{W}(\frak{g},\frak{h})}(\frak{z}_0)$ of
$\frak{z}_0$ in the analytic Weyl group $\bold{W}(\frak{g},\frak{h})$.
Moreover, $\bold{N}_{\bold{S}_0}(\frak{h}_0)$ is contained
in $\bold{Z}_{\bold{G}}(\frak{z}_0)$ and
$\bold{N}_{\bold{S}_0}(\frak{h}_0)\cap
\bold{H}=\bold{H}_0$.
Thus $(\FormulaKA)$ follows from the general homomorphism theorems
of groups.\par\smallskip
If $\frak{h}$ is maximally compact, then $\frak{h}_0$ is
a maximally compact Lie subalgebra of $\frak{s}_0$. Then (see e.g.
[Kn, Proposition 7.90, p.488]) the Cartan subalgebra
$\bold{H}_0$ is connected and
$(\FormulaKC)$ follows from $(\FormulaKB)$.
\par
If $\frak{h}$ is maximally noncompact, then, by
[Kn, Proposition 7.90, p.488] the Cartan subalgebra
$\bold{H}$ intersects
every connected component of $\bold{G}_+$. This implies that
the map $\pi_0(\bold{H}/\bold{H}_0)@>>>\pi(\bold{G}_+)$ 
in $(\FormulaKB)$ is onto and hence $(\FormulaKD)$ holds true and
$\bold{Z}_{\bold{W}(\frak{g},\frak{h})}(\frak{z}_0)
\cong\bold{W}(\frak{s}_0,\frak{h}_0)$.\enddimo
We describe now $\pi_0(\bold{G}_+)$ in terms of characters\,:
\prop{Let $\frak{h}$ be a maximally noncompact Cartan
subalgebra of $\frak{g}_+$. 
Fix a Weyl chamber $C\in\frak{C}(\Cal{R})$ such that
$\Cal{B}(C)=\{\alpha_1,\hdots,\alpha_{\ell}\}$ contains a basis
$\Cal{B}_0(C)=\{\alpha_1,\hdots,\alpha_m\}$ of simple roots of
$\Cal{Q}^r\cap\bar{\Cal{Q}}^{\,r}$. Let
$E=\langle \alpha_{m+1},\hdots,\alpha_{\ell}
\rangle^\perp$. Then\,:
\form{
\pi_0(\bold{G}_+)\simeq
\pi_0\left(\{\chi\in\Hom^{\sigma}(\Pi(\Cal{R}),\Bbb{C}^*)\,|\,
\chi(\mu)=1\,,\;\forall \mu\in\Pi(\Cal{R})\cap E\}\right)\,.
}}\edef\RFMLICD{\number\q.\number\t}
\dimo \edef\PRPII{\number\q.\number\x}
 We denote by $\varpi_0$ the 
orthogonal projection onto $\langle \alpha_{m+1},\hdots,\alpha_{\ell}
\rangle^\perp$. Let $\Cal{B}^*(C)=\{\omega_1,\hdots,\omega_{\ell}\}$
be the adjoint basis of $\Cal{B}(C)$. Then we have
$\varpi_0(\eta)=\sum_{i=1}^m{(\eta|\alpha_i^{\vee})\omega_i}$\,.
We can identify $\Cal{R}'=\Cal{R}(\Hat{\frak{s}}_0,\Hat{\frak{h}}_0)$
with $\varpi_0(\Cal{Q}^r\cap\bar{\Cal{Q}}^{\,r})$. The advantage of this
identification is that the corresponding weight space $\Pi(\Cal{R}')$
becomes a subspace of $\Pi(\Cal{R})$, and moreover\,:
\form{\Pi(\Cal{R}')\,=\,\varpi_0(\Pi(\Cal{R}))\,=\,
\Pi(\Cal{R})\cap \langle \alpha_{m+1},\hdots,\alpha_{\ell}
\rangle^\perp\,.}
In this way,
for each character $\chi\in\Hom^{\sigma}(\Pi(\Cal{R}),\Bbb{C}^*)$,
the composition 
$\chi\circ\varpi_0$ is still a character 
of $\Hom^{\sigma}(\Pi(\Cal{R}),\Bbb{C}^*)$, 
whose restriction to $\Pi(\Cal{R}')$
is a character
in $\Hom^{\sigma}(\Pi(\Cal{R}'),
\Bbb{C}^*)$, and we have\,:
\form{\bold{H}_0=
\{h_{\chi}\,|\, \chi\in\Hom^{\sigma}(\Pi(\Cal{R}'),
\Bbb{C}^*)\}=\{h_{\chi\circ\varpi}\,|\,\chi\in\Hom^{\sigma}(\Pi(\Cal{R}),
\Bbb{C}^*)\}\,.}
Thus, setting\,:
\form{\bold{H}_1=\{h_{\chi}\,|\, \chi\in\Hom^{\sigma}(\Pi(\Cal{R}),
\Bbb{C}^*)\,,\; \chi(\eta)=1\;\forall\eta\in\Pi(\Cal{R}')\}\,,}
the map\,:
\form{\bold{H}\ni h_{\chi} @>>> h_{\chi\circ\varpi}^{-1}\circ h_{\chi}
\in \bold{H}_1}
yields, by passing to the quotient, the isomorphism
$\bold{H}/\bold{H}_0@>{\sim}>> \bold{H}_1$, which implies  (\RFMLICD).\enddimo
\cor{Let $(\frak{g},\frak{q})$ be an effective parabolic totally real
$CR$ algebra. Let
$\frak{h}$ be a maximally noncompact 
Cartan subalgebra of $\frak{g}$ contained in $\frak{g}_+$ and let
$C\in\frak{C}(\Cal{R},\Cal{Q})$ be S-fit for $(\frak{g},\frak{q})$
and hence S-adapted to the conjugation $\sigma$. 
Then $\pi_0(\bold{G}_+)$ is isomorphic to $\Bbb{Z}_2^e$, where
$e$ is the number of roots in $\Phi_C\cap\rootre$.}
\edef\CORGH{\number\q.\number\x}
\dimo 
By inspecting the Satake diagrams of the simple real Lie algebras,
we obtain that $\chi(\omega_i)$ is a positive real number when
$\omega_i$ is a real weight corresponding to a complex root 
$\alpha_i\in\Cal{B}(C)$. Then the statement follows from
Proposition {\PRPII}.\par
Another independent proof can be given, using [DKV] and [Wi]. Because
$\bold{G}$ is connected, we
have an exact sequence\,:
\form{\CD
\cdots @>>>\pi_1(\bold{G}/\bold{G}_+) @>>>\pi_0(\bold{G}_+)@>>>\bold{1}\,.
\endCD}\edef\FRMLISE{\number\q.\number\t}
Take $\Phi'_C=\Phi_C\setminus\rootre$. Then also 
$\frak{q}'=\frak{q}_{\Phi'_C}$ is the complexification of a real parabolic
subalgebra $\frak{g}'_+$ of $\frak{g}$. By [Wi] and (\number\q.\number\t),
the stabilizer $\bold{G}_+'$ of $\frak{q}'$ in $\bold{G}$ is connected.
By Proposition {\PRPII}, this implies that 
$\chi(\omega_i)>0$ on the real $\omega_i\in\Cal{B}^*(C)$ corresponding
to a complex root $\alpha_i$ in $\Cal{B}(C)$. Again, 
the statement follows from
Proposition {\PRPII}.
\enddimo
\par\edef\CRLIJ{\number\q.\number\x}

For each $\omega\in\Pi(\Cal{R})$ there is, modulo isomorphisms, a unique 
irreducible finite dimensional complex linear representation 
$\rho_{\omega}:\Hat{\frak{g}}@>>>\frak{gl}_{\Bbb{C}}(V_{\omega})$, for which 
$\omega$ is an extremal weight. For each weight $\eta$ we denote by\,:
\form{
V_{\omega}^{\eta}=\{v\in V_{\omega}\,|\, \rho_{\omega}(H)(v)=\eta(H)\,v\,,\;
\forall H\in\Hat{\frak{h}}\}}
the eigenspace of $V_{\omega}$ corresponding to the
weight $\eta$\,.
Let\,: 
\form{\Hat{\frak{g}}_{\omega}=\{Z\in\Hat{\frak{g}}\,|\, \rho_{\omega}(
V_{\omega}^{\omega})\subset V_{\omega}^{\omega}\}\,.}
Since the eigenspace
$V_{\omega}^{\omega}$ of $\omega$ is one-dimensional, for each
$Z\in\Hat{\frak{g}}_{\omega}$ there is a unique complex number,
that we shall denote by $\omega(Z)$, such that\,:
\form{\rho_{\omega}(Z)(v)=\omega(Z)\,v\,,\qquad
\forall Z\in\Hat{\frak{g}}_{\omega}\,,
\;\forall v\in V_{\omega}^{\omega}\,.}
This agrees with the natural definition of $\omega(Z)$
by the duality pairing when
$Z\in\Hat{\frak{h}}$.
\prop{Let $\omega_1\,,\,\omega_2\in\Pi(\Cal{R})$ and $\omega=\omega_1+\omega_2$.
Then\,: 
\form{
\Hat{\frak{g}}_{\omega_1}\cap\Hat{\frak{g}}_{\omega_2}\subset
\Hat{\frak{g}}_{\omega}} and\,:
\form{\omega(Z)=\omega_1(Z)+\omega_2(Z)\,,\;\forall Z\in 
\Hat{\frak{g}}_{\omega_1}\cap\Hat{\frak{g}}_{\omega_2}\,.}}
\edef\FRMKD{\number\q.\number\t}
\advance\t by-1
\edef\FRMKC{\number\q.\number\t}
\advance\t by1
\dimo There is an injective homomorphism of $\Hat{\frak{g}}$-modules
$V_{\omega}@>>>V_{\omega_1}\otimes V_{\omega_2}$ that maps
$V_{\omega}^{\omega}$ into $V_{\omega_1}^{\omega_1}\otimes
V_{\omega_2}^{\omega_2}$. For $Z\in\Hat{\frak{g}}$, we have
$\left(\rho_{\omega_1}\otimes\rho_{\omega_2}\right)(v_1\otimes v_2)=
\rho_{\omega_1}(Z)(v_1)\otimes v_2 +
v_1\otimes\rho_{\omega_2}(Z)(v_2)$ for all $v_1\in V_{\omega_1}$,
$v_2\in V_{\omega_2}$. This implies
$(\FRMKC)$ and $(\FRMKD)$\,.\enddimo
We also have\,:
\lem{Let $\omega\in\Pi(\Cal{R}(\Hat{\frak{g}},\Hat{\frak{h}}))$ and
$\Hat{\frak{h}}'$ another 
Cartan subalgebra of $\Hat{\frak{g}}$ contained
in $\Hat{\frak{g}}_{\omega}$. 
Let $V_{\omega}$ be an irreducible $\Hat{\frak{g}}$-module
with extremal weight $\omega$. Then\,:
\roster
\item"$(i)$" There is a weight $\omega'\in\Pi(\Cal{R}(\Hat{\frak{g}},
\Hat{\frak{h}}'))$ that is an extremal weight for $V_{\omega}$\,;
\item"$(ii)$" we can choose 
$\omega'\in\Pi(\Cal{R}(\Hat{\frak{g}},\Hat{\frak{h}}'))$ 
in such a way that 
$V_{\omega}^{\omega}=V_{\omega}^{\omega'}$\,;
\item"$(iii)$" $\Hat{\frak{g}}_{\omega}=\Hat{\frak{g}}_{\omega'}$ and
$\omega(Z)=\omega'(Z)\quad\forall Z\in\Hat{\frak{g}}_{\omega}\,.$\qed
\endroster}\edef\LMMAKB{\number\q.\number\x} 
Since we took a simply connected $\Hat{\bold{G}}$,
for each $\omega\in\Pi(\Cal{R})$, the representation 
$\rho_{\omega}:\Hat{\frak{g}}@>>>\frak{gl}_{\Bbb{C}}(V_{\omega})$
lifts to a representation 
$\rho_{\omega}:\Hat{\bold{G}}@>>>\bold{GL}_{\Bbb{C}}(V_{\omega})$.
The stabilizer of the line $V_{\omega}^{\omega}$ in $\Hat{\bold{G}}$
contains the Cartan subgroup $\Hat{\bold{H}}$ of $\Hat{\bold{G}}$ and
therefore
is the closed connected subgroup $\Hat{\bold{G}}_{\omega}$ of
$\Hat{\bold{G}}$ with Lie algebra $\Hat{\frak{g}}_{\omega}$. 
The map $\omega:\Hat{\frak{g}}_{\omega}@>>>\Bbb{C}$ also lifts to
a character $\varphi_{\omega}:\Hat{\bold{G}}_{\omega} @>>>\Bbb{C}^*$,
with $\varphi_{\omega}(\exp(Z))=\exp(\omega(Z))$ for all 
$Z\in\Hat{\frak{g}}_{\omega}$\,.\par
We observe that every parabolic $\frak{q}$ containing $\Hat{\frak{h}}$
is equal to some
$\Hat{\frak{g}}_{\omega}$, for a suitable choice of 
$\omega\in \Pi(\Cal{R})$. This choice is not unique. It can be done
in the following way. First we choose a Weyl chamber 
$C\in\frak{C}(\Cal{R},\Cal{Q})$. Let $\frak{q}=\frak{q}_{\Phi_C}$ for
$\Phi_C\subset\Cal{B}(C)$. Let $\Cal{B}^*(C)$ be the adjoint basis
in $\Pi(\Cal{R})$ and 
$\Phi^*(C)=\{\omega_1,\hdots,\omega_k\}$ 
the set of elements of $\Cal{B}^*(C)$
that vanish on $\Cal{B}(C)\setminus\Phi_C$. Then for any $k$-tuple of
strictly positive integers $(t_1,\hdots,t_k)$ we have\,:
\form{\frak{q}=\Hat{\frak{g}}_{t_1\omega_1+\cdots+t_k\omega_k}=
\bigcap_{1\leq i\leq k}{\Hat{\frak{g}}_{\omega_i}}\,.}
\prop{Let $\bold{Q}$ be the parabolic subgroup of $\Hat{\bold{G}}$
of the parabolic complex Lie subalgebra $\frak{q}$. Then, with
the notation above\,:
\form{\bold{Q}=\bigcap_{1\leq i\leq k}{\Hat{\bold{G}}_{\omega_i}}
\quad\text{and}\quad \Hom_{\Bbb{C}}(\bold{Q},\Bbb{C}^*)=\{\varphi_{\omega}\,|\,
\omega\in\Pi(\Cal{R})\cap{\langle{\Cal{Q}^r}\rangle}_{
\Bbb{R}}^\perp\}}
is the free Abelian group generated by $\varphi_{\omega_1}$, $\hdots$,
$\varphi_{\omega_k}$\,.}
\dimo The inclusion $\{\varphi_{\omega}\,|\,
\omega\in\Pi(\Cal{R})\cap\langle{\Cal{Q}^r}\rangle_{
\Bbb{R}}^\perp\}\subset\Hom_{\Bbb{C}}(\bold{Q},\Bbb{C}^*)$ is a consequence
of the previous discussion. \par
To prove the opposite inclusion, we fist observe that a character
in $\Hom_{\Bbb{C}}(\bold{Q},\Bbb{C}^*)$ restricts to a character in
$\Hom_{\Bbb{C}}(\Hat{\bold{H}},\Bbb{C}^*)$, and hence is of the form
$\varphi_{\omega}$ for some $\omega\in\Pi(\Cal{R})$. Moreover, 
being connected, $\bold{Q}$ admits a Levi decomposition
$\bold{Q}=\bold{Q}_{\roman{ss}}\bold{Q}_{\roman{rad}}$ and
$\varphi_{\omega}(g)=1$ for $g\in\bold{Q}_{\roman{ss}}$($=$ the Levi subgroup of
$\bold{Q}$). This yields $(\omega|\alpha^{\vee})=0$ for
$\alpha\in\Cal{Q}^r$.
\enddimo 
\lem{If $\omega\in\Pi(\Cal{R})$ is a real weight, then 
$\varphi_{\omega}$ is real valued on $\Hat{\bold{G}}_{\omega}\cap\bold{G}$.}
\dimo From $\varphi_{\omega}(\exp(Z))=\exp(\omega(Z))$ for 
$Z\in\Hat{\frak{g}}_{\omega}$ we obtain that $\varphi_{\omega}(g)$ is real
when $g$ belongs to a neighborhood of the identity in 
$\Hat{\bold{G}}_{\omega}\cap\bold{G}$. Since $\varphi_{\omega}(g)$ is
an algebraic function of $g$, the statement follows because
any neighborhood of the identity is Zariski dense in
$\Hat{\bold{G}}_{\omega}\cap\bold{G}$.\enddimo
Let $(\vartheta,\frak{h})$ be a Cartan pair adapted to the effective
parabolic $CR$ algebra $(\frak{g},\frak{q})$, choose
 $C\in\frak{C}(\Cal{R},\Cal{Q})$ and set\,:
$$\Phi^*_C =\{\omega_1,\hdots,\omega_k\}=
\Cal{B}^*(C)\cap\left(\Cal{Q}^r\cap\bar{\Cal{Q}}^{\,r}\right)^\perp\,.$$
Then $\bold{G}_+\subset\bigcap_{1\leq i\leq k}{\Hat{\bold{G}}_{\omega_i}}$,
and we can define the map\,:
\form{\boldsymbol{\varphi}:\bold{G}_+\ni g @>>> (\varphi_{\omega_1}(g),\hdots,
\varphi_{\omega_k}(g))\in\left(\Bbb{C}^*\right)^k\,.}
We have\,:\edef\FRMKG{\number\q.\number\t}
\thm{The map $(\FRMKG)$ yields, by passing to the quotients, a group
isomorphism\,: 
$\pi_0(\bold{G}_+) @>{\boldsymbol{\varphi}_*}>>\pi_0(\boldsymbol{\varphi}(\bold{G}_+))$.}
\dimo By using Lemma {\LMMAKB}, we can restrain to proving the Theorem
in the case $\frak{h}$ is maximally noncompact among the admissible
Cartan subalgebras of $\frak{g}$ contained in $\frak{g}_+$. 
Moreover, by passing to the weakest $CR$ model, we can as well assume
that $\Cal{Q}^r=\bar{\Cal{Q}}^{\,r}$. 
By Proposition {\PRPII}, we obtain a commutative diagram\,:
\form{\CD
\pi_0(\bold{H}/\bold{H}_0) @>{\sim}>> \pi_0(\bold{G}_+)\\
@V{\wr}VV        @VV{\boldsymbol{\varphi}_*}V\\
\pi_0(\boldsymbol{\varphi}(\bold{H})) @>{\sim}>>\pi_0(\boldsymbol{\varphi}(\bold{G}_+))
\endCD}
showing that also $\boldsymbol{\varphi}_*$ is an isomorphism. \enddimo
In particular, when $(\frak{g},\frak{q})$ is totally real,
i.e. $\frak{q}=\bar{\frak{q}}$, we obtain,
by using Corollary {\CRLIJ}\,:
\thm{Let $(\frak{g},\frak{q})$ be a totally real parabolic $CR$ 
algebra, $\frak{h}$ a maximally noncompact Cartan subalgebra of
$\frak{g}$ contained in $\frak{g}_+$, $C\in\frak{C}(\Cal{R},\Cal{Q})$
an S-fit (and S-adapted) Weyl chamber. Set\,:
$$
\Cal{B}(C)=\{\alpha_1,\hdots,\alpha_{\ell}\}\,,\quad
\{\alpha_1,\hdots,\alpha_\mu\}=\Phi_C\cap\rootre\,,\quad
\Cal{B}^*(C)=\{\omega_1,\hdots,\omega_{\ell}\}\,.
$$
Then $\varphi_{\omega_i}(g)\in\Bbb{R}^*$ for $1\leq i\leq\mu$.
The map\,:
\form{\boldsymbol{\varphi}^\flat\colon\bold{G}_+\ni g
@>>>\left(\frac{\varphi_{\omega_1}(g)}{|\varphi_{\omega_1}(g)|},
\hdots,\frac{\varphi_{\omega_\mu}(g)}{|\varphi_{\omega_\mu}(g)|}\right)
\in\{-1,1\}^{\mu}=\Bbb{Z}_2^\mu} 
defines, by passing to the quotient, an isomorphism\,:
\form{\boldsymbol{\varphi}^\flat_*\colon
\pi_0(\bold{G}_+)\ni [g] @>>>\boldsymbol{\varphi}^{\flat}(g)
\in\Bbb{Z}_2^\mu\,.\hfill\qed}}

\advance\t by-1
\edef\FRMLQ{\number\q.\number\t}
\advance\t by1
\edef\TRMKE{\number\q.\number\x}
\thm{We keep the notation and assumptions of Theorem {\TRMKE}.
For each $\alpha\in\Phi_C\cap\rootre$, let 
$\bold{S}_{\alpha}$ be the simple analytic real subgroup
of $\bold{G}$, with
Lie algebra $\frak{s}_{\alpha}=
\frak{g}\cap\left(\Hat{\frak{g}}^{\alpha}\oplus
\Hat{\frak{g}}^{-\alpha}\oplus [\Hat{\frak{g}}^{\alpha},
\Hat{\frak{g}}^{-\alpha}]\right)$.
Let  
$t^*_{\alpha}$ be a generator of $\pi_1(\bold{S}_{\alpha})\simeq\Bbb{Z}$,
and let $t_{\alpha}$ be its image in $\pi_1(\bold{G}/\bold{G}_+)$\,.
Then $\{t_{\alpha}\,|\, \alpha\in\Phi_C\cap\rootre\}$ is a set of
generators of $\pi_1(\bold{G}/\bold{G}_+)$, and the 
homomorphism 
$\boldsymbol{\delta}:\pi_1(\bold{G}/\bold{G}_+)@>>>\pi_0(\bold{G}_+)$ of
the exact homotopy sequence
$(\FRMLISE)$ of the principal bundle $\bold{G}@>>>\bold{G}/\bold{G}_+$
is given by\,:
\form{\boldsymbol{\delta}\colon \pi_1(\bold{G}/\bold{G}_+)\ni t_{\alpha_i} @>>> 
(1,\hdots,\underset{i}\to{\underbrace{-1}},
\hdots,1)\in\Bbb{Z}^\mu_2\simeq\pi_0(\bold{G}_+)\,.}}
\edef\FRMLKK{\number\q.\number\t}\edef\TRMKF{\number\q.\number\x}
\dimo By using results from [DKV], in [Wi, \S{2}] it is shown that
the fundamental group of $\bold{G}/\bold{G}_+$ is generated by
the images of the generators of the $\pi_1(\bold{S}_{\alpha})$'s,
for simple roots $\alpha\in\Phi_C\cap\rootre$. In fact the 
generalization of the Bruhat
decomposition in [DKV]
yields a cell decomposition of $\bold{G}/\bold{G}_+$
where the open $1$-cells correspond to roots 
$\alpha\in\Phi_C\cap\rootre$.
Let $\rho_{\alpha}:\frak{sl}(2,\Bbb{R})@>>>\frak{s}_{\alpha}$ be the
representation with\,:
$$\pmatrix
0&0\\
1&0
\endpmatrix @>{\rho_{\alpha}}>> X_{-\alpha}\,,\;
\pmatrix
0&-1\\
0&0
\endpmatrix @>{\rho_{\alpha}}>> X_{\alpha}\,,\;
\pmatrix
1&0\\
0&-1
\endpmatrix @>{\rho_{\alpha}}>> H_{\alpha}\,.
$$
We still denote by $\rho_{\alpha}$ the corresponding group representation
$\bold{SL}(2,\Bbb{R})@>>> \bold{S}_{\alpha}$.
The open $1$-cell corresponding to $\alpha$ can be parametrized by\,:
$$
\left.\left\{g(t)=\rho_{\alpha}\left[\left(\smallmatrix
\tan t&1\\
-1&0\endsmallmatrix\right)
\right]\,\right|\, -{\pi}/{2}<t<{\pi}/{2}\right\}\,.
$$
Since $g_+(t)=\rho_{\alpha}\left[
\left(\smallmatrix
\cos t &-\sin t\\
0&{1}/{\cos t}\endsmallmatrix\right)
\right]\in\bold{G}_+$ for all $|t|<{\pi}/{2}$, the closure in
$\bold{G}/\bold{G}_+$ of this $1$-cell is the 
loop that is the image in $\bold{G}/\bold{G}_+$ of\,:
$$\left[-{\pi}/{2},{\pi}/{2}\right]\ni t @>>>
g_0(t)=g(t)g_+(t)=\rho_{\alpha}\left[
\left(\smallmatrix
\sin t&\cos t\\
-\cos t&\sin t
\endsmallmatrix\right)
\right]\,,$$
that coincides with a maximal compact torus of $\bold{S}_{\alpha}$.
We note that $g_0(\pi/2)=e$, while $g_0(-\pi/2)=\exp(i\pi H_{\alpha})$.
Hence $\varphi_{\omega_j}(g_0(\pi/2))=1$, and\,:
$$\varphi_{\omega_j}(g_0(-\pi/2))=\exp(i\pi\omega_j(H_{\alpha}))=
\exp(i\pi(\omega_j|\alpha^{\vee}))=(-1)^{(\omega_j|\alpha^{\vee})}\,.$$
This proves $(\FRMLKK)$. \enddimo


\se{The fundamental group of parabolic $CR$ manifolds}
In this section we give an explicit combinatorial description of 
the fundamental group of parabolic $CR$ manifolds.
We keep the notation of the previous sections.
By using 
Proposition {\ProposizioneLJ} and Theorem {\TRMKE}, we obtain\,:
\thm{Let $(\frak{g},\frak{q})$ be an effective parabolic $CR$ algebra
and let $M=M(\frak{g},\frak{q})$ be the corresponding homogeneous
$CR$ manifold. Then there exists a totally real parabolic $CR$ algebra
$(\frak{g},\frak{q}')$ such that 
$$\bold{G}_+=\{g\in\bold{G}\,|\,\roman{Ad}_{\Hat{\frak{g}}}(g)(\frak{q})
\subset\frak{q}\}\subset
\bold{G}_+'=\{g\in\bold{G}\,|\,\roman{Ad}_{\Hat{\frak{g}}}(g)(\frak{q}')
\subset\frak{q}'\}$$
and the $\bold{G}$-equivariant map
$$f:M =\bold{G}/\bold{G}_+@>>> M'=M(\frak{g},\frak{q}')=
\bold{G}/\bold{G}_+'$$
has simply connected complex fibers. With
$F=\bold{G}'_+/\bold{G}_+$, we obtain exact sequences\,:
\form{\CD
\bold{1} @>>>\pi_0(\bold{G}_+)
@>>>\pi_0(\bold{G'}_+)@>>>\pi_0(F)@>>>\bold{1}\endCD\,,}
 \form{\CD
\bold{1}@>>>\pi_1(M)@>>>\pi_1(M')@>>>\pi_0(F)@>>>\bold{1}\,.
\endCD}
The induced map in homotopy $f_*:\pi_1(M)@>>>\pi_1(M')$
is injective and $f_*(\pi_1(M))$ is a normal subgroup with
finite index in $\pi_1(M')$.}
\dimo \edef\TEOMHA{\number\q.\number\x}
\edef\FRMLHB{\number\q.\number\t}
\advance\t by-1
\edef\FRMLHA{\number\q.\number\t}
\advance\t by1
Let $\frak{q}'$ be the parabolic $\frak{q}_m$ of 
Proposition {\ProposizioneLJ}. We proved that the typical fiber
$F=\bold{G}'_+/\bold{G}_+$ is simply connected. Thus
$(\FRMLHA)$ and $(\FRMLHB)$ are consequences of  Serre's long
exact sequence for the homotopy groups of a fiber bundle. 
By Theorem {\TRMKE}, $\pi_0(\bold{G}'_+)$ is 
a finite Abelian group.
Hence $\pi_0(\bold{G}_+)$ is normal in $\pi_0(\bold{G}'_+)$ and 
$\pi_0(F)\simeq \pi_0(\bold{G}'_+)/\pi_0(\bold{G}_+)$ has a natural
structure of Abelian group. Then also the maps in $(\FRMLHB)$ are
group homomorphisms and therefore $f_*(\pi_1(M))$ is a normal
subgroup of $\pi_1(M')$.
\enddimo
\cor{If $\frak{g}$ is a simple real Lie algebra of the 
complex type or of real type
$\roman{A\,I\! I}$, $\roman{A\,I\!I\! Ia}$,
$\roman{A\,I\! V}$, $\roman{B\,I\! I}$, $\roman{C\,I\! I}$,
$\roman{D\,I\! I}$, $\roman{D\,I\! I\! Ib}$, $\roman{E\,I\! I\! I}$,
$\roman{E\,I\! V}$, $\roman{F\,I\! I}$,
then all orbits $M=M(\frak{g},\frak{q})$ are simply connected.}

\dimo The multiplicities of the simple real roots of $\Hat{\frak{g}}$ are
always different from one, hence $M'=M(\frak{g},\frak{q}')$ is simply 
connected, forcing $M$ to be simply connected.\enddimo

Note that the $M(\frak{g},\frak{q}')$ for a
totally real parabolic $CR$ algebra $(\frak{g},\frak{q}')$ are
the {\it real flag manifolds}, and that
there is a precise formula to compute their fundamental
groups
(see e.g. [DKV], [Wi]).
Let $\frak{h'}$ be a maximally noncompact Cartan subalgebra of $\frak{g}$
contained in $\frak{g}'_+$, fix an S-fit (and S-adapted) 
$C\in\frak{C}(\Cal{R}(\Hat{\frak{g}},\Hat{\frak{h}}'),\Cal{Q}')$ 
and let $\{\alpha_1,\hdots,\alpha_m\}=
\Cal{B}(C)\cap\rootre$. Then $\pi_1(M')$ is described by generators
$t_{\alpha_1},\hdots,t_{\alpha_m}$ that satisfy the relations\,:
\form{t_{\alpha_i}
t_{\alpha_j}=t_{\alpha_j}
t_{\alpha_i}^{(-1)^{(\alpha_i|\alpha_j^{\lor})}}\quad\text{for}\quad
1\leq i,j\leq m\,,\qquad t_{\alpha_i}
=e\quad\text{if}\quad \alpha_i\notin\Phi'_C\,.}
By using Theorem {\TRMKF} we can now give a description of
$\pi_1(M)$\,:
\thm{Let $(\frak{g},\frak{q})$ be an effective parabolic $CR$ algebra
and $(\vartheta,\frak{h})$ an adapted Cartan pair for $(\frak{g},\frak{q})$
with $\frak{h}$ maximally noncompact in $\frak{g}_+$. 
Let $\bold{H}$ be the corresponding
Cartan subgroup of $\bold{G}$. \par
Let $(\frak{g},\frak{q}')$ be
the totally real parabolic $CR$ algebra of Theorem {\TEOMHA},
$(\vartheta,\frak{h}')$ an adapted Cartan pair for $(\frak{g},\frak{q}')$
with $\frak{h}'$ 
maximally noncompact in $\frak{g}$.
\par
Fix a Weyl chamber $C'\in\frak{C}(\Cal{R}',\Cal{Q}')$,
where $\Cal{R}'=\Cal{R}(\Hat{\frak{g}},\Hat{\frak{h}}')$ and
$\Cal{Q}'$ is the parabolic set of $\frak{q}'$ in $\Cal{R}'$, 
that is S-fit 
to $(\frak{g},\frak{q}')$ and S-adapted.
With $\Cal{B}(C')=\{\alpha_1,\hdots,\alpha_{\ell}\}$ and
$\Cal{B}^*(C')=\{\omega_1,\hdots,\omega_\ell\}$ (defined by
$(\omega_i|\alpha_j^\vee)=\delta_{i,j}$), if
$\{\omega_1,\hdots,\omega_{\mu}\}=\Cal{B}^*(C')\cap
\left[{\Cal{Q}'}^{\,r}\right]^\perp$,
we have $\{\alpha_1,\hdots,\alpha_{\mu}\}=\Phi_{C'}(\Cal{Q}')
\cap{\Cal{R}'_{\roman{
re}}}$\,.
Let us consider the maps\par\smallskip\centerline{
$\boldsymbol{\varphi}^{\flat}\colon\bold{G}'_+@>>>\Bbb{Z}_2^{\mu}$ of $(\FRMLQ)$
\qquad and  \qquad
$\boldsymbol{\delta}:\pi_1(M')@>>>\Bbb{Z}_2^{\mu}$ of $(\FRMLKK)$.}\par\smallskip
Then we have\,:
\form{\qquad\qquad\pi_1(M)\,=\, \boldsymbol{\delta}^{-1}\left(\boldsymbol{\phi}^{\flat}
(\bold{H})\right)\,.\qquad\qquad\qed}}
\edef\TRMHC{\number\q.\number\x}
\cor{Let $(\frak{g},\frak{q})$, be an effective parabolic $CR$ algebra,
$M=M(\frak{g},\frak{q})$ the corresponding parabolic $CR$ manifold.
Let
$M'=M(\frak{g},\frak{q}')$ 
and $f:M@>>>M'$ be defined as in Theorem {\TEOMHA}.
If there is a Cartan subalgebra $\frak{h}$ adapted to both
$(\frak{g},\frak{q})$ and $(\frak{g},\frak{q}')$ and maximally
noncompact in $\frak{g}_+'$, such that moreover
$\bold{H}/\bold{Z}(\bold{G})$ is connected, then
$f_*:\pi_1(M)@>>>\pi_1(M')$ is an isomorphism.}
\edef\CRLHB{\number\q.\number\x}
\dimo For the isotropy subgroup $\bold{G}'_+$ of $M'$ we have, by
$(\FormulaKD)$, that
$\pi_0(\bold{G}'_+)\simeq\pi_0
\left(\bold{H}/\left(\bold{H}\cap\bold{S}'_0
\right)\right)$, where
$\bold{S}'_0$ is an analytic semisimple subgroup of $\bold{G}'_+$.
Since $\bold{Z}(\bold{G})\subset\bold{G}_+$, the inclusion
$\bold{H}\hookrightarrow \bold{G}'_+$ defines, passing to the quotients,
a surjective map $\pi_0(\bold{H}/\bold{Z}(\bold{G}))@>>>\pi_0
\left(\bold{G}'_+/\bold{G}_+\right)$. 
Thus the fiber of $f\colon M@>>>M'$ is connected and,
by Theorem {\TEOMHA}, $f_*:\pi_1(M)@>>>\pi_1(M')$ is an isomorphism.
\enddimo
\cor{With the notation of Corollary {\CRLHB},
if $\frak{g}$ is a simple Lie algebra of real type
$\roman{A I\! I\! Ib}$, $\roman{D I\! I\! Ia}$, then
the map $f_*\colon\pi_1(M)@>>>\pi_1(M')$ is an isomorphism.}
\dimo In fact in these cases
$\bold{H}/\bold{Z}(\bold{G})$ is connected 
for every choice of $\frak{h}$, and we can apply
Corollary {\CRLHB}
to obtain an isomorphism of the fundamental group
of $M(\frak{g},\frak{q})$ with the fundamental group of
a totally real $M(\frak{g},\frak{q}')$.\enddimo


\se{Examples}
In this section we apply our results to several examples.

\esemp{Consider the 
semisimple real Lie algebra
$\frak{g}=\frak{sl}(3,\Bbb{R})$ of real $3\times 3$ matrices with zero
trace. Let $e_1,e_2,e_3$ be the canonical basis of $\Bbb{R}^3\subset\Bbb{C}^3$ 
and $(\epsilon_1,\epsilon_2,\epsilon_3)$ the basis of $\Bbb{C}^3$ given by\,:
\par\smallskip\centerline{$\epsilon_1=e_1+ie_2\,,\qquad
\epsilon_2=e_1-ie_2\,,\qquad\epsilon_3=e_3\,.$}\par\smallskip\noindent
Let $\frak{q}$ be the complex Borel subalgebra of
complex matrices $Z\in\frak{sl}(3,\Bbb{C})$ such that
$Z(\langle\epsilon_1\rangle)\subset\langle \epsilon_1\rangle$ and
$Z(\langle \epsilon_1,\epsilon_2\rangle)\subset
\langle \epsilon_1,\epsilon_2\rangle$. Let $\frak{h}\subset\frak{q}$ 
be the Cartan  subalgebra of traceless matrices that are diagonal 
in the basis $(\epsilon_1,\epsilon_2,\epsilon_3)$. The corresponding 
Cartan subgroup\,:
\par\smallskip\centerline{$\bold{H}=\{\roman{diag}(\lambda,\bar\lambda,
|\lambda|^{-2})\in\bold{SL}(3,\Bbb{R})
\,|\,\lambda\neq 0\}$}\par\smallskip\noindent
is connected, hence also $\bold{G}_+$ is connected.
There exists
a unique Weyl chamber $C\in\frak{C}(\Cal{R}(\hat{\frak{g}},\hat{\frak{h}}))$ 
adapted to $(\frak{g},\frak{q})$, which is both
S-fit and V-fit.
The corresponding diagram is\,:
$$\xymatrix @M=0pt @R=2pt @!C=15pt{
\circledast\ar@{-}[r]
&\oplus\\
\alpha_1&\alpha_2\\
\times&\times}$$
and we see that  $M=M(\frak{g},\frak{q})$ is weakly degenerate.
The basis of the weakly nondegenerate reduction is the totally real
parabolic $CR$ manifold  $M'=M(\frak{g},\frak{q}')$,
where $\frak{q}'=\{Z\in\frak{sl}(3,\Bbb{C})\,|\,
Z(\langle \epsilon_1,\epsilon_2\rangle)\subset
\langle \epsilon_1,\epsilon_2\rangle\}$. Its diagrams, with respect
to the Cartan subalgebras $\frak{h}$ and $\frak{h}'$, 
where $\frak{h}'$ is the maximally noncompact Cartan subalgebra 
of traceless diagonal matrices in the basis $(e_1,e_2,e_3)$\,,
are\,:
$$\xymatrix @M=0pt @R=2pt @!C=15pt{
\circledast\ar@{-}[r]
&\oplus\\
\alpha_1&\alpha_2\\
&\times}
\quad\simeq\quad
\xymatrix @M=0pt @R=2pt @!C=15pt{
\bbianca\ar@{-}[r]
&\bbianca\\
\alpha_1&\alpha_2\\
&\times}$$
By [Wi] the fundamental group of $M'$ is $\pi_1(M')=\Bbb{Z}_2$.
On the other hand, by Corollary {\CORGH} the fiber of the weakly
nondegenerate reduction has two connected components.
Hence the exact sequence $(\FRMLHA)$ implies that $M$ is simply connected.}

\esemp{Let us compute the fundamental group of $M(\frak{g},\frak{q})$ in
the case of Example {\ESMPKZ}. 
 Since in this case
$\frak{q}$ is Borel and the Cartan subalgebra
$\frak{h}$ is maximally compact, we know that the isotropy
$\bold{G}_+=\{X\in\frak{g}\,|\roman{ad}_{\Hat{\frak{g}}}(X)
(\frak{q})=\frak{q}\}$ is connected (see [Kn, Prop.7.90]).
Consider the fiber $F$ over the point
$(\langle e_1,e_4\rangle,\langle e_1,e_2,e_4,e_5\rangle)$. 
We can verify that $F$ has $4$ connected components
and $\pi_0(F)\simeq \Bbb{Z}_2\times\Bbb{Z}_2$. The fundamental
group of $M(\frak{g},\frak{q}_3)$ can be computed using [Wi].
We have $\pi_1(M(\frak{g},\frak{q}_3))\simeq \Bbb{Z}_2\times
\Bbb{Z}_2$ and thus, from the exact sequence\,:
$$\bold{1}@>>>\pi_1(M(\frak{g},\frak{q}))@>>>\pi_1(M(\frak{g},\frak{q}_3))
\simeq\Bbb{Z}_2\times\Bbb{Z}_2 @>>>\pi_0(F)\simeq \Bbb{Z}_2\times\Bbb{Z}_2
@>>>\bold{1}$$
we obtain that $M(\frak{g},\frak{q})$ is simply connected.}

\esemp{\edef\ESEHA{\number\q.\number\y}
Consider 
the complex flag manifold $\frak{M}$
of $\bold{SO}(5,\Bbb{C})$, consisting of the complex 
projective lines contained in the quadric $\{z_2^2+2z_0z_4+2z_1z_3=0\}
\subset\Bbb{CP}^4$. We identify the Lie algebra $\frak{so}(5,\Bbb{C})$
to the matrix algebra\,:
$$\frak{so}(5,\Bbb{C})\simeq\Hat{\frak{g}}\,=\left\{
Z\in\frak{sl}(5,\Bbb{C})\,|\, {^t\! Z}S_5+S_5Z=0\,\right\}\quad
\text{for}\quad S_5=\left(\smallmatrix
&&&&1\\\vspace{-1\jot}
&&&1\\\vspace{-1\jot}
&&1\\\vspace{-1\jot}
&1\\\vspace{-1\jot}
1\endsmallmatrix\right)\,.$$
We consider the real form $\frak{g}\simeq\frak{so}(2,3)$ of
$\frak{so}(5,\Bbb{C})$ defined by\,:
$$\frak{so}(2,3)\simeq{\frak{g}}\,=\left\{
Z\in\Hat{\frak{g}}\,|\, Z^*K+KZ=0 \right\}\quad
\text{for}\quad K=\left(\smallmatrix
&&&&1\\
&1&\;\,0&\;0\\
&0&-1\,&\;0\\
&0&\;\,0&\;1\\
1\endsmallmatrix\right)\,.$$
Let $M$ be the orbit of the projective line
corresponding to the plane $\ell_2=\langle e_1,e_2\rangle$ of $\Bbb{C}^5$
by the action of the analytic subgroup $\bold{G}$
with Lie algebra $\frak{g}$\,:
$M$ is the submanifold 
of the Grassmannian of the complex $2$-planes
in $\Bbb{C}^5$, consisting of those 
that are totally isotropic for the symmetric form $S_5$ and
degenerate, with signature $(+,0)$, with respect to the Hermitian symmetric
form $K$. Denoting by $\frak{q}$ the stabilizer of $\langle e_1,e_2\rangle$
in $\bold{SO}(5,\Bbb{C})$, we have
$M=M(\frak{g},\frak{q})$. Take the Cartan subalgebra 
$\frak{h}$ of $\frak{g}$
consisting of the diagonal matrices. With $e_1,e_2$ also denoting the
value of the first and the second diagonal entry, we note that 
the conjugation $\sigma$ is defined in $\frak{h}_{\Bbb{R}}^*$ by
$\sigma(e_i)=-(-1)^ie_i$. Take the Weyl chambers 
$C,C'\in\frak{C}(\Cal{R},\Cal{Q})$ associated
to the basis $\Cal{B}(C)=\{e_1-e_2,e_2\}$ and $\Cal{B}(C')=\{e_2-e_1,e_1\}$.
Then $C$ is S-fit and $C'$ V-fit for $(\frak{g},\frak{q})$.
We can describe $M$ by the cross-marked diagrams\,: 
\medskip
$$
\xymatrix @M=0pt @R=2pt @!C=5pt{
\oplus\ar@{=}[rr]&>&\circledast\\
\alpha_1&&\alpha_2\\
&&\times}
\qquad\qquad 
\xymatrix @M=0pt @R=2pt @!C=5pt{
\ominus\ar@{=}[rr]&>&\bbianca\\
\alpha_1'&&\alpha_2'\\
&&\times}
$$
\medskip
From the first diagram we see that $(\frak{g},\frak{q})$ is fundamental,
since $\Phi_C=\{\alpha_2\}$ and 
$\bar\alpha_1=\alpha_1+2\alpha_2\succ_C\alpha_2$; from the second
we see that $(\frak{g},\frak{q})$ is weakly non-degenerate, because
$\Phi_{C'}=\{\alpha'_2\}$ and
$\bar\alpha_2'=\alpha_2'\succ_{C'}0$.
The weakest $CR$ model of $(\frak{g},\frak{q})$ is the
parabolic $CR$ algebra $(\frak{g},\frak{v})=(\frak{g},
\frak{q}_{\Phi_C^\sharp})$, with $\Phi_C^\sharp=\{\alpha_1,\alpha_2\}$.
The weakly non-degenerate reduction of  $(\frak{g},\frak{v})$ is
the totally real parabolic $CR$ algebra $(\frak{g},\frak{q}_{\Psi_C})$,
with $\Psi_C=\{\alpha_1\}$\,:
\medskip
$$
\xymatrix @M=0pt @R=2pt @!C=5pt{
\oplus\ar@{=}[rr]&>&\circledast\\
\alpha_1&&\alpha_2\\
\times &&\times} \quad{\longrightarrow}\quad
\xymatrix @M=0pt @R=2pt @!C=5pt{
\oplus\ar@{=}[rr]&>&\circledast\\
\alpha_1&&\alpha_2\\
\times } 
$$
\medskip
By composition we obtain the  
$\bold{G}$-equivariant projection\,:\par\centerline{
$M(\frak{g},\frak{q})@>{\sim}>>
M(\frak{g},\frak{v})@>>>M'=M(\frak{g},\frak{q}_{\Psi_C=\{\alpha_1\}})$}
\par\noindent of $M$ 
onto a totally real parabolic $CR$ manifold $M'$. This
projection associates to each $\ell_2\in M$ the isotropic line
$\ell_2\cap\ell_2^\perp$, where $\perp$ is taken with respect to
the Hermitian symmetric form $K$. The fiber over $\ell_1=\langle e_1\rangle$
consists of the planes generated by the columns of the matrices
$$\left(\matrix
1&0\\\vspace{-1\jot}
0&z_0\\\vspace{-1\jot}
0&z_1\\\vspace{-1\jot}
0&z_2\\\vspace{-1\jot}
0&0
\endmatrix\right)\qquad\text{with}\qquad\cases
(z_0:z_1:z_2)\in\Bbb{CP}^2\\
2z_0z_2+z_1^2=0\\
z_0\bar{z}_0+z_2\bar{z}_2>z_1\bar{z}_1\,.
\endcases
$$
Then we see that the fiber is biholomorphic to 
$\Bbb{CP}^1\setminus\Bbb{RP}^1$, that is the disjoint union of two disks
in $\Bbb{C}$. Thus the fiber $F$ of the projection $M@>>>M'$ has two connected
components and $\pi_0(F)\simeq\Bbb{Z}_2$.
Note that, by [Wi], $\pi_1(M')\simeq\Bbb{Z}$. Thus the
Serre's exact sequence\,:
$$\CD
\bold{1}@>>>\pi_1(M)@>>>\pi(M')\simeq\Bbb{Z}@>>>\pi_0(F)\simeq\Bbb{Z}_2
@>>>\bold{1}\endCD$$
shows that $\pi_1(M)\simeq 2\Bbb{Z}\simeq\Bbb{Z}$\,.}\par
\esemp{Let $(\epsilon_i)_{1\leq i\leq 4}$ be the canonical basis of
$\Bbb{C}^4$, and $\frak{g}\simeq\frak{sl}(4,\Bbb{R})$ consist
of the elements of $\frak{sl}(4,\Bbb{C})$ that have real entries
in the basis $(\epsilon_i)_{1\leq i\leq 4}$. We introduce the basis
$$e_1=\epsilon_1+i\epsilon_2\,,\;
e_2=\epsilon_1-i\epsilon_2\,,\;
e_3=\epsilon_3+i\epsilon_4\,,\;
e_4=\epsilon_3-i\epsilon_4\,.
$$
We take the complex flag manifold $\frak{M}$ whose points are
the pairs $(\ell_1,\ell_3)$ of a complex line $\ell_1$ of $\Bbb{C}^4$
and a complex $3$-plane $\ell_3$ 
with $\ell_1\subset\ell_3\subset\Bbb{C}^4$, and consider the 
$\bold{G}$-orbit $M$ that contains the point 
$\frak{o}=\langle e_1\rangle\subset\langle e_1,e_2,e_3\rangle$\,:
we have $M=M(\frak{g},\frak{q})$ where $\frak{q}$ is the stabilizer
in $\Hat{\frak{g}}\simeq\frak{sl}(4,\Bbb{C})$ of $\frak{o}$. 
Consider the Cartan subalgebra $\frak{h}$ of 
the elements of $\frak{g}$ that are diagonal matrices in the basis
$(e_i)_{1\leq i\leq 4}$. With $e_i(H)$ also denoting the value of the
$i$-th entry of $H\in\frak{h}_{\Bbb{R}}$, we note that
$\Cal{B}(C)=\{\alpha_i=e_i-e_{i+1}\,|\, 1\leq i\leq 3\}$
is the system of simple roots for an S-fit Weyl chamber
$C\in\frak{C}(\Cal{R},\Cal{Q})$, and
$\Cal{B}(C)=\{\alpha'_1=e_1-e_3\,,\;\alpha'_2=e_3-e_2\,,\;\alpha'_3=
e_2-e_4\}$ 
is the system of simple roots for a V-fit Weyl chamber. The corresponding
cross-marked diagrams are\,:
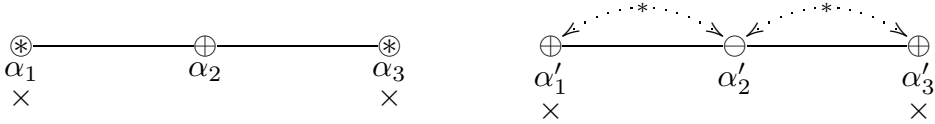
$$
\xymatrix @M=0pt @R=2pt @!C=5pt{
\circledast\ar@{-}[rr]
&&\oplus\ar@{-}[rr]
&&\circledast\\
\alpha_1&&\alpha_2&&\alpha_{3}\\
\times&&&&\times} \qquad\qquad
\xymatrix @M=0pt @R=2pt @!C=5pt{
\oplus\ar@{-}[rr]\ar@/^1pc/@{<.>}[rr]|{\ast}
&&\ominus\ar@{-}[rr]\ar@/^1pc/@{<.>}[rr]|{\ast}
&&\oplus\\
\alpha'_1&&\alpha'_2&&\alpha'_{3}\\
\times&&&&\times} 
$$ 
Since $\bar\alpha_2=\alpha_1+\alpha_2+\alpha_3$, from the first
we see that $(\frak{g},\frak{q})$ is fundamental, while the second
shows that $(\frak{g},\frak{q})$ is weakly nondegenerate. 
Its weakest $CR$ model $(\frak{g},\frak{v})$ has a totally real
weakly nondegenerate reduction
$(\frak{g},\frak{q}')$, corresponding to diagrams\,:
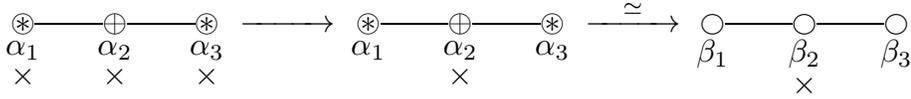
$$\CD
\xymatrix @M=0pt @R=2pt @!C=5pt{
\circledast\ar@{-}[r]
&\oplus\ar@{-}[r]
&\circledast\\
\alpha_1&\alpha_2&\alpha_{3}\\
\times&\times&\times}
@>>>
\xymatrix @M=0pt @R=2pt @!C=5pt{
\circledast\ar@{-}[r]
&\oplus\ar@{-}[r]
&\circledast\\
\alpha_1&\alpha_2&\alpha_{3}\\
&\times}
@>{\simeq}>>
\xymatrix @M=0pt @R=2pt @!C=5pt{
\bbianca\ar@{-}[r]
&\bbianca\ar@{-}[r]
&\bbianca\\
\beta_1&\beta_2&\beta_{3}\\
&\times}
\endCD
$$ 
where the last diagram is obtained by utilizing the Cartan subalgebra
of real diagonal matrices of $\frak{g}$ with respect to the canonical
basis $(\epsilon_i)_{1\leq i\leq 4}$. Using [Wi], we obtain
that $\pi_1(M(\frak{g},\frak{q}'))\simeq\Bbb{Z}_2$. 
The isotropy subgroup $\bold{G}_+$ of $M(\frak{g},\frak{q}')$ is
isomorphic to the group of matrices of the form $\pmatrix 
A&B\\
0&C\endpmatrix$ with $A,B,C$ real $2\times 2$ matrices with
$\roman{det}(A)\cdot\roman{det}(C)=1$, and hence has two connected 
components. The isotropy subgroup $\bold{G}_+$ of
$M(\frak{g},\frak{q})$ is connected\,: indeed $\bold{G}_+=\bold{N}\ltimes
\bold{H}$ for an Euclidean $\bold{N}=\exp(\frak{n})$ and a 
Cartan subgroup $\bold{H}=\bold{Z}_{\bold{G}}(\frak{h})$ that is
connected because $\frak{h}$ is maximally noncompact (cf. [Kn, Proposition
7.90, p.488]. Thus the fiber $\bold{G}'_+/\bold{G}_+$
has two connected components.
Thus, from the exact sequence\,:
$$\CD
\bold{1}@>>>\pi_1(M(\frak{g},\frak{q}))@>>>
\underbrace{\pi_1(M(\frak{g},\frak{q}'))}
@>>>\underbrace{\pi_0(\bold{G}'_+/\bold{G}_+)}
@>>>\bold{1}\\\vspace{-1\jot}
@.@.\simeq\Bbb{Z}_2 @.\simeq\Bbb{Z}_2 \endCD$$
we obtain that $M(\frak{g},\frak{q})$ is simply connected.}

\esemp{
Let $(\epsilon_i)_{1\leq i\leq 6}$ be the canonical basis of $\Bbb{C}^6$.
Let $\bold{G}\simeq \bold{SL}(6,\Bbb{R})$, with Lie algebra $\frak{g}$,
be the subgroup of $\bold{SL}(6,\Bbb{C})$ consisting of
the matrices with real entries. Consider the basis
$$e'_1=\epsilon_1+i\epsilon_4\,,\; e'_2=\epsilon_2\,,\;
e'_3=\epsilon_3+i\epsilon_6\,,\;
e'_4=\epsilon_1-i\epsilon_4\,,\;
e'_5=\epsilon_5\,,\;
e'_6=\epsilon_3+i\epsilon_6\,.$$
Let $\frak{q}$ be the stabilizer of 
$\langle e_1',e_2'\rangle\subset\langle e_1',e_2',e_3',e_4'\rangle$
in $\frak{sl}(6,\Bbb{C})$ and consider the parabolic $CR$ algebra
$(\frak{g},\frak{q})$. 
The matrices in 
$\frak{g}$ that are diagonal in the basis $e'_1,\hdots,e'_6$
form a Cartan subalgebra $\frak{h}$ 
of $\frak{g}$ adapted to $(\frak{g},\frak{q})$. Then $\frak{h}_{\Bbb{R}}$
consists of the matrices of 
$\Hat{\frak{g}}=\frak{sl}(6,\Bbb{C})$ that are real and diagonal
in the basis $e'_1,\hdots,e'_6$. We identify $e'_i$ to the evaluation
function of the $i$-th diagonal term of $H\in\frak{h}_{\Bbb{R}}$.
Then the $\alpha_i'=e'_i-e'_{i+1}$ are the simple root of a
$C'\in\frak{C}(\Cal{R},\Cal{Q})$ that is V-fit for $(\frak{g},\frak{q})$.
We have the cross-marked diagram for $(\frak{g},\frak{q})$\,:
$$
\xymatrix @M=0pt @R=2pt @!C=5pt{
\ominus\ar@{-}[r]\ar@/^1pc/@{<.>}[rr]|{\ast}
&\oplus\ar@{-}[r]&\ominus\ar@{-}[r]
\ar@/^1pc/@{<.>}[rr]|{\ast}&\oplus
\ar@{-}[r]&\ominus\\
\alpha_1'&\alpha_2'&\alpha_{3}'&\alpha_4'&\alpha_5'\\
&\times&&\times&}
$$
Since $\Phi_{C'}=\{\alpha_2',\alpha_4'\}$ and 
$\bar\alpha'_2=\alpha'_2+\alpha'_3+\alpha'_4+\alpha'_5\succ_{C'}0$,
$\bar\alpha'_4=
\alpha'_1+\alpha'_2+\alpha'_3+\alpha'_4\succ_{C'}0$, the $CR$ algebra
$(\frak{g},\frak{q})$ is weakly nondegenerate.
We obtain an S-fit Weyl chamber for $(\frak{g},\frak{q})$
by reordering the basis $e'_1,\hdots,e'_6$.
Set\,:
$$e_1=e'_2\,,\; e_2=e_1'\,,\;
e_3=e'_4\,,\; 
e_4=e'_3\,,\; 
e_5=e'_6\,,\; 
e_6=e'_5\,.
$$ 
Then $\alpha_i=e_i-e_{i+1}$ ($1\leq i\leq 5$) is the basis $\Cal{B}(C)$
of an S-fit
 Weyl chamber $C\in\frak{C}(\Cal{R},\Cal{Q})$ yielding the
cross-marked diagram\,:
$$
\xymatrix @M=0pt @R=2pt @!C=5pt{
\oplus\ar@{-}[r]
&\grigia\ar@{-}[r]
&\oplus\ar@{-}[r]
&\grigia\ar@{-}[r]
& \oplus\\
\alpha_1&\alpha_2&\alpha_{3}&\alpha_4&\alpha_5\\
&\times&&\times&}
$$
Since $\Phi_C=\{\alpha_2,\alpha_4\}$ and $\bar\alpha_3=\alpha_2+\alpha_3+
\alpha_4\succ_C\alpha_2,\alpha_4$, the $CR$ algebra $(\frak{g},\frak{q})$
is also fundamental. The weakest $CR$ model $(\frak{g},\frak{v})$ of
 $(\frak{g},\frak{q})$ is obtained by taking the complex Borel
subalgebra of $\Hat{\frak{g}}$ associated to the chamber $C$. We have
$\frak{v}=\frak{q}_{\Phi_C^\sharp}$ with 
$\Phi_C^\sharp=\{\alpha_1,\alpha_2,\alpha_3,\alpha_4,\alpha_5\}$.
The weakly nondegenerate reduction of $(\frak{g},\frak{q})$ is
the totally real $CR$ algebra $(\frak{g},\frak{q}_{\Psi_C})$ with
$\Psi_C=\{\alpha_1,\alpha_3,\alpha_5\}$\,:
$$
\xymatrix @M=0pt @R=2pt @!C=5pt{
\oplus\ar@{-}[r]
&\grigia\ar@{-}[r]
&\oplus\ar@{-}[r]
&\grigia\ar@{-}[r]
& \oplus\\
\alpha_1&\alpha_2&\alpha_{3}&\alpha_4&\alpha_5\\
\times&\times&\times&\times&\times}
\quad{\longrightarrow}\quad
\xymatrix @M=0pt @R=2pt @!C=5pt{
\oplus\ar@{-}[r]
&\grigia\ar@{-}[r]
&\oplus\ar@{-}[r]
&\grigia\ar@{-}[r]
& \oplus\\
\alpha_1&\alpha_2&\alpha_{3}&\alpha_4&\alpha_5\\
\times&&\times&&\times}
$$
By composition we obtain a $\bold{G}$-equivariant fibration\,:\par
\centerline{
$M=M(\frak{g},\frak{q})@>>>M'=M(\frak{g},\frak{q}_{\Psi_C})$,
with $(\ell_2\subset\ell_4)
@>>>(\ell_2\cap\bar{\ell}_2\subset\ell_4\cap\bar{\ell}_4\subset
\ell_4+\bar{\ell}_4)$.}\par
The fiber $F$ over $\langle e_1\rangle\subset\langle e_1,e_2,e_3\rangle\subset
\langle e_1,e_2,e_3,e_4,e_5\rangle$ consists of the 
pairs $(\ell_2\subset\ell_4)$ where
$\ell_2$ is a complex $2$-plane with 
$\langle e_1\rangle\subset\ell_2\subset\langle e_1,e_2,e_3\rangle$
and $\ell_4$ is a complex $4$-plane with
$\langle e_1,e_2,e_3\rangle\subset\ell_4\subset\langle e_1,e_2,e_3,e_4,e_5
\rangle$, with $\ell_2\neq\bar{\ell}_2$ and $\ell_4\neq\bar{\ell}_4$.
Thus $\pi_0(F)\simeq \Bbb{Z}_2\times\Bbb{Z}_2$. On the other hand, by
[Wi], we have 
$\pi_1(M')\simeq\Bbb{Z}_2\times\Bbb{Z}_2\times
\Bbb{Z}_2$. From the exact sequence
$$\CD \bold{1}@>>>\pi_1(M)@>>>\pi_1(M')\simeq\Bbb{Z}_2^3@>>>
\pi_0(F)\simeq\Bbb{Z}_2^2@>>>\bold{1}\endCD$$
we obtain that $\pi_1(M)\simeq\Bbb{Z}_2$.}
\par
\esemp{Let $\Hat{\frak{g}}\simeq\frak{so}(5,\Bbb{C})$ be as in Example {\ESEHA}.
We take now $\frak{g}\simeq\frak{so}(2,3)$, defined by\,:
$$\frak{g}=\left\{Z\in\Hat{\frak{g}}\,|\, Z^*K+KZ=0\right\}\qquad
\text{with}\quad K=\left(\smallmatrix
&&\left[\smallmatrix
1\\\vspace{-1\jot}
&1\endsmallmatrix\right]\\
&1\\
\left[\smallmatrix
1\\\vspace{-1\jot}
&1\endsmallmatrix\right]\endsmallmatrix\right)\,.$$
Let $\frak{q}$ be the stabilizer of the line $\langle e_1\rangle
\subset\Bbb{C}^5$ and consider the orbit $M=M(\frak{g},\frak{q})$. 
Our $M$ is one of the two connected components of the
manifold $M^+$ 
the nonreal lines $\ell_1$, 
contained in the quadric cone
$\{2z_0z_4+2z_1z_3+z_2^2=0\}\subset\Bbb{C}^5$, for which
$\left(\ell_1+\bar{\ell}_1\right)$ is
a totally isotropic complex $2$-plane
for the Hermitian symmetric form $K$. 
The involution $\ell_1 @>>>\bar{\ell}_1$
interchanges the two connected components of $M^+$.
\par
The diagonal matrices in $\frak{g}$ define a Cartan subalgebra of
$\frak{g}$ adapted to $(\frak{g},\frak{q})$. As usual we denote again by
$e_i$ the value of the $i$-th entry in the diagonal of $\frak{h}_{\Bbb{R}}$.
Then we define a V-fit $C'$ and an S-fit $C$ by taking
$\Cal{B}(C')=\{\alpha'_1=e_1+e_2,\;\alpha'_2=-e_2\}$ and
$\Cal{B}(C)=\{\alpha_1=e_1-e_2,\;\alpha_2=e_2\}$. The 
associated cross-marked diagrams are\,:
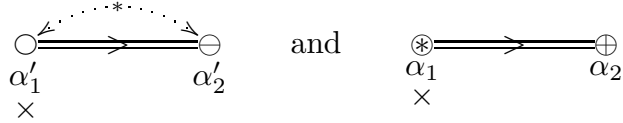
$$
\xymatrix @M=0pt @R=2pt @!C=5pt{
{\bbianca}\ar@{=}[rr]\ar@/^1pc/@{<.>}[rr]|{\ast} 
&>&\ominus\\
\alpha_1'&&\alpha_2'\\
\times}\qquad\text{and}\qquad
\xymatrix @M=0pt @R=2pt @!C=5pt{
{\circledast}\ar@{=}[rr]&>&\oplus\\
\alpha_1&&\alpha_2\\
\times}
$$
From the first, as $\frak{q}=\frak{q}_{\Phi_{C'}}$ with
$\Phi_{C'}=\{\alpha_1'\}$ and $\bar\alpha_1'=\alpha_1'\succ_{C'}0$,
we see that $(\frak{g},\frak{q})$ is weakly nondegenerate; from
the second, as 
$\frak{q}=\frak{q}_{\Phi_{C}}$ with
$\Phi_{C}=\{\alpha_1\}$ and 
$\bar\alpha_2=\alpha_1+\alpha_2\succ_C\alpha_1$,
we see that $(\frak{g},\frak{q})$ is also fundamental.
The weakest $CR$ model of $(\frak{g},\frak{q})$ is 
$(\frak{g},\frak{q}_{\Psi_C^\sharp})$ with 
$\Psi_C^\sharp=\{\alpha_1,\alpha_2\}$.
The basis of its weakly nondegenerate 
reduction is the totally real
$CR$ algebra $(\frak{g},\frak{q}_{\Psi_C})$ with $\Psi_C=\{\alpha_2\}$.
We can represent these maps by the diagram\,:
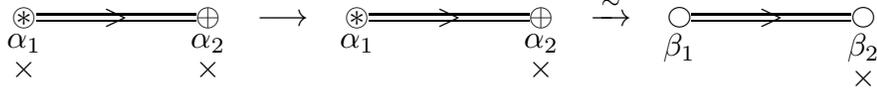
$$
\xymatrix @M=0pt @R=2pt @!C=5pt{
{\circledast}\ar@{=}[rr]&>&\oplus\\
\alpha_1&&\alpha_2\\
\times&&\times}\quad{\longrightarrow}\quad
\xymatrix @M=0pt @R=2pt @!C=5pt{
{\circledast}\ar@{=}[rr]&>&\oplus\\
\alpha_1&&\alpha_2\\
&&\times}\quad{@>\sim>>}\quad
\xymatrix @M=0pt @R=2pt @!C=5pt{
{\bbianca}\ar@{=}[rr]&>&\bbianca\\
\beta_1&&\beta_2\\
&&\times}
$$
\noindent where the last is the cross-marked Satake diagram of 
$(\frak{g},\frak{q}_{\Psi_C})$, i.e. the
cross-marked diagram for an S-fit and S-adapted $C$ and
a maximally noncompact $\frak{h}$. We have, by [Wi, Theorem 1.1],
$\pi_1(M(\frak{g},\frak{q}_{\Psi_C}))\simeq\Bbb{Z}_2$. 
We observe that the stabilizer in
the connected component of the identity
of $\bold{SO}(2,3)$ of a totally isotropic $2$-plane
$\ell_2\subset\Bbb{R}^5$ keeps its orientation and 
is connected. Thus the 
the fiber $F$ of the projection 
$M(\frak{g},\frak{q})@>>>M(\frak{g},\frak{q}_{\Psi_C})$
is connected because 
the isotropy subgroup of $M(\frak{g},\frak{q}_{\Psi_C})$ is connected.
Therefore 
$\pi_1(M(\frak{g},\frak{q}))\simeq\pi_1(M(\frak{g},\frak{q}_{\Psi_C}))
\simeq
\Bbb{Z}_2$.}\par
\esemp{Let ${\frak{g}}$ be a simple real Lie algebra of type
$\roman{F}\roman{I}$ (split real form). We fix a Cartan subalgebra 
$\frak{h}$ of $\frak{g}$, such that the conjugation $\sigma$ 
defined in $\Hat{\frak{g}}$ by the real form $\frak{g}$ restricts in 
$\frak{h}^*_{\Bbb{R}}=\Bbb{R}^4=\langle e_1,e_2,e_3,e_4\rangle_{\Bbb{R}}$
to the linear involution that is defined on the canonical basis by\,:
$\sigma(e_1)=-e_3\,,\;\sigma(e_2)=e_4\,,\;\sigma(e_3)=-e_1\,,\;
\sigma(e_4)=e_2\,.$
The vectors $\alpha_1=e_2-e_3$,
$\alpha_2=e_3-e_4$, $\alpha_3=e_4$ and
$\alpha_4=\frac{1}{2}\left(e_1-e_2-e_3-e_4\right)$ are a basis of simple roots
of $\Cal{R}=\Cal{R}(\Hat{\frak{g}},\Hat{\frak{h}})$.
We consider the parabolic $CR$ manifold $M(\frak{g},\frak{q})$ that
corresponds to the cross-marked diagram\,:
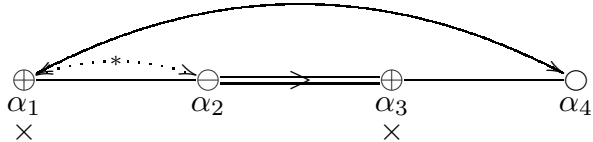
\bigskip
$$
\xymatrix @M=0pt @R=2pt @!C=5pt{
\oplus\ar@{-}[rr]\ar@/^.5pc/@{<.>}[rr]|{\ast} \ar@/^2pc/@{<->}[rrrrrr] 
&&\ominus\ar@{=}[rr]&>
&\oplus\ar@{-}[rr]
&&\bbianca\\
\alpha_1&&\alpha_2&&\alpha_3&&\alpha_4\\
\times&& &&\times}
$$
This is a representation of $(\frak{g},\frak{q})$ in a V-fit Weyl chamber
$C\in\frak{C}(\Cal{R},\Cal{Q})$.
Since $\bar\alpha_1\succ_C 0$, $\bar\alpha_3\succ_C 0$, 
by Theorem {\TEOEI} the $CR$ algebra $(\frak{g},\frak{q})$ is
weakly nondegenerate.\par
Consider the Weyl chamber $C_1$ obtained from $C$ by the symmetry
$s_{\alpha_2}$. We have
$\Cal{B}(C_1)=\left\{e_2-e_4,e_4-e_3,e_3,\frac{1}{2}(e_1-e_2-e_3-e_4)
\right\}=\{\alpha'_1,\alpha'_2,\alpha'_3,\alpha'_4\}$ 
and the cross-marked diagram
for $(\frak{g},\frak{q})$ for the chamber $C_1\in\frak{C}(\Cal{R},\Cal{Q})$
is\,:
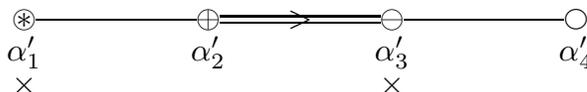
$$
\xymatrix @M=0pt @R=2pt @!C=5pt{
\circledast\ar@{-}[rr]
&&\oplus\ar@{=}[rr]&>
&\ominus\ar@{-}[rr]
&&\bbianca\\
\alpha'_1&&\alpha'_2&&\alpha'_3&&\alpha'_4\\
\times&& &&\times}
$$
This diagram is S-fit. Since $\bar\alpha'_2=2\alpha'_1+3\alpha'_2+4\alpha'_3
+2\alpha'_4$, by Theorem {\TeorEF} the parabolic $CR$ algebra
$(\frak{g},\frak{q})$ is also fundamental.
\par
Our next aim is to construct the fibration of Proposition {\ProposizioneLJ}.
First we observe that the weakest $CR$ model of $(\frak{g},\frak{q})$
is $(\frak{g},\frak{v}_1)$ with $\frak{v}_1=\frak{q}_{\{\alpha_1',\alpha_2',
\alpha_3'\}}$. Its weakly nondegenerate reduction is
$(\frak{g},\frak{q}_2)$ with $\frak{q}_2=\frak{q}_{\{\alpha_2'\}}$.
\par
To compute its weakest $CR$ model, we need to find an S-fit Weyl
chamber $C_2$ for $(\frak{g},\frak{q}_2)$. This is provided by
the basis of simple roots $\Cal{B}(C)=\{\beta_1,\beta_2,\beta_3,\beta_4\}$
with
$\beta_1=e_2-e_4$, $\beta_2=e_1-e_2$, $\beta_3=\frac{1}{2}(e_1-e_2-e_3-e_4)$,
$\beta_4=-\frac{1}{2}(e_1-e_2+e_3-e_4)$, that is obtained from
$\Cal{B}(C_1)$ by the rotation  
$s_{\frac{1}{2}(e_1-e_2+e_3-e_4)}\circ s_{e_3}$ of $\alweyl$.
The corresponding cross-marked diagram for $(\frak{g},\frak{q}_2)$ is\,:
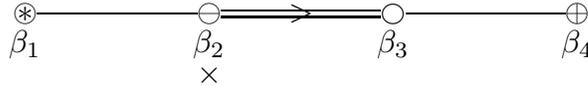
$$
\xymatrix @M=0pt @R=2pt @!C=5pt{
\circledast\ar@{-}[rr]
&&\ominus\ar@{=}[rr]&>
&\bbianca\ar@{-}[rr]
&&\oplus\\
\beta_1&&\beta_2&&\beta_3&&\beta_4\\
&&\times &&}
$$
Since $\bar\beta_4=\beta_1+2\beta_2+2\beta_3+\beta_4$, the weakest $CR$ model
of $(\frak{g},\frak{q}_2)$ is $(\frak{g},\frak{v}_2)$ with
$\frak{v}_2=\frak{q}_{\{\beta_2,\beta_4\}}$. The $CR$ algebra
$(\frak{g},\frak{v}_2)$ has the weakly nondegenerate reduction
$(\frak{g},\frak{q}_{\{\beta_4\}})$. 
The element $A=(-1,-1,+1,-1)$ of $\frak{h}_{\Bbb{R}}$ 
defines the parabolic set of $\frak{v}_2$ and $\bar{A}=A$ shows then
that $(\frak{g},\frak{v}_2)$ is totally real. Thus, by choosing a
maximally noncompact Cartan subalgebra $\frak{h}'$ adapted
to $(\frak{g},\frak{v}_2)$, we can associate to $(\frak{g},\frak{v}_2)$
its cross-marked Satake diagram as a totally real parabolic minimal
$CR$ algebra\,:
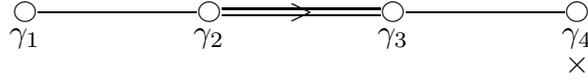
$$
\xymatrix @M=0pt @R=2pt @!C=5pt{
\bbianca\ar@{-}[rr]
&&\bbianca\ar@{=}[rr]&>
&\bbianca\ar@{-}[rr]
&&\bbianca\\
\gamma_1&&\gamma_2&&\gamma_3&&\gamma_4\\
&& && &&\times}
$$
Then the isotropy subgroup $\bold{G}_+'$ of $M(\frak{g},\frak{v}_2)$
has two connected components, and the fundamental group
$\pi_1(M(\frak{g},\frak{v}_2))$ is isomorphic to $\Bbb{Z}_2$ and is
generated by any simple path joining the two connected components of
$\bold{G}'_+$. By using Lemma {\LEMAIA} and 
$(\FormulaKD)$ of Theorem {\TRMAZ}, we find that $\bold{H}$ has four
connected components and $\boldsymbol{\varphi}^\flat(\bold{H})$ 
has two connected components. Hence Theorem {\TRMHC}
yields $\pi_1(M)\simeq\Bbb{Z}_2$.}





\Refs
\widestnumber\key{FHW}

\ref\key AF
\by A.Andreotti, G.Fredricks
\paper Embeddability of real analytic Cauchy-Riemann manifolds
\jour Ann. Scuola Norm. Sup. Pisa (4)
\vol 6
\yr 1979
\pages 285-304
\endref

\ref\key AFN
\by A.Andreotti, G.Fredricks, M.Nacinovich
\paper On the absence of Poincar\'e lemma in Tangential Cauchy-Riemann
complexes
\jour Ann. Scuola Norm. Sup. Pisa
\yr 1981
\vol 8
\pages 365-404
\endref

\ref\key AMN
\by A.Altomani, C.Medori, M.Nacinovich
\paper The $CR$ structure of minimal orbits in complex flag manifolds
\jour J. Lie Theory
\vol 16
\yr 2006
\pages 483-530
\endref

\ref\key Ar
\by S.Araki
\paper On root systems and an infinitesimal classification
of irreducible symmetric spaces
\jour J. Math. Osaka City Univ.
\vol A 13
\yr 1962
\pages 1-34
\endref

\ref\key BER
\by M.S.Baouendi, P.Ebenfelt, L.Preiss Rothschild
\book Real submanifolds in complex space and their mappings
\publ Princeton University Press
\yr 1999
\publaddr Princeton, NJ
\endref

\ref\key BlG
\by T.Bloom, I.Graham
\paper A geometric characterization of points of type $m$ on real
   submanifolds of ${\Bbb{C}}^{n}$,
\jour J. Differential Geometry
\vol 12
\yr 1977
\pages 171-182
\endref

\ref\key BL
\by R.Bremigan, J.Lorch
\paper Orbit duality for flag manifolds
\jour Manuscripta Math.
\vol 109
\yr 2002
\pages 233-261
\endref
  


\ref\key BT1
\by A.Borel, J.Tits
\paper Groupes r\'eductifs
\jour Publ. Math. IHES
\vol 27
\yr 1965
\pages 55-151
\endref

\ref\key BT2
\bysame
\paper Compl\'ements \`a l'article {``Groupes r\'eductifs''}
\jour Publ. Math. IHES
\vol 41
\yr 1972
\pages 253-276
\endref

\ref\key B1
\by N.Bourbaki
\book Groupes et alg\`ebres de Lie, Ch. 4,5 et 6
\publ Masson
\yr 1981
\publaddr Paris
\endref

\ref\key B2
\bysame
\book Groupes et alg\`ebres de Lie, Ch. 7 et 8
\publ Masson
\yr 1990
\publaddr Paris
\endref

\ref\key CS
\by L.Casian, R.J.Stanton
\paper Schubert cells and representation theory
\jour Invent. Math.
\vol 137
\yr 1999
\pages 461-539
\endref

\ref\key Ch
\by C.Chevalley
\book Th\'eorie des groupes de Lie
\publ Hermann
\yr 1955
\publaddr Paris
\endref

\ref\key DKV
\by J.J.Duistermaat, J.A.C.Kolk, V.S.Varadarajan
\paper Functions, flows and oscillatory integrals on flag manifolds and
conjugacy classes in real semisimple Lie groups
\jour Compositio Math.
\yr 1983
\vol 49
\pages 309-398
\endref

\ref\key Fe
\by G.Fels
\paper Locally homogeneous finitely nondegenerate $CR$ manifolds
\moreref preprint, arXiv:\linebreak math.CV/0606032, Jun 2006
\pages 1-21
\endref

\ref\key FHW
\by G.Fels, A.T.Huckleberry, J.A.Wolf
\book Cycle spaces of flag domains. A complex geometric viewpoint
\bookinfo Progress in Mathematics
\vol 245
\publ Birkh\"auser 
\publaddr Boston
\yr 2006
\endref

\ref\key FeK
\by G.Fels, W.Kaup
\paper Classification of Levi degenerate homogeneous $CR$ manifolds
of dimension $5$
\moreref  arXiv.org:math/0610375 2006-10-11
\pages 1-49
\endref

\ref \key GM
\by S.Gindikin,
T.Matsuki
\paper Stein extensions of Riemannian symmetric spaces and dualities of
   orbits on flag manifolds
\jour Transform. Groups
\vol 8
\yr 2003
\pages 333-376
\endref

\ref\key HN
\by C.D.Hill, M.Nacinovich
\paper Pseudoconcave $CR$ manifolds
\inbook Complex analysis and geometry
\bookinfo Lect. Notes in Pure and Appl. Math.
\vol 173
\pages 275-297
\publ Marcel and Dekker
\publaddr New York
\yr 1996
\endref

\ref\key HW
\by A.T.Huckleberry, J.A.Wolf
\paper Schubert varieties and cycle spaces
\jour Duke Math. J.
\vol 120
\yr 2003
\pages 229-249
\endref

\ref\key Ka
\by M.Kashiwara
\paper $D$-modules and representation theory of Lie groups
\jour Ann. Inst. Fourier
\vol 43
\yr 1993
\pages 1597-1618
\endref

\ref\key KZ
\by W.Kaup, D.Zaitsev
\paper On the {CR}-structure of compact group orbits associated with
bounded symmetric domains
\jour Invent. Math.
\vol 153
\yr 2003
\pages 45-104
\endref

\ref\key Kn
\by A.Knapp
\book Lie groups beyond an introduction (2nd ed.)
\publaddr Boston
\yr 2002 
\publ Birkh\"auser
\endref

\ref\key LN
\by A.Lotta, M.Nacinovich
\paper On a class of symmetric $CR$ manifolds
\jour Adv. in Math.
\yr 2005
\vol 191
\pages 114-146
\endref


\ref\key MN1
\by C.Medori, M.Nacinovich 
\paper Levi-Tanaka algebras and homogeneous $CR$ manifolds
\jour Compositio Math.
\yr 1997
\vol 109
\pages 195-250
\endref

\ref\key MN2
\bysame
\paper Classification of semisimple Levi-Tanaka algebras
\jour Ann. Mat. Pura Appl.
\yr 1998
\vol 174
\pages 285-349
\endref

\ref\key MN3
\bysame
\paper Standard $CR$ manifolds of codimension $2$
\jour Transformation groups
\yr 2001
\vol 6
\pages 53-78
\endref

\ref\key MN4
\bysame
\paper Algebras of infinitesimal $CR$ automorphisms
\jour Journal of Algebra
\yr 2005
\vol 287
\pages 234-274
\endref


\ref\key Mo
\by G.D.Mostow
\paper On maximal subgroups of real Lie groups
\jour Ann. Math.
\yr 1961
\vol 74
\pages 503-517
\endref

\ref\key Po
\by M.M.Postnikov
\yr 1982
\book Lie groups and Lie algebras
\publ Nauka
\publaddr Moskow
\moreref French translation MIR 1985
\endref

\ref\key OV
\eds A.L. Onishchik, E.B. Vinberg
\book Lie groups and Lie algebras III
\publaddr Berlin
\yr 1994
\publ Springer
\bookinfo Encyclopaedia of Mathematical Sciences \vol 43
\endref

\ref\key S1
\by M.Sugiura
\jour J. Math. Soc. Japan
\vol 11
\yr 1959
\pages 374-434
\paper Conjugate classes of Cartan subalgebras in real semisimple Lie algebras
\endref

\ref\key S2
\bysame
\jour J. Math. Soc. Japan
\vol 23
\yr 1971
\pages 379-383
\paper Correction to my paper:
Conjugate classes of Cartan subalgebras in real semisimple Lie algebras
\endref

\ref\key Wa
\by G.Warner
\book Harmonic analysis on semi-simple Lie groups I
\yr 1972
\publ Springer
\publaddr Berlin
\endref

\ref\key Wi
\by M.Wiggerman
\paper The fundamental group of a real flag manifold
\jour Indag.Matem.,N.S.
\vol 9
\pages 141-153
\yr 1998
\endref

\ref\key W
\by J.A.Wolf
\paper The action of real semisimple Lie groups on a complex
manifold I: Orbit structure and holomorphic arc components
\jour Bull. Amer. Math. Soc.
\vol 75
\yr 1969
\pages 1121-1137
\endref


\endRefs
\enddocument